\documentclass[12pt]{article}
\usepackage{latexsym, amssymb, enumerate, amsmath}
\usepackage{graphicx} 
\begin{document}
\newcommand{\vare}{\varepsilon} 
\newtheorem{tth}{Theorem}[section]
\newtheorem{dfn}[tth]{Definition}
\newtheorem{lem}[tth]{Lemma}
\newtheorem{prop}[tth]{Proposition}
\renewcommand{\quad}{\hspace*{1em}}
\begin{center}
{\Large {\bf Bifurcations of 
Wavefronts on an  $r$-corner II:
Semi-local classification} }
\vspace*{0.4cm}\\
{\large Takaharu Tsukada}
\footnote{Higashijujo 3-1-16 
Kita-ku, Tokyo 114-0001
JAPAN. e-mail : tsukada@math.chs.nihon-u.ac.jp}
\vspace*{0.2cm}\\
{\large  College of Humanities \& Sciences, Department of Mathematics,\\
 Nihon University}\end{center}
\begin{abstract}
We give a classification of generic bifurcations of 
intersections of wavefronts generated by different points of a 
hypersurface with or without boundaries.
\end{abstract}
\section{Introduction}
\quad
In \cite{tPleg:cite} we investigated the theory of 
{\it reticular Legendrian unfoldings} 
which describes bifurcations of wavefronts on a hypersurface germ
with a boundary, a corner, or an $r$-corner in a smooth manifold
and studied the stabilities and a generic classification of such bifurcations.

In this paper we give a semi-local classification of bifurcations of wavefronts
generated by a hypersurface with an $r$-corner.
The wavefronts generated around several points on the initial hypersurface
intersect.
For example we consider the below figures in which the initial wavefront with 
a boundary in a plane and wavefronts generated by the initial wavefront and the boundary
to normal directions  are described respectively.
\begin{figure}[htbp]
  \begin{center}
    \includegraphics[width=15cm,height=4.5cm]{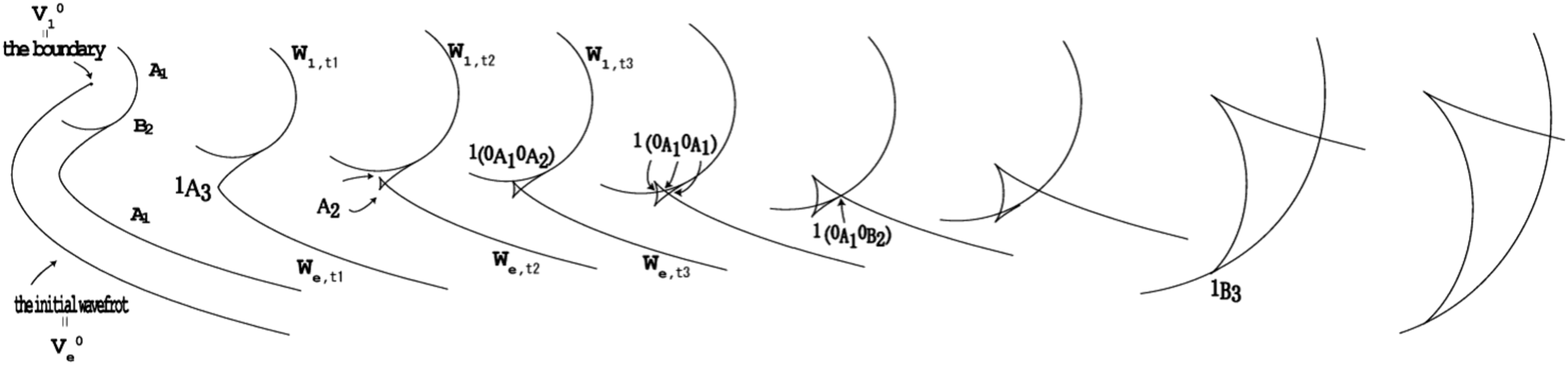}
  \end{center}
\caption{The initial wavefront $V^0$ with a boundary and 
generated wavefronts {\rm(} $e=\emptyset$, $t_1<t_2<t_3${\rm)}}\label{pre:fig}
\end{figure}  
The generated wavefronts bifurcate and typical shapes of bifurcations occur.
The ${}^1({}^0A_1{}^0A_2)$, ${}^1({}^0A_1{}^0A_1)$,
and ${}^1({}^0A_1{}^0B_2)$ fronts are intersections of fronts generated 
by different points of the initial wavefront.
The purpose of this paper is the classification of such generic bifurcations of
intersections of wavefronts.

As the bifurcation theory of wavefronts are described by 
reticular Legendrian unfoldings, we shall introduce the notion of 
{\it multi-reticular Legendrian unfoldings} in order to describe 
the theory of generic bifurcations of intersections of wavefronts generated by 
different points of a hypersurface with an $r$-corner.

The main result of this paper is a classification of generic bifurcations of 
intersections of wavefronts generated
by a hypersurface without boundary or with boundaries in a manifold with dimension
$2$ and $3$.
The classification list of generating families of generic multi-reticular Legendrian unfoldings
are given before Theorem \ref{genericclass}.
We also draw all figures of such bifurcations at the last of this paper.

This paper consists of four sections.
In section \ref{unfold:sec} we investigate the theory of map germs 
with respect to {\it the reticular $t$-$({\cal P}$-${\cal K})_{(m)}$-equivalence relation}
which plays important roles as generating families of multi-reticular Legendrian unfoldings.
In section \ref{RetLegunf:sec}, we review the theory of reticular Legendrian unfoldings which is developed in \cite{tPleg:cite}.
In section \ref{MrLu:sec} we introduce the notion of multi-reticular Legendrian unfoldings and investigate their equivalence relation, generating families, and stabilities.
In section \ref{generic:sec} we reduce our investigation to finitely dimensional jet spaces and give a generic classification of multi-reticular Legendrian unfoldings.
\section{Stabilities of unfoldings}\label{unfold:sec}
\quad
In this section we investigate the theory of map germs  with respect to 
{\it the reticular $t$-$({\cal P}$-${\cal K})_{(m)}$-equivalence relation}
which is proved by almost parallel methods of \cite{tPKfunct}.

\quad 
We denote by ${\cal E}(r;k_1,r;k_2)$ the set of all germs at $0$ of
smooth maps ${\mathbb H}^r\times {\mathbb R}^{k_1} \rightarrow 
{\mathbb H}^r\times {\mathbb R}^{k_2}$ and set ${\mathfrak M}(r;k_1,r;k_2)=
\{ f\in {\cal
E}(r;k_1,r;k_2)|f(0)=0 \}$.
We denote ${\cal E}(r;k_1,k_2)$ for ${\cal E}(r;k_1,0;k_2)$ and 
denote ${\mathfrak M}(r;k_1,k_2)$ for ${\mathfrak M}(r;k_1,0;k_2)$.

 If $k_2=1$ we write simply ${\cal E}(r;k)$ for 
${\cal E}(r;k,1)$
and ${\mathfrak M}(r;k)$ for ${\mathfrak M}(r;k,1)$. 
Then ${\cal E}(r;k)$ is an ${\mathbb R}$-algebra in the usual
way and ${\mathfrak M}(r;k)$ is its unique maximal ideal. 
We also denote by ${\cal E}(k)$ for 
${\cal E}(0;k)$
and ${\mathfrak M}(k)$ for ${\mathfrak M}(0;k)$.
We remark that ${\cal E}(r;k,p)$ is 
an ${\cal E}(r;k)$-module generated by $p$-elements..

We
denote by $J^l(r+k,p)$ the set of $l$-jets at $0$ of germs in ${\cal
E}(r;k,p)$. There are natural projections:
\[ \pi_l:{\cal E}(r;k,p)\longrightarrow J^l(r+k,p),
\pi^{l_1}_{l_2}:J^{l_1}(r+k,p)\longrightarrow J^{l_2}(r+k,p)\ (l_1 > l_2).  \]
We write $j^lf(0)$ for $\pi_l(f)$ for each $f\in {\cal E }(r;k,p)$.

Let $(x,y)=(x_1,\cdots,x_r,y_1,\cdots,y_k)$ be a fixed
coordinate system of $({\mathbb H}^r\times {\mathbb R}^k,0)$. 
We denote by 
${\cal B}(r;k)$ the group of
diffeomorphism germs $({\mathbb H}^r\times {\mathbb R}^{k},0)\rightarrow 
({\mathbb H}^r\times {\mathbb R}^{k},0)$ of the form:
\[ \phi(x,y)=(x_1\phi_1^1(x,y),\cdots,x_r\phi_1^r(x,y),\phi_2^1(x,y),\cdots,\phi_2^k(x,y)
). \]

We denote by 
${\cal B}_n(r;k+n)$ the group of
diffeomorphism germs $({\mathbb H}^r\times {\mathbb R}^{k+n},0)\rightarrow 
({\mathbb H}^r\times {\mathbb R}^{k+n},0)$ of the form:
\begin{eqnarray*}
\lefteqn{\phi(x,y,u)=(x_1\phi_1^1(x,y,u),\cdots,x_r\phi_1^r(x,y,u),}& &\\
& & \hspace{3cm}\phi_2^1(x,y,u),\cdots,\phi_2^k(x,y,u)
,\phi_3^1(u),\ldots,\phi_3^n(u)).
\end{eqnarray*}
We denote $\phi(x,y,u)=(x\phi_1(x,y,u),\phi_2(x,y,u),\phi_3(u))$, 
$\frac{\partial f_0}{\partial y}=(\frac{\partial f_0}{\partial y_1},$ $\cdots,\frac{\partial f_0}{\partial y_k})$, and denote other notations analogously.
\begin{lem}\label{gw1.8:cor}{\rm (cf. \cite[Corollary 1.8]{spsing})}
Let $B$ be a submodule of ${\cal E}(r;k+n+{m}',m)$,
$A_1$ be a finitely generated 
${\cal E}(m')$-submodule of ${\cal E}(r;k+n+{m}',m)$ generated $d$-elements, and 
$A_2$ be a finitely generated ${\cal E}(n+{m}')$ submodule of ${\cal E}(r;k+n+{m}',m)$.
Suppose 
\begin{eqnarray*}
{\cal E}(r;k+n+{m}',m)=B+A_2+A_1+{\mathfrak M}(m'){\cal E}(r;k+n+{m}',m)
\hspace{1cm}\\
+{\mathfrak M}(n+{m}')^{d+1}{\cal E}(r;k+n+{m}',m). 
\end{eqnarray*}
Then 
\[{\cal E}(r;k+n+{m}',m)=B+A_2+A_1,\]
\[ {\mathfrak M}(n+{m}')^d{\cal E}(r;k+n+{m}',m)\subset B+A_2+{\mathfrak M}(m'){\cal E}(r;k+n+{m}',m).\]
\end{lem}

We say that $f_0=(f_{0,1},\ldots,f_{0,m})(x,y),g_0=(g_{0,1},\ldots,g_{0,m})(x,y)\in{\cal E}(r;k+n)$ are {\it reticular ${\cal K}_{(m)}$-equivalent} if
there exist $\Phi_i\in{\cal B}_n(r;k+n)$ and a unit 
$\alpha_i\in {\cal E}(r;k+n)$
such that $g_{0,i}=\alpha_i\cdot f_{0,i}\circ \Phi_i$ for $i=1,\ldots,m$.
We denote $\Phi=(\Phi_1,\ldots,\Phi_{m}),\ \alpha=(\alpha_1,\ldots,\alpha_{m})$ and we call $(\Phi,\alpha)$ a reticular ${\cal K}_{(m)}$-isomorphism from $f$ to $g$. 
We remark that $f_0$ and $g_0$ are reticular ${\cal K}_{(m)}$-equivalent if 
and only if $f_{0,i}$ and $g_{0,i}$ are reticular ${\cal K}$-equivalent
for $i=1,\ldots,m$.
\begin{lem}
Let $f_0(x,y)\in {\mathfrak M}(r;k,m)$ and $z=j^lf_0(0)$. 
Let $O_{r{\cal K}}^l(z)$ be the submanifold of 
$J^l(r+k,m)$ 
consist of the image by $\pi_l$ of the orbit under the 
reticular ${\cal K}_{(m)}$-equivalence of $f_0$. Put $z=j^lf_0(0)$. Then 
\begin{eqnarray*}
T_z(O^l_{r{\cal K}}(z)) = \pi_l((\langle f_{0,1},x\frac{\partial f_{0,1}}{\partial x}
\rangle_{{\cal E}(r;k)}+{\mathfrak M}(r;k)
\langle\frac{\partial f_{0,1}}{\partial y}
\rangle)\hspace{2cm}\\ 
\times\cdots\times
(\langle f_{0,m},x\frac{\partial f_{0,m}}{\partial x}
\rangle_{{\cal E}(r;k)}+{\mathfrak M}(r;k)
\langle\frac{\partial f_{0,m}}{\partial y}
\rangle)
). 
\end{eqnarray*}
\end{lem}

We say that a map germ $f_0=(f_{0,1},\ldots,f_{0,m})(x,y)\in {\mathfrak M}(r;k,m)$ is {\it reticular ${\cal K}_{(m)}$-$l$-determined} 
if all map germ in ${\mathfrak M}(r;k,m)$ which has the same $l$-jet of $f_0$ is 
reticular ${\cal K}_{(m)}$-equivalent to $f_0$. 
\begin{lem}\label{findetc:lm}{\rm (Cf., \cite[Lemma 2.3]{tPKfunct})}
Let $f_0=(f_{0,1},\ldots,f_{0,m})(x,y)\in {\mathfrak M}(r;k,m)$ and let 
\begin{eqnarray*}
{\mathfrak M}(r;k)^{l+1}\subset {\mathfrak M}(r;k)(\langle f_{0,i},x \frac{\partial f_{0,i}}{\partial
x_1}\rangle +{\mathfrak M}(r;k)\langle
\frac{\partial f_{0,i}}{\partial y}\rangle )
+{\mathfrak M}(r;k)^{l+2}
\end{eqnarray*}
for $i=1,\ldots,m$. 
Then $f_0$ is reticular ${\cal K}_{(m)}$-$l$-determined. 
Conversely if $f_0(x,y)\in {\mathfrak M}(r;k,m)$ is reticular ${\cal K}_{(m)}$-$l$-determined, then 
\[ {\mathfrak M}(r;k)^{l+1}\subset \langle f_{0,i},x \frac{\partial f_{0,i}}{\partial x_1}\rangle_{ {\cal E}(r;k) }  +{\mathfrak M}(r;k)\langle
\frac{\partial f_{0,i}}{\partial y}\rangle \ \mbox{ for } i=1,\ldots,m. \]
\end{lem}   

We say that $f=(f_1,\ldots,f_m)(x,y,u),g=(g_1,\ldots,g_m)(x,y,u)
\in{\cal E}(r;k+n,m)$ are {\it reticular $({\cal P}$-${\cal K})_{(m)}$-equivalent} if
there exist $\Phi_i\in{\cal B}_n(r;k+n)$ of the form:
\[ \Phi_i(x,y,u)=(x\phi^i_1(x,y,t,u),\phi^i_2(x,y,t,u),\phi_3(u)) \]
and a unit $\alpha=(\alpha_1,\ldots,\alpha_m)(x,y,u)\in {\cal E}(r;k+n,m)$
such that $g_i=\alpha_i\cdot f_i\circ \Phi_i$ for $i=1,\ldots,m$.
We write $\Phi=(\Phi_1,\ldots,\Phi_m)$ and call $(\Phi,\alpha)$ 
a reticular $({\cal P}$-${\cal K})_{(m)}$-isomorphism from $f$ to $g$. \\

We say that a map germ $f=(f_{1},\ldots,f_{m})(x,y,u)\in {\mathfrak M}(r;k+n,m)$ is {\it reticular $({\cal P}$-${\cal K})_{(m)}$-$l$-determined} 
if all map germ in ${\mathfrak M}(r;k+n,m)$ which has the same $l$-jet of $f$ is 
reticular $({\cal P}$-${\cal K})_{(m)}$-equivalent to $f$. 
\begin{lem}
Let $f=(f_1,\ldots,f_m)(x,y,u)\in {\cal E}(r;k+n,m)$ and $z=j^lf(0)$. 
We write $O^l_{r{\cal P}\mbox{-}{\cal K}_{(m)}}(z)$ the 
submanifold of the image of the orbit of $f$ under the reticular 
$t$-$({\cal P}$-${\cal K})_{(m)}$ equivalence relation.
Then 
\begin{eqnarray*}
O^l_{r{\cal P}\mbox{-}{\cal K}_{(m)}}(z)=\pi_l((\langle f_1,x\frac{\partial f_1}{\partial x}\rangle_{{\cal E}(r;k+n)}+{\mathfrak M}(r;k+n)\langle
\frac{\partial f_1}{\partial y}
\rangle)\hspace{1cm}\\
\times\cdots\times
(\langle f_{m},x\frac{\partial f_{m}}{\partial x}\rangle_{{\cal E}(r;k+n)}+
{\mathfrak M}(r;k+n)\langle
\frac{\partial f_{m}}{\partial y}
\rangle)+{\mathfrak M}(n)\langle \frac{\partial f}{\partial u}\rangle).
\end{eqnarray*}
\end{lem}
\begin{lem}\label{nmqdet}{\rm (Cf., \cite[Lemma 3.10]{tPKfunct})}
Let $f=(f_1,\ldots,f_m)(x,y,u)\in {\mathfrak M}(r;k+n,m)$ and $l$ be a non-negative integer.
If 
\begin{eqnarray}
{\mathfrak M}(r;k+n,m)^{l}\subset (\langle f_1,x \frac{\partial f_1}{\partial
x}\rangle_{{\cal E}(r;k+n)} 
+{\mathfrak M}(r;k+n)\langle
\frac{\partial f_1}{\partial y}\rangle)\times\cdots\times\nonumber\\
(\langle f_m,x \frac{\partial f_m}{\partial
x}\rangle_{{\cal E}(r;k+n)} 
+{\mathfrak M}(r;k+n)\langle
\frac{\partial f_m}{\partial y}\rangle) 
+{\mathfrak M}(n)\langle
\frac{\partial f}{\partial u}\rangle\nonumber\\
+{\mathfrak M}(n) {\mathfrak M}(r;k+n,m)^l,\label{qqq:ren}
\end{eqnarray}
then $f$ is reticular 
$({\cal P}$-${\cal K})_{(m)}$-$l$-determined. 
\end{lem}
%
%

In convenience, we denote an unfolding of  a function germ 
$f(x,y,u)$ $\in {\mathfrak M}(r;k+n,m)$ by
$F(x,y,t,u)\in {\mathfrak M}(r;k+{m}'+n,m)$.

Let 
$F(x,y,t,u)\in {\mathfrak M}(r;k+{m}'_1+n,m)$ and 
$G(x,y,s,u)\in {\mathfrak M}(r;k+{m}'_2+n,m)$ be
unfoldings of $f(x,y,u)\in {\mathfrak M}(r;k+n,m)$.

{\it A reticular $t$-$({\cal P}$-${\cal K})_{(m)}$-$f$-morphism from $F$ to $G$} is a pair $(\Phi,\alpha )$,
where $\Phi=(\Phi_1,\ldots,\Phi_{m})$
for $\Phi_i\in {\mathfrak M}(r;k+{m}'_2+n,r;k+{m}'_1+n)$ and $\alpha=(\alpha_1,\ldots,\alpha_{m})$ is a unit of ${\cal E}(r;k+{m}'_2+n,m)$
satisfying the following conditions:\\ 
(1) $\Phi_i$ can be written in the form:
\[
\Phi_i(x,y,s,u)=(x\phi^i_1(x,y,s,u),\phi^i_2(x,y,s,u),\phi_3(s),\phi_4(s,u)),\]
(2) $\Phi_i|_{{\mathbb H}^{r}\times{\mathbb R}^{k+n}}=id_{{\mathbb H}^{r}
\times{\mathbb R}^{k+n}}$, 
$\alpha_i|_{{\mathbb H}^{r}\times{\mathbb R}^{k+n}}\equiv 1$,\\
(3) $G_i(x,y,s,u)=\alpha_i(x,y,s,u) \cdot F_i\circ \Phi_i(x,y,s,u)$ for all 
$(x,y,s,u)\in ({\mathbb H}^{r}\times{\mathbb R}^{k+{m}'_2+n},0)$.

If there exists a reticular $t$-$({\cal P}$-${\cal K})_{(m)}$-$f$-morphism from 
$F$ to $G$,
 we say that
{\it $G$ is reticular $t$-$({\cal P}$-${\cal K})_{(m)}$-$f$-induced from $F$}. 
If $m_1=m_2$ and $\Phi$ is
invertible, we call $(\Phi,\alpha )$ {\it a reticular $t$-$({\cal P}$-${\cal K})_{(m)}$-$f$-isomorphism from $F$ to $G$}
 and we say that $F$ is reticular $t$-$({\cal P}$-${\cal K})_{(m)}$-$f$-equivalent 
to $G$.\vspace{5mm}

We
say that $F(x,y,t,u),G(x,y,t,u)\in {\cal E}(r;k+{m}'+n,m)$
are {\it reticular $t$-$({\cal P}$-${\cal K})_{(m)}$-equivalent} if 
there exist diffeomorphism germs 
\[ \Phi_i:({\mathbb H}^{r}\times {\mathbb R}^{k+{m}'+n},0)\rightarrow
({\mathbb H}^{r}\times {\mathbb R}^{k+{m}'+n},0)\ \ (i=1,\ldots,m)\]
of the form
\[ \Phi_i(x,y,t,q,z)=(x\phi^i_1(x,y,t,u),\phi^i_2(x,y,t,u),\phi_3(t),
\phi_4(t,u))\]
and a unit $\alpha=(\alpha_1,\ldots,\alpha_{m})\in {\cal E}(r;k+{m}'+n,m)$
such that $G_i=\alpha_i\cdot F_i\circ \Phi_i$ for $i=1,\ldots,m$.
We call $(\Phi,\alpha )$ {\it a
reticular $t$-$({\cal P}$-${\cal K})_{(m)}$-isomorphism from $F$ to $G$}. \\

We say that function germs $F_0(x_1,\ldots,x_{r_1},y_1,\cdots,y_{k_1},t,q,z)\in 
{\mathfrak M}$ $(r_1;k_1+1+n+1)$ 
and $F_1(x_1,\ldots,x_{r_2},y_1,\cdots,y_{k_2},t,q,z)\in {\mathfrak M}(r_2;k_2+1+n+1)$ 
are {\it stably 
 reticular $t$-${\cal P}$-${\cal K}$-equivalent} 
if $F_0$ and $F_1$ are reticular $t$-${\cal P}$-${\cal K}$-equivalent 
after additions of linear forms of $x$ of which all coefficients are not zero and non-degenerate quadratic forms in the variables $y$. \\

We  define the notion of {\it 
stable reticular $t$-$({\cal P}$-${\cal K})_{(m)}$-equivalence} 
in the same way as the 
stable reticular $t$-${\cal P}$-${\cal K}$-equivalence.

\begin{dfn}{\rm 
We define stabilities of unfoldings.
Let $F=(F_1,$ $\ldots,F_m)(x,y,t,u)\in {\mathfrak M}(r;k+{m}'+n,m)$ be an 
unfolding of $f=(f_1,$ $\ldots,f_m)(x,y,u)\in {\mathfrak M}(r;k+n,m)$.\\

We say that $F$ is {\it reticular $t$-$({\cal P}$-${\cal K})_{(m)}$-stable}
if the following condition holds:
For any neighbourhood $U$ of $0$ in ${\mathbb H}^{r}\times {\mathbb R}^{k+{m}'+n}$ 
and any representative $\tilde{F}\in 
C^\infty (U,{\mathbb R}^m)$, there exists a neighbourhood $N_{{\tilde{F}}}$
of $\tilde{F}$ in $C^\infty$-topology such that 
for any $\tilde{G}\in N_{{\tilde{F}}}$
there exist $(0,y^i,t^0,q^0,z^0)\in U$ for $i=1,\ldots,m$ such that 
$G=(G_1,\ldots,G_{m})$ and $F$ are 
reticular $t$-$({\cal P}$-${\cal K})_{(m)}$-equivalent
where $G_i\in {\mathfrak M}(r;k+{m}'+n)$ is defined by 
$G_i(x,y,t,u)=\tilde{G}(x,y+y^i,t+t^0,u+u^0)-\tilde{G}(0,y,
t^0,u^0)$.\\

We say that $F$ is {\it a reticular
$t$-$({\cal P}$-${\cal K})_{(m)}$-versal unfolding of $f$} if any unfolding of $f$ is reticular $t$-$({\cal P}$-${\cal K})_{(m)}$-$f$-induced
from $F$. \\

We say that $F$ is {\it a reticular
$t$-$({\cal P}$-${\cal K})_{(m)}$-universal unfolding of} $f$ if $m$ is minimal in 
reticular
$t$-$({\cal P}$-${\cal K})_{(m)}$-versal unfoldings of $f$.\\

We say that $F$ is {\it reticular $t$-$({\cal P}$-${\cal K})_{(m)}$-infinitesimally versal} if 
\begin{eqnarray*}
{\cal E}(r;k+n,m)=
\langle f_1,x\frac{\partial f_1}{\partial x},\frac{\partial f_1}{\partial y}
\rangle_{{\cal E}(r;k+n)}\times\cdots\times\hspace{2cm}\\
\langle f_{m},x\frac{\partial f_{m}}{\partial x},\frac{\partial f_{m}}{\partial y}
\rangle_{{\cal E}(r;k+n)}
+\langle \frac{\partial f}{\partial u}\rangle_{{\cal E}(n)}+
\langle \frac{\partial F}{\partial t}|_{t=0}\rangle_{{\mathbb R}}
\end{eqnarray*}

We say that $F$ is {\it reticular $t$-$({\cal P}$-${\cal K})_{(m)}$-infinitesimally stable} if 
\begin{eqnarray*}
{\cal E}(r;k+{m}'+n,m)=
\langle F_1,x\frac{\partial F_1}{\partial x},\frac{\partial F_1}{\partial y}
\rangle_{{\cal E}(r;k+{m}'+n)}\times\cdots\times\hspace{2cm}\\
\langle F_{m},x\frac{\partial F_{m}}{\partial x},\frac{\partial F_{m}}{\partial y}
\rangle_{{\cal E}(r;k+{m}'+n)}
+\langle \frac{\partial F}{\partial u}\rangle_{{\cal E}(m'+n)}+
\langle \frac{\partial F}{\partial t}\rangle_{{\cal E}(m')}
\end{eqnarray*}

We say that $F$ is
{\it reticular $t$-$({\cal P}$-${\cal K})_{(m)}$-homotopically stable} if for any smooth path-germ $({\mathbb 
R},0)\rightarrow {\cal E}(r;k+{m}'+n,m),\tau\mapsto F_\tau=(F_{1,\tau},\ldots,F_{m,\tau})$ with $F_0=F$, there exists a
smooth path-germs $({\mathbb R},0)\rightarrow {\cal B}(r;k+{m}'+n)\times {\cal
E}(r;k+{m}'+n),\tau\mapsto (\Phi^i_\tau,\alpha^i_\tau)$ with $(\Phi^i_0,\alpha^i_0)=(id,1)$ 
and $\Phi^i_\tau$ has the form
\[ \Phi^i_\tau(x,y,t,u)=(x\phi_\tau^{i,1}(x,y,t,u),\phi_\tau^{i,2}(x,y,t,u),
\phi_\tau^3(t),\phi_\tau^4(t,u)) \]
such that each
$(\Phi^i_\tau,\alpha^i_\tau)$ is a 
reticular $t$-${\cal P}$-${\cal K}$-isomorphism and $F_{i,\tau}=\alpha^i_\tau\cdot
F_{i,0}\circ \Phi^i_\tau$ for 
$\tau\in ({\mathbb R},0)$ and $i=1,\dots,m$. 
}\end{dfn}

Let $U$ be a neighbourhood of $0$ in ${\mathbb H}^{r}\times {\mathbb R}^{k+{m}'+n}$ 
and let $\tilde{F}=(\tilde{F}^1,\ldots,\tilde{F}^m):U\rightarrow {\mathbb R}^m$ be a smooth map and 
$l$ be a non-negative integer. 
We choose a neighbourhood  $U'$ of  $0$ in 
${\mathbb R}^{km+m'+n}$ such that 
$(0,y^i,t,u)\in U$ for any $(y^1,\ldots,y^{m},t,u)\in U$ and $i$.
We define the smooth map germ
\[ j^{l}_1\tilde{F}:U'\longrightarrow J^{l}(r+k+n)  \]
as the follow: For $(y^1,\ldots,y^{m},t,u)\in U$ we set 
$j^{l}_1\tilde{F}(y^1,\ldots,y^{m},t,u)$ by 
the $l$-jet of the map germ 
$(\tilde{F}^1_{(y^1,t,u)},\ldots,\tilde{F}^m_{(y^m,t,u)})\in {\mathfrak M}(r;k+n,m)$ at $0$, 
where  $\tilde{F}^i_{(y^i,t,u)}$ is given by 
$\tilde{F}^i_{(y^i,t,u)}(x',y',u')=\tilde{F}(x',y^i+y',t,u+u')-\tilde{F}^i(0,y^i,t,u)$\\

Let $F(x,y,t,u)\in {\mathfrak M}(r;k+{m}'+n,m)$ be an unfolding of 
$f(x,y,u)\in {\mathfrak M}(r;k+n,m)$. 
Let $l$ be a non-negative integer and $z=j^{l}f(0)$. 
We say that $F$
is {\it reticular $t$-$({\cal P}$-${\cal K})_{(m)}$-$l$-transversal} if 
$j^{l}_1\tilde{F}$ at $0$ is transversal to
$O_{r{\cal P}\mbox{-}{\cal K}_{(m)}}^{l}(z)$
for a representative $\tilde{F}\in C^\infty (U,{\mathbb R}^m)$ of $F$.
\begin{lem}{\rm (Cf., \cite[Lemma 3.4]{tPKfunct})}
Let $F(x,y,t,u)\in {\mathfrak M}(r;k+{m}'+n,m)$ be an unfolding of 
$f(x,y,u)\in {\mathfrak M}(r;k+n,m)$. 
Then $F$ is reticular $t$-$({\cal P}$-${\cal K})_{(m)}$-$l$-transversal unfolding of $f$ if and only if 
\begin{eqnarray*}
{\cal E}(r;k+n,m)=
\langle f_1,x\frac{\partial f_1}{\partial x},\frac{\partial f_1}{\partial y}
\rangle_{{\cal E}(r;k+n)}\times\cdots\times\hspace{1cm}\\
\langle f_{m},x\frac{\partial f_{m}}{\partial x},\frac{\partial f_{m}}{\partial y}
\rangle_{{\cal E}(r;k+n)}
+\langle \frac{\partial f}{\partial u}\rangle_{{\cal E}(n)}+
\langle \frac{\partial F}{\partial t}|_{t=0}\rangle_{{\mathbb R}}\\
+{\mathfrak M}(r;k+n)^{l+1}{\cal E}(r;k+n,m).
\end{eqnarray*}

\end{lem}
\begin{tth}\label{univunf:tth}{\rm (Uniqueness of universal unfoldings)}{\rm (Cf., \cite[Theorem 3.13]{tPKfunct})}
Let $F(x,y,u,$ $t),G(x,y,u,t)\in{\mathfrak M}(r;k+n+m',m)$ be unfoldings of 
 $f\in {\mathfrak M}(r;k+n,m)$.
If $F$ and $G$ are reticular $t$-$({\cal P}$-${\cal K})_{(m)}$-universal,
then $F$ and $G$ are reticular 
$t$-${\cal P}$-${\cal K}$-$f$-equivalent.
\end{tth}
\begin{tth}\label{mthft:th}{\rm (Cf., \cite[Theorem 3.14]{tPKfunct})}
Let $F(x,y,t,u)\in {\mathfrak M}(r;k+{m}'+n,m)$ be an unfolding of 
$f(x,y,u)\in {\mathfrak M}(r;k+n,m)$ 
and let $f$ is an unfolding of $f_0(x,y)\in {\mathfrak M}(r;k,m)$.  
Then 
following are equivalent. \\ 
{\rm (1)} There exists a non-negative number $l$ such that 
$f_0$ is reticular 
${\cal K}_{(m)}$-$l$-determined and  $F$ is reticular
$t$-$({\cal P}$-${\cal K})_{(m)}$-$l'$-transversal for $l'\geq lm+l+{m}'$,\\
{\rm (2)} $F$ is reticular $t$-$({\cal P}$-${\cal K})_{(m)}$-stable,\\
{\rm (3)} $F$ is reticular $t$-$({\cal P}$-${\cal K})_{(m)}$-versal,\\ 
{\rm (4)} $F$ is reticular
$t$-$({\cal P}$-${\cal K})_{(m)}$-infinitesimally versal,\\ 
{\rm (5)} $F$ is reticular
$t$-$({\cal P}$-${\cal K})_{(m)}$-infinitesimally stable,\\ 
{\rm (6)} $F$ is reticular
$t$-$({\cal P}$-${\cal K})_{(m)}$-homotopically stable.
\end{tth}
\section{Reticular Legendrian unfoldings}\label{RetLegunf:sec}
\quad
We review the theory of reticular Legendrian unfoldings.
Let $J^1({\mathbb R}^n,{\mathbb R})$ be the $1$-jet bundle of functions 
in $n$-variables which may be considered as ${\mathbb R}^{2n+1}$ with 
a natural coordinate system $(q,z,p)=(q_1,\ldots,$ $q_n,z,p_1,\ldots,p_n)$,
where $(q_1,\ldots,q_n)$ is a coordinate system of ${\mathbb R}^n$.
We equip the contact structure on $J^1({\mathbb R}^n,{\mathbb R})$
defined by the canonical $1$-form $\theta=dz-\sum_{i=1}^np_idq_i$.
We have a natural projection 
$\pi:J^1({\mathbb R}^n,{\mathbb R})\rightarrow {\mathbb R}^n\times {\mathbb R}$
by $\pi(q,z,p)=(q,z)$.

We remark that: By the identification of the coordinate system $(u_1,\ldots,u_{n+1})$ 
and $(q_1,\ldots,q_n,z)$ on a manifold, the dimension of a manifold on which wavefronts 
propagate is higher $1$ than $n$.

We also consider the $1$-jet bundle 
$J^1({\mathbb R}\times{\mathbb R}^n,{\mathbb R})$
and the canonical $1$-form $\Theta$ on that space.
Let $(t,q_1,\ldots,q_n)$ be the canonical coordinate system on 
${\mathbb R}\times{\mathbb R}^n$ and 
\[ (t,q_1,\ldots,q_n,z,s,p_1,\ldots,p_n) \]
be the corresponding coordinate system on 
$J^1({\mathbb R}\times{\mathbb R}^n,{\mathbb R})$.
Then the canonical $1$-form $\Theta$ is given by 
\[ \Theta=dz-\sum_{i=1}^np_idq_i-sdt=\theta-sdt. \]
We define the natural projection 
$\Pi:J^1({\mathbb R}\times{\mathbb R}^n,{\mathbb R})\rightarrow
({\mathbb R}\times{\mathbb R}^n)\times {\mathbb R}$
by $\Pi(t,q,z,s,p)=(t,q,z)$.\\

Let   
$\tilde{L}^{r,0}_\sigma=\{(t,q,z,s,p) \in J^1({\mathbb R}\times 
{\mathbb R}^n,{\mathbb R})|
q_\sigma=p_{I_r-\sigma}=q_{r+1}=\cdots=q_n=s=z=0,q_{I_r-\sigma}\geq 0 \}$
for $\sigma\subset I_r=\{ 1,2,\cdots,r \}$ and let 
${\mathbb L}^r=\{(t,q,z,s,p) \in J^1({\mathbb R}\times 
{\mathbb R}^n,{\mathbb R})|
q_1p_1=\cdots=q_rp_r=q_{r+1}=\cdots=q_n=s=z=0,q_{I_r}\geq 0 \}$ be a 
representative as a germ 
of the union of $\tilde{L}^{r,0}_\sigma$ for all $\sigma\subset I_r$.

We recall that $\tilde{L}^{r,0}_\sigma$ is a normalisation of the 
Legendrian submanifold 
consists of the particles generated
by the {\it $\sigma$-corner} of hypersurface germ  $V=\{ (u_1,\ldots,u_n,0)
\in ({\mathbb R}^{n+1},0) |
u_{I_r}\geq 0\}$ with an $r$-corner to conormal direction, where $\sigma$-corner is defined by 
$V\cap \{u_\sigma=0\}$ and $u$ is a local coordinate system of 
 a manifold in which
wavefronts propagate.
\begin{dfn}{\rm 
Let $C:(J^1({\mathbb R}\times {\mathbb R}^n,{\mathbb R}),0)\rightarrow 
(J^1({\mathbb R}\times {\mathbb R}^n,{\mathbb R}),w)\ ($ $\Pi(w)=0)$ be 
a contact diffeomorphism germ.
We say that $C$ is a {\it ${\cal P}$-contact diffeomorphism germ} if $C$ has the 
form:
\begin{eqnarray}
C(t,q,z,s,p)=(t,q_t(q,z,p),z_t(q,z,p),h(t,q,z,s,p),p_t(q,z,p)). \label{P-contaform}
\end{eqnarray}
We say a map germ ${\cal L}:({\mathbb L}^r,0)\rightarrow 
(J^1({\mathbb R}\times {\mathbb R}^n,{\mathbb R}),w)\ (\Pi(w)=0)$ 
{\it a reticular Legendrian unfolding} if there exists a ${\cal P}$-contact diffeomorphism germ $C$ such that ${\cal L}$ is the restriction of $C$ to ${\mathbb L}^r$.

We call $\{ {\cal L}(\tilde{L}^{r,0}_\sigma) \}_{\sigma\subset I_r}$ 
{\it the unfolded contact regular $r$-cubic configuration of } ${\cal L}$.
}\end{dfn}

In order to study perestroikas of wavefronts of reticular Legendrian 
unfoldings, 
we introduce the following equivalence relation.
Let ${\cal L}_i:({\mathbb L}^r,0)\rightarrow 
(J^1({\mathbb R}\times {\mathbb R}^n,{\mathbb R}),w_i), (i=0,1)$ 
be reticular Legendrian unfoldings.
We say that ${\cal L}_1$ and ${\cal L}_2$ are 
{\it ${\cal P}$-Legendrian equivalent} 
if there exist a contact diffeomorphism germ 
\[ K:(J^1({\mathbb R}\times {\mathbb R}^n,{\mathbb R}),w_1)\rightarrow 
(J^1({\mathbb R}\times {\mathbb R}^n,{\mathbb R}),w_2) \]
of the form 
\begin{eqnarray}
K(t,q,z,s,p)=(\phi_1(t),\phi_2(t,q,z),\phi_3(t,q,z),\hspace{3cm}\nonumber\\
\phi_4(t,q,z,s,p),
\phi_5(t,q,z,s,p))\label{PLequi}
\end{eqnarray}
and {\it a reticular $r$-diffeomorphism} 
$\Psi$ on 
$(J^1({\mathbb R}\times {\mathbb R}^n,{\mathbb R}),0)$
such that 
$K\circ {\cal L}_1={\cal L}_2\circ \Psi$,
where a reticular $r$-diffeomorphism is defined by a contact diffeomorphism 
$\Psi$ on 
$(J^1({\mathbb R}\times {\mathbb R}^n,{\mathbb R}),0)$
such that $t\circ\Psi$ depends only on $t$ and $\Psi$ preserves 
$\tilde{L}^{r,0}_\sigma$ for all $\sigma\subset I_r$. \\

We can construct generating families of reticular Legendrian unfoldings.
A function germ $F(x,y,t,q,z)\in{\mathfrak M}(r;k+1+n+1)$ is called
{\it ${\cal P}$-$C$-non-degenerate} if 
$\frac{\partial F}{\partial x}(0)=\frac{\partial F}{\partial y}(0)=0,
\frac{\partial F}{\partial z}(0)\neq 0$ and 
\[ x_1,\cdots,x_r,t,F,\frac{\partial F}{\partial x_1},\cdots,\frac{\partial F}{\partial x_r},\frac{\partial F}{\partial y_1},\cdots,\frac{\partial F}{\partial y_k} \]
are independent on $({\mathbb H}^k\times {\mathbb R}^{k+1+n+1},0)$.

Let ${\cal L}$ be a reticular Legendrian unfolding.
A ${\cal P}$-$C$-non-degenerate function germ 
$F(x,y,$ $t,q,z)\in {\mathfrak M}(r;k+1+n+1)$
is called {\it a generating family} of ${\cal L}$ if 
$F$ is a generating family of the unfolded contact 
regular $r$-cubic configuration $\{ {\cal L}(\tilde{L}^{r,0}_\sigma) 
\}_{\sigma \subset I_r}$, that is
\begin{eqnarray*}
 {\cal L}(\tilde{L}^{r,0}_\sigma)=
\{ (t,q,z,\frac{\partial F}{\partial t}/(-\frac{\partial F}{\partial z}),
\frac{\partial F}{\partial q}/(-\frac{\partial F}{\partial z}))\in
(J^1({\mathbb R}\times {\mathbb R}^n,{\mathbb R}),w)|\\
x_\sigma=F=\frac{\partial F}{\partial x_{I_r-\sigma}}=
\frac{\partial F}{\partial y}=0,x_{I_r-\sigma}\geq 0\} 
\end{eqnarray*}
for all $\sigma\subset I_r$.\\

By \cite[Theorem 0.8.5]{tPleg:cite} we have that:
\begin{tth}\label{UCgf:th}
{\rm (1)} For any reticular Legendrian unfolding ${\cal L}:({\mathbb L}^r,0)
$ $\rightarrow (J^1({\mathbb R}\times {\mathbb R}^n,{\mathbb R}),w)$, 
there exists a function germ $F(x,y,t,q,z)\in{\mathfrak M}(r;k+1+n+1)$ 
which is a generating family of ${\cal L}$.\\
{\rm (2)} For any ${\cal P}$-$C$-non-degenerate function germ 
$F(x,y,t,q,z)\in{\mathfrak M}(r;k+1+n+1)$  there exists a reticular
 Legendrian unfolding ${\cal L}:({\mathbb L}^r,0)$ $\rightarrow 
(J^1({\mathbb R}\times {\mathbb R}^n,{\mathbb R}),w)$ 
of which $F$
is a generating family.\\
{\rm (3)} Two reticular Legendrian unfolding are ${\cal P}$-Legendrian 
equivalent if and only if their generating families are stably reticular 
$t$-${\cal P}$-${\cal K}$-equivalent.
\end{tth}

\section{Multi-reticular Legendrian unfolding}\label{MrLu:sec}
\quad 
Let ${\cal L}_i:({\mathbb L}^{r_i},0)\rightarrow 
(J^1({\mathbb R}\times {\mathbb R}^n,{\mathbb R}),w_i)(i=1,\ldots,m)$ 
be reticular Legendrian unfoldings with $\Pi(w_i)=0$ where 
$w_1,\ldots,w_m$ are distinct.
Then we call ${\cal L}=({\cal L}_1,\ldots,{\cal L}_m)$ {\it a multi-reticular 
Legendrian unfolding}.

Let $({\cal L}_1,\ldots,{\cal L}_m)$ and $({\cal L}'_1,\ldots,{\cal L}'_m)$
be multi-reticular Legendrian unfoldings.
We say that they are reticular ${\cal P}_{(m)}$-Legendrian equivalent
if there exist contact diffeomorphism germs
\[ K_i:(J^1({\mathbb R}\times {\mathbb R}^n,{\mathbb R}),w_i)
\rightarrow (J^1({\mathbb R}\times {\mathbb R}^n,{\mathbb R}),w'_i)
\ \ \ (i=1,\ldots ,m)\]
of the form
\begin{eqnarray*}
\lefteqn{K_i(t,q,z,s,p)=}\\
& &(\phi_1(t),\phi_2(t,q,z),\phi_3(t,q,z),
\phi^i_4(t,q,z,s,p),\phi_5^i(t,q,z,s,p))
\end{eqnarray*}
and reticular ${\cal P}_{r_i}$-diffeomorphisms $\Psi_i$ on $({\mathbb L}^{r_i},0)$
such that $K_i\circ {\cal L}_i={\cal L}'_i\circ \Psi_i$ for 
$i=1,\ldots,m$.

Let $({\cal L}_1,\ldots,{\cal L}_m)$ be a multi-reticular Legendrian unfoldings 
and $F_i\in {\mathfrak M}(r_i;k_i+1+n+1)$ be generating families of 
${\cal L}_i$.
We call $F=(F_1,\ldots,F_m)$ {\it a multi-generating family} of 
$({\cal L}_1,\ldots,{\cal L}_m)$.\\

We say that a map germs $(F_1,\ldots,F_m)
\in{\mathfrak M}(r;k+1+n+1)$
is {\it ${\cal P}_{(m)}$-$C$-non-degenerate} if all $F_i$ are 
${\cal P}$-$C$-non-degenerate and 
$(\frac{\partial F_1}{\partial t}/(-\frac{\partial F_1}{\partial z}),
\frac{\partial F_1}{\partial q}/(-\frac{\partial F_1}{\partial z})),\ldots,
(\frac{\partial F_m}{\partial t}/(-\frac{\partial F_m}{\partial z}),$ $
\frac{\partial F_m}{\partial q}/(-\frac{\partial F_m}{\partial z}))$
are distinct.\\

By Theorem 0.8.5 in \cite{tPleg:cite}  we have that:
\begin{prop}
$(1)$ For any multi-reticular Legendrian unfolding ${\cal L}$,
there exists a multi-generating family of ${\cal L}$,\\
$(2)$ For any ${\cal P}_{(m)}$-$C$-non-degenerate map germs 
$F=(F_1,\ldots,F_m)\ (F_i\in {\mathfrak M}(r_i;k_i+1+n+1))$, 
there exists a multi-reticular Legendrian unfolding of which 
$F$ is a multi-generating family,\\
$(3)$ Let $F=(F_1,\ldots,F_m)$ and $F'=(F'_1,\ldots,F'_m)$ be multi-generating 
families of multi-reticular Legendrian unfoldings 
$({\cal L}_1,\ldots,{\cal L}_m)$ and $({\cal L}'_1,$ $\ldots,{\cal L}'_m)$
respectively.
Then $({\cal L}_1,\ldots,{\cal L}_m)$ and $({\cal L}'_1,\ldots,{\cal L}'_m)$ 
are ${\cal P}_{(m)}$-Legendrian equivalent if and only if
$F$ and $F'$ are stably reticular $t$-$({\cal P}$-${\cal K})_{(m)}$-equivalent.
\end{prop} 

\begin{dfn}{\rm
We define stabilities of multi-reticular Legendrian unfoldings.
Let ${\cal L}=({\cal L}_1,\ldots,{\cal L}_m)$ be
a multi-reticular Legendrian unfolding.\\
{\bf Stability}: We say that ${\cal L}$ is {\it stable}  if the 
following condition holds:
Let $C^{0,i}\in 
C_T(J^1({\mathbb R}\times {\mathbb R}^n,{\mathbb R}),0)\ (i=1,\ldots,m)$ 
be ${\cal P}$-contact embedding germs 
such that $C^{0,i}|_{{\mathbb L}^{r_i}}={\cal L}_i$ and 
$\tilde{C}^{0,i}\in C_T(U_i,J^1({\mathbb R}\times {\mathbb R}^n,{\mathbb R}))$ be
a representative of $C^{0,i}$.
Then there exist open neighbourhoods $N_{\tilde{C}^{0,i}}$ of $\tilde{C}^{0,i}$ 
in $C^\infty$-topology for $i=1,\ldots,m$
such that for any $\tilde{C}^i\in N_{\tilde{C}^{0,i}}$,
there exist points $x_i=(T^i,0,\ldots,0,P^i_{r+1},\ldots,P^i_n)\in U_i$ such that 
the multi-reticular Legendrian unfolding 
$({\cal L}'_{x_1},\ldots,{\cal L}'_{x_m})$ 
and 
${\cal L}$ are 
${\cal P}_{(m)}$-Legendrian equivalent,
where reticular Legendrian unfoldings ${\cal L}'_{x_i}$ are chosen that 
the reticular Legendrian unfolding 
$\tilde{C}^i|_{{\mathbb L}^{r_i}}:
({\mathbb L}^{r_i},x_i) \rightarrow 
(J^1({\mathbb R}\times {\mathbb R}^n,{\mathbb R}),\tilde{C}^i(x_i))$
and ${\cal L}'_{x_i}$
are ${\cal P}$-Legendrian equivalent.
We remark that the definition of stability is not depend on choices of 
${\cal L}'_{x_i}$.\\
{\bf Homotopical stability}: We say that ${\cal L}$ is {\it homotopically stable} if 
for any reticular ${\cal P}$-contact deformations 
$\bar{C^i}=\{ C^i_\tau\}$ of ${\cal L}_i$,
there exist 
 one-parameter families of ${\cal P}$-Legendrian 
  equivalences $\bar{K^i}=\{ K^i_\tau \}$ on 
$(J^1({\mathbb R}\times {\mathbb R}^n,{\mathbb R}),w_i)$
with $K^i_0=id$ of the form 
\begin{eqnarray}
\lefteqn{K^i_\tau(t,q,z,s,p)=}\nonumber\\ & & 
(\phi^1_\tau (t),\phi^2_\tau (t,q,z),\phi^3_\tau (t,q,z),
\phi^{4,i}_\tau (t,q,z,s,p),\phi^{5,i}_\tau (t,q,z,s,p))\label{Khomo:eqn}
\end{eqnarray}
and 
a one-parameter deformation of
reticular ${\cal P}_{r_i}$-diffeomorphisms $\bar{\Psi}^i=\{ \Psi^i_\tau \}$
such that $C^i_\tau=K^i_\tau\circ C^i_0\circ \Psi^i_\tau$ for $t$ around
$0$ and $i=1,\ldots,m$.\\
{\bf Infinitesimal stability}:
We say that ${\cal L}$ is  
 {\it infinitesimally stable} if for any
extension $C^i$ of ${\cal L}_i$ and any infinitesimal ${\cal P}$-contact 
transformation $v_i$ of $C^i$,
there exist infinitesimal reticular ${\cal P}_{r_i}$-diffeomorphisms $\xi_i$ and 
infinitesimal ${\cal P}$-Legendrian equivalences $\eta_i$ at $w_i$
of 
the form
\begin{eqnarray*}
\eta_i(t,q,z,s,p)=a_1(t)\frac{\partial}{\partial t}+
a_2(t,q,z)\frac{\partial}{\partial q}+a_3(t,q,z)\frac{\partial}{\partial z}
\hspace{2cm}\nonumber \\
+
a^i_4(t,q,z,s,p)\frac{\partial}{\partial s}+a^i_5(t,q,z,s,p)
\frac{\partial}{\partial p}  \label{etaform}
\end{eqnarray*}
such  that $v_i=C^i_*\xi_i+\eta_i\circ C^i$ for $i=1,\ldots,m$.
}\end{dfn}

We denote $C_T(U,J^1({\mathbb R}\times {\mathbb R}^n,
{\mathbb R}))^{(m)}= C_T(U,J^1({\mathbb R}\times {\mathbb R}^n,
{\mathbb R}))\times\cdots\times  C_T(U,J^1({\mathbb R}\times {\mathbb R}^n,
{\mathbb R}))$  ($m$-products) and denote other notations analogously.
We call an element $C=(C_1,\ldots,C_m)\in C_T(J^1({\mathbb R}\times {\mathbb R}^n,
{\mathbb R}),0)^{(m)}$ {\it a ${\cal P}_{(m)}$-contact diffeomorphism germ}
if $\Pi\circ C(0)=0$ and $C_1(0),\ldots,C_m(0)$ are distinct.
\begin{tth}[Multi ${\cal P}$-Contact transversality theorem]\label{mp_contra:th}
Let $Q_i,$ $i=1,2,\ldots$ are submanifolds of 
$J^{l}_{C^\Theta_T}(U,J^1({\mathbb R}\times {\mathbb R}^n,
{\mathbb R}))^{(m)}$.
Then the set 
\begin{eqnarray*}
T=\{ C=(C_1,\ldots,C_m)\in C^\Theta_T(U,
J^1({\mathbb R}\times {\mathbb R}^n,{\mathbb R}))^{(m)}| \hspace{1cm}\\
j^lC \mbox{ is transversal to } Q_i \mbox{ for all } i\in {\mathbb N}\}
\end{eqnarray*}
is a residual set in 
$C^\Theta_T(U,
J^1({\mathbb R}\times {\mathbb R}^n,{\mathbb R}))^{(m)}$
\end{tth}
This is directly proved by ${\cal P}$-Contact transversality theorem 
(\cite[Theorem 0.4.5]{tPleg2:cite}).\\

We denote the ring ${\cal E}(1+n+n,m)$ 
on the coordinates $(t,q,p)$
by ${\cal E}_{t,q,p(m)}$ and denote other notations analogously.
\begin{tth}\label{staLeg:th}
Let ${\cal L}=({\cal L}_1,\ldots,{\cal L}_m)$ 
be a multi-reticular Legendrian unfolding with a 
multi-generating family $F=(F_1,\ldots,F_m)$.
Let $C=(C_1,\ldots,C_m)\in  C^\Theta_T(
J^1({\mathbb R}\times {\mathbb R}^n,{\mathbb R}),0)^{(m)}$ be 
an extension of ${\cal L}$.
Then the following are equivalent.\\
{\rm (u)} $F$ is a reticular $t$-$({\cal P}$-${\cal K})_{(m)}$-stable unfolding of 
$F|_{t=0}$.\\
{\rm (hs)} ${\cal L}$ is homotopically stable.\\
{\rm (is)} ${\cal L}$ is infinitesimally stable.\\ 
{\rm (a)} ${\cal E}_{t,q,p(m)}=
B^{r_1}_0\times\cdots\times B^{r_m}_0+
\langle 1,p_1\circ C',\ldots,p_n\circ C'\rangle_{(\Pi\circ C')^*{\cal E}_{t,q,z}}+
\langle s\circ C'\rangle_{{\cal E}_t}$,
where $C'=C|_{z=s=0}$, $B^{r_i}_0=\langle q_1p_1,\ldots,q_{r_i}p_{r_i},
q_{r_i+1},\ldots,q_n\rangle_{{\cal E}_{t,q,p}}$, and 
$(\Pi\circ C')^*{\cal E}_{t,q,z}$ be the ${\cal E}_{t,q,z}$-submodule of 
${\cal E}_{t,q,z(m)}$ such that $a.(f_1,$ $\ldots,f_m)=(a(\Pi\circ C'_1(f_1)),
\ldots,a(\Pi\circ C'_m(f_m))$ for 
$a\in {\cal E}_{t,q,z},(f_1,\ldots,f_m)\in{\cal E}_{t,q,p(m)}$.
\end{tth}
We remark that sufficiently near multi-reticular Legendrian unfoldings of 
stable one are stable by the condition (a).
\section{Genericity}\label{generic:sec}
\quad 
In order to give a generic classification of multi-reticular Legendrian 
unfoldings,
we reduce our investigation to a finitely dimensional jet space 
of ${\cal P}$-contact diffeomorphism germs.
\begin{dfn}\label{l,l+1detLeg}{\rm
Let ${\cal L}=({\cal L}_1,\ldots,{\cal L}_m)$ be a multi-reticular Legendrian unfolding.
We say that ${\cal L}$ is $l$-determined if the following condition holds:
For any extension $C_i\in C_T(J^1({\mathbb R}\times{\mathbb R}^n,{\mathbb R}),0)$ of ${\cal L}_i$, the multi-reticular Legendrian unfolding 
$(C'_1|_{{\mathbb L}^{r_1}},\ldots,C'_m|_{{\mathbb L}^{r_m}})$ 
and ${\cal L}$ are ${\cal P}_{(m)}$-Legendrian 
equivalent for all $C'_i\in C_T(J^1({\mathbb R}\times{\mathbb R}^n,{\mathbb R}),0)$ satisfying that  $j^lC_i(0)=j^lC'_i(0)$ for $i=1,\ldots,m$.
}\end{dfn}

As Lemma \ref{sta:lm}, we may consider the following other definition
of finitely determinacy of reticular Legendrian maps:\\
(1) The definition given by replacing $C_T(J^1({\mathbb R}\times{\mathbb R}^n,{\mathbb R}),0)$
to $C^\Theta_T(J^1({\mathbb R}$ $\times{\mathbb R}^n,{\mathbb R}),0)$.\\
(2) The definition given by replacing $C_T(J^1({\mathbb R}\times{\mathbb R}^n,{\mathbb R}),0)$
to $C^Z_T(J^1({\mathbb R}\times{\mathbb R}^n,{\mathbb R}),0)$.\\
(3) The definition given by replacing $C_T(J^1({\mathbb R}\times{\mathbb R}^n,{\mathbb R}),0)$
to $C^{\Theta,Z}_T(J^1({\mathbb R}$ $\times{\mathbb R}^n,{\mathbb R}),0)$.\\
Then the following holds by Proposition 0.5.2 in \cite{tPleg2:cite}:
\begin{prop}
Let ${\cal L}$ be a multi-reticular Legendrian unfolding.
Then \\
{\rm (A)} If ${\cal L}$ is $l$-determined of the original definition, then ${\cal L}$ is $l$-determined of the 
definition {\rm (1)}.\\
{\rm (B)} If ${\cal L}$ is $l$-determined of the definition {\rm (1)}, then ${\cal L}$ is $l$-determined of the definition {\rm (3)}.\\
{\rm (C)} If ${\cal L}$ is $(l+1)$-determined of the definition {\rm (3)}, then ${\cal L}$ is $l$-determined of the definition {\rm (2)}.\\
{\rm (D)} If ${\cal L}$ is $l$-determined of the definition {\rm (2)}, then 
${\cal L}$ is $l$-determined of the original definition.
\end{prop}

\begin{tth}\label{n+1det:Leg}
Let ${\cal L}=({\cal L}_1,\ldots,{\cal L}_m)$ be a multi
reticular Legendrian unfolding.
If ${\cal L}$ is infinitesimally stable then ${\cal L}$ is $(n+5)$-determined.
\end{tth}

Let  $J^l(2n+3,2n+3)^{(m)}$ be the set of multi-$l$-jets of map germs  from $(J^1({\mathbb R}\times{\mathbb R}^n,{\mathbb R}),0)$ 
to $J^1({\mathbb R}\times{\mathbb R}^n,{\mathbb R})$ and 
$tC^l(n,m)$ be the Lie group in $J^l(2n+3,2n+3)^{(m)}$ 
consists of multi-$l$-jets of ${\cal P}$-contact embedding  germs.
Let $L^l(2n+3)^{(m)}$ be the Lie group consists of multi-$l$-jets of diffeomorphism germs on 
$(J^1({\mathbb R}\times{\mathbb R}^n,{\mathbb R}),0)$.

We  consider the Lie subgroup $rtLe^l(n,m)$ of $L^l(2n+3)^{(m)}\times L^l(2n+3)$ consists of 
 multi-$l$-jets of reticular ${\cal P}_{r_i}$-diffeomorphisms on the  source space and $l$-jets of ${\cal P}$-Legendrian 
equivalences of $\Pi$ at $0$: 
\begin{eqnarray*}
\lefteqn{
 rtLe^l(n,m)}\\
& & =\{ (j^l\Psi_1(0),\ldots,j^l\Psi_m(0),j^lK(0))\in 
L^l(2n+3)^{(m)}\times L^l(2n+3)\ | \\
& & \Psi_i
 \mbox{ is a reticular} \mbox{ ${\cal P}_{r_i}$-diffeomorphism on }
(J^1({\mathbb R}\times{\mathbb R}^n,{\mathbb R}),0), \\
& & K \mbox{ is a ${\cal P}$-Legendrian equivalence of } 
\Pi \}.
\end{eqnarray*}

The group $rtLe^l(n,m)$ acts on $J^l(2n+3,2n+3)$ and $tC^l(n,m)$ 
is invariant 
under this action.
Let $C=(C_1,\ldots,C_m)$ be  a ${\cal P}_{(m)}$-contact diffeomorphism germ from 
$(J^1({\mathbb R}\times{\mathbb R}^n,{\mathbb R}),0)$ to 
$J^1({\mathbb R}\times{\mathbb R}^n,{\mathbb R})$
and set  
$z_i=j^lC_i(0)$, ${\cal L}_i=C_i|_{{\mathbb L}^{r_i}}$,
${\cal L}=({\cal L}_1,\ldots,{\cal L}_m)$, $z=(z_1,\ldots,z_m,)$.
We denote  the orbit $rtLe^l(n,m)\cdot z$ by $[z]$.
Then 
\begin{eqnarray*}
[z]=\{ j^lC'(0)\in tC^l(n,m)\ | \ {\cal L}  \mbox{ and } 
(C'_1|_{{\mathbb L}^{r_1}},\ldots,C'_m|_{{\mathbb L}^{r_m}})\hspace{1cm}\\ 
\mbox{ are ${\cal P}_{(m)}$-Legendrian equivalent} \}. 
\end{eqnarray*}

We denote by $VI_C$  the vector space consists of infinitesimal ${\cal P}_{(m)}$-contact transformation germs of $C$
and  denote by $VI^0_C$ the subspace of $VI_C$ consists of germs which
 vanish on $0$.
We denote by $VL_{J^1({\mathbb R}\times{\mathbb R}^n,{\mathbb R})}$ by the vector space consists of 
infinitesimal ${\cal P}$-Legendrian 
equivalences on $\Pi$  and denote by $VL^0_{J^1({\mathbb R}\times{\mathbb R}^n,{\mathbb R})}$
by the subspace of $VL_{J^1({\mathbb R}\times{\mathbb R}^n,{\mathbb R})}$ consists of germs which vanish at $0$.

We denote by $V^0_{{\mathbb L}^{r}}$ the vector space consists of infinitesimal reticular $r$-diffeomorphisms 
on $(J^1({\mathbb R}\times{\mathbb R}^n,{\mathbb R}),0)$ which vanishes at $0$ and set $V^0_{{\mathbb L}}=V^0_{{\mathbb L}^{r_1}}\times V^0_{{\mathbb L}^{r_m}}$.
We have that by Lemma \ref{infsta:t-Leglem}:
\begin{eqnarray*}
VI^0_C=\{ (v_1,\ldots,v_m)|v_i:(J^1({\mathbb R}\times{\mathbb R}^n,{\mathbb R}),0)\rightarrow T(J^1({\mathbb R}\times{\mathbb R}^n,{\mathbb R})),\hspace{0.5cm}\\
v_i=X_{f_i}\circ C_i \mbox{ for some }f\in {\mathfrak M}^2_{t,q,z,p}\},
\end{eqnarray*}
\begin{eqnarray*}
VL^0_{J^1({\mathbb R}\times{\mathbb R}^n,{\mathbb R})}=\{\eta\in X(J^1({\mathbb R}\times{\mathbb R}^n,{\mathbb R}),0)\ | \
\eta=X_H
\mbox{ for some} \hspace{1cm}\\
\mbox{ ${\cal P}$-fiver preserving  function germ } H\in {\mathfrak M}^2_{J^1({\mathbb R}\times{\mathbb R}^n,{\mathbb R})}\}, 
\end{eqnarray*}
\begin{eqnarray*}
V^0_{{\mathbb L}}=\{ (\xi_1,\ldots,\xi_m)| \xi_i\in X(J^1({\mathbb R}\times{\mathbb R}^n,{\mathbb R}),0),\ 
\xi_i=X_{g_i}  \mbox{ for some }\\g_i \in B^{r_i}\cap {\mathfrak M}^2_{J^1({\mathbb R}\times{\mathbb R}^n,{\mathbb R})}\}. 
\end{eqnarray*}

\vspace{2mm}
We define the homomorphism $tC:{\mathfrak M}_{J^1({\mathbb R}\times{\mathbb R}^n,{\mathbb R})}
VI_{{\mathbb L}}\rightarrow VI^0_{C}$ by 
$tC(v)=(C_{1*}v_1,$ $\ldots,C_{m*}v_m)$ and define the homomorphism 
$wC:VL^0_{J^1({\mathbb R}\times{\mathbb R}^n,{\mathbb R})}\rightarrow 
VI^0_{C}$ by 
$wC(\eta)=(\eta\circ C_1,\ldots,\eta\circ C_m)$.

\vspace{2mm}
We denote $VI^{l}_C$ the subspace of $VI_C$ consists of infinitesimal
${\cal P}$-contact transformation germs of $C$ whose $l$-jets are $0$:
\[  VI^l_C=\{ (v_1,\ldots,v_m)\in VI_C\ | \ j^lv_i(0)=0\}. \]

For $\tilde{C}=(\tilde{C}_1,\ldots,\tilde{C}_m)\in C_T(U,J^1({\mathbb R}\times{\mathbb R}^n,{\mathbb R}))^{(m)}$,
we define
 the continuous map 
$j^l_0\tilde{C}:U\rightarrow tC^l(n,m)$ by $x$ to the $l$-jet of 
$(\tilde{C}_{1x},\ldots,\tilde{C}_{mx})$.
\begin{tth}\label{stabletrans_tleg2:th}
Let ${\cal L}=({\cal L}_1,\ldots,{\cal L}_m)$ be a reticular Legendrian unfolding. 
Let $C_i$ be an extension of ${\cal L}_i$ and $l\geq (n+2)^2$.
We set $C=(C_1,\ldots,C_m)$.
Then the followings are equivalent:\\
{\rm (s)} ${\cal L}$  is stable.\\
{\rm (t)} $j^l_0C$ is transversal to $[j^l_0C(0)]$.\\
{\rm (a')} ${\cal E}_{t,q,p(m)}=
B^{r_1}_0\times\cdots\times B^{r_m}_0+
\langle 1,p_1\circ C',\ldots,p_n\circ C'\rangle_{(\Pi\circ C')^*{\cal E}_{t,q,z}}+
\langle s\circ C'\rangle_{{\cal E}_t}+{\mathfrak M}_{t,q,p(m)}^l$,
where $C'=C|_{z=s=0}$ and $B^{r_i}_0=\langle q_1p_1,\ldots,q_{r_i}p_{r_i},
q_{r_i+1},$ $\ldots,q_n\rangle_{{\cal E}_{t,q,p}}$,\\
{\rm (a)} ${\cal E}_{t,q,p(m)}=
B^{r_1}_0\times\cdots\times B^{r_m}_0+
\langle 1,p_1\circ C',\ldots,p_n\circ C'\rangle_{(\Pi\circ C')^*{\cal E}_{t,q,z}}+
\langle s\circ C'\rangle_{{\cal E}_t}$,\\
{\rm (is)} ${\cal L}$ is infinitesimally stable,\\
{\rm (hs)} ${\cal L}$ is homotopically stable,\\
{\rm (u)} A multi-generating family $F$ of ${\cal L}$ is reticular 
$t$-$({\cal P}$-${\cal K})_{(m)}$-stable unfolding of $F|_{t=0}$.
\end{tth}

Let ${\cal L}=({\cal L}_1,\ldots,{\cal L}_m)$ be a stable multi-reticular Legendrian unfolding.
We say that ${\cal L}$ is {\it simple} if 
there exists a representative $\tilde{C}\in C_T(
U,J^1({\mathbb R}\times{\mathbb R}^n,{\mathbb R}))^{(m)}$ of
 a extension of ${\cal L}$ such that 
$\{ \tilde{C}_x| x\in  U\}$ is covered by 
finite orbits $[C_1],\ldots,[C_m]$ for  
 some multi-${\cal P}$-contact diffeomorphism germs $C_1,\ldots,C_m\in 
C_T(U,J^1({\mathbb R}\times{\mathbb R}^n,{\mathbb R}))^{(m)}$. 

\begin{prop}
Let ${\cal L}=({\cal L}_1,\ldots,{\cal L}_m)$ be a stable multi-reticular Legendrian unfolding.
Then ${\cal L}$ is simple if and only if all reticular Legendrian unfoldings
${\cal L}_i$ are simple for $i=1,\ldots,m$.
\end{prop}

In order to classify generic multi-reticular Legendrian unfoldings,
we classify stable multi-unfoldings $F=(F_1,\ldots,F_m)\ (n\leq 2,m\geq 2)$ of 
$f=(f_{1},\ldots,f_{m})$ and $f_0=(f_{0,1},\ldots,f_{0,m})$ 
satisfying the condition: the reticular $({\cal P}$-${\cal K})_{(m-1)}$-codimension of 
$(f_1,\ldots,\check{f_i},\ldots,,f_m)\ =0$ for any $i$.\\

Let a stable multi-unfoldings $F=(F_1,\ldots,F_m)\ (n\leq 2,m\geq 2)$ of 
$f=(f_{1},\ldots,f_{m})$ and $f_0=(f_{0,1},\ldots,f_{0,m})$ 
satisfying the condition be given.
Then each $F_i$ is a reticular $t$-${\cal P}$-${\cal K}$-stable unfolding of 
$f_i$, there exist monomials $\varphi_{i,1},\ldots,\varphi_{i,\mu_i}
\in {\mathfrak M}(r_i;k_i)$ such that they consist a basis of 
$Q_{f_{0,i}}=\langle f_{0,i},x\frac{\partial f_{0,i}}{\partial x},
\frac{\partial f_{0,i}}{\partial y}\rangle_{{\cal E}(r_i;k_i)}$ and 
$\varphi_{i,1}$ has the maximal degree, and $\varphi_{i,\mu_i}=1$.
Then we have that $\mu_1+\cdots \mu_m\leq n+2$.
Since $\sum_{i=1}^m\mu_i\leq n+2\leq 4$, we have that all $f_{0,i}$ are
simple singularities.
Therefore $f_0$ is stably 
reticular ${\cal K}_{(m)}$-equivalent to one of the multi-germs in the following 
list:\\
$n=1$;\\
\ \ \ $m=1;\ y^2,\ y^3,\ y^4,\ x^2,\ x^3,\ \pm xy+y^3$,\\
\ \ \ $m=2;\ (y^2,y^2),\ (y^2,y^3),\ (y^2,x^2)$,\\
\ \ \ $m=3;\ (y^2,y^2,y^2)$,\\
$n=2$;\\
\ \ \ $m=1;\ y^2,\ y^3,\ y^4,\ y^5,\ y_1^2\pm y_2^2,\ 
x^2,\ x^3,\ x^4,\ \pm xy+y^3,\ xy+y^4,\ x^2+y^3$,\\
\ \ \ $m=2; (y^2,y^2),\ (y^2,y^3),\ (y^2,y^4),\ (y^3,y^3),\ (y^2,x^2),\ 
(y^2,x^3),\ (y^2,\pm xy+y^3),\ (x^2,x^2)$,\\
\ \ \ $m=3;\ (y^2,y^2,y^2),\ (y^2,y^2,y^3),\ (y^2,y^2,x^2)$,\\
\ \ \ $m=4;\ (y^2,y^2,y^2,y^2)$.\\

We construct a reticular $({\cal P}$-${\cal K})_{(m)}$-versal unfolding for each multi-germ by the usual method.
Then the corresponding list is as follows:\\
$n=1;\\
\ \ \ \ \ (1)\ \ (y^2+u_{1,1},y^2+u_{2,1})\\
\ \ \ \ \ (2)\ \ (y^2+u_{1,1},y^3+u_{2,1}y+u_{2,2})\\
\ \ \ \ \ (3)\ \ (y^2+u_{1,1},x^2+u_{2,1}x+u_{2,2})\\
\ \ \ \ \ (4)\ \ (y^2+u_{1,1},y^2+u_{2,1},y^2+u_{3,1})$\\
$n=2;\\
\ \ \ \ \ (6)\ \ (y^2+u_{1,1},y^2+u_{2,1})\\
\ \ \ \ \ (7)\ \ (y^2+u_{1,1},y^3+u_{2,1}y+u_{2,2})\\
\ \ \ \ \ (8)\ \ (y^2+u_{1,1},y^4+u_{2,1}y^2+u_{2,2}y+u_{2,3})\\
\ \ \ \ \ (9)\ \ (y^3+u_{1,1}y+u_{1,2},y^3+u_{2,1}y+u_{2,2})\\
\ \ \ \  (10)\ \  (y^2+u_{1,1},x^2+u_{2,1}x+u_{2,2})\\
\ \ \ \ (11)\ \ (y^2+u_{1,1},x^3+u_{2,1}x^2+u_{2,2}x+u_{2,3})\\
\ \ \ \ (12)\ \ (y^2+u_{1,1},\pm xy+y^3+u_{2,1}y^2+u_{2,2}y+u_{2,3})\\
\ \ \ \ (13)\ \ (x^2+u_{1,1}x+u_{1,2},x^2+u_{2,1}x+u_{2,2})\\
\ \ \ \ (14)\ \ (y^2+u_{1,1},y^2+u_{2,1},y^2+u_{3,1})\\
\ \ \ \ (15)\ \ (y^2+u_{1,1},y^2+u_{2,1},y^3+u_{3,1}y+u_{3,2})\\
\ \ \ \ (16)\ \ (y^2+u_{1,1},y^2+u_{2,1},x^2+u_{3,1}x+u_{3,2})\\
\ \ \ \ (17)\ \ (y^2+u_{1,1},y^2+u_{2,1},y^2+u_{3,1},y^2+u_{4,1})$.\\

In the case that the reticular $({\cal P}$-${\cal K})_{(m)}$-codimension of 
$f$ $=0$, the $F$ is stably reticular $t$-$({\cal P}$-${\cal K})_{(m)}$-equivalent 
to $G=(G_1,\ldots,G_m)$, 
where $G_i(x,y,t,q,z)=f_{0,i}(x,y)+u_{i,1}\varphi_1(x,y)+\cdots +
u_{i,\mu_i}\varphi_{i,\mu_i}(x,y)-z$ for $i=1,\ldots,m-1$,
$G_m(x,y,t,q,z)=f_{0,i}(x,y)+u_{m,1}\varphi_1(x,y)+\cdots +
u_{m,\mu_{m-1}}\varphi_{m,\mu_{m-1}}(x,y)-z$, and  
$(q_1,\ldots,q_n,z)=(u_{1,1},\ldots,u_{1,\mu_1},\ldots,u_{m,1},$ $\ldots,
u_{m,\mu_{m-1}},u_1,\ldots,u_\mu,z)$.

In the case that the reticular $({\cal P}$-${\cal K})_{(m)}$-codimension of 
$f$ $=1$, $F$ is stably reticular $t$-$({\cal P}$-${\cal K})_{(m)}$-equivalent to 
$(F'_1,\ldots,F'_m)$, where\\
$(1)$ $F'_1=f_{0,1}(x,y)+(t+a(q,z))\varphi_1(x,y)+u_{1,1}\varphi_2(x,y)\cdots +
u_{1,\mu_{1}-1}\varphi_{1,\mu_{1}}(x,$ $y)-z$,\\
$(2)$ $F'_i=f_{0,i}(x,y)+u_{i,1}\varphi_1(x,y)+\cdots +
u_{i,\mu_{i}}\varphi_{i,\mu_{i}}(x,y)-z$ for $1<i<m$,\\
$(3)$  $F'_m=f_{0,m}(x,y)+u_{m,1}\varphi_{m,1}(x,y)+\cdots +
u_{m,\mu_{m}-1}\varphi_{m,\mu_{m}-1}(x,y)-z$,\\
$(4)$ $(q_1,\ldots,q_n,z)=(u_{1,1},\ldots,u_{1,\mu_1-1},
u_{2,1},\ldots,u_{2,\mu_2},\ldots,u_{m,1},\ldots,
u_{m,\mu_{m-1}},$ $u_1,\ldots,u_\mu,z)$,\\
$(5)$ $\frac{\partial a}{\partial  u_{1,j}}(0)=0$ 
for $j=1,\ldots,\mu_1-1$, $\frac{\partial a}{\partial  z}(0)=0$,
and $(\frac{\partial^2 a}{\partial  u_i\partial u_j}(0))_{i,j=1,\ldots,\mu}$ is 
non-degenerate.

We denote the linear part of $a$ by $v=v_2+\cdots +v_{m}$,
where $v_i$ depends only on $u_{i,1},\ldots,u_{i,\mu_i}$.
Then we may reduce $a$ to 
\begin{eqnarray*}
a=t\pm u_{2,1}\pm \cdots \pm u_{m,1}\pm u_1^2\pm \cdots u_\mu^2.
\end{eqnarray*}

Then $F$ is stably reticular $t$-$({\cal P}$-${\cal K})_{(m)}$-equivalent to one of the following list:\\
$n=1;\vspace{1mm}\\
\ \ \ \ m=2;\vspace{1mm}\\
\ \ \ \ \ \ {}^0({}^0A_1{}^0A_1):\ \ (y^2+q-z,y^2-z);\\  
\ \ \ \ \ \ {}^1({}^0A_1{}^0A_1):\ \ 
    (y^2+t\pm q^2-z,y^2-z);\\ 
\ \ \ \ \ \ {}^1({}^0A_1{}^0A_2): \ (y^2+t\pm q-z,y^3+qy-z);\\   
\ \ \ \ \ {}^1({}^0A_1{}^0B_2):\ \ (y^2+t\pm q-z,x^2+qx-z).
\vspace{1mm}\\  
\ \ \  m=3;\vspace{1mm}\\
\ \ \ \ \ {}^1({}^0A_1{}^0A_1{}^0A_1):\ \
 (y^2+t-z,y^2+q-z,y^2-q-z)$.\vspace{1mm}\\  
$n=2;\vspace{1mm}\\
\ \ \  m=2;\vspace{1mm}\\
\ \ \ \ \ {}^1({}^0A_1{}^0A_1):\ \ (y^2+t\pm q_1^2\pm q_2^2-z,y^2-z);\\
\ \ \ \ \ \ {}^0({}^0A_1{}^0A_2):\ \ 
(y^2+t\pm q_1-z,y^3+q_1y-q_2-z);\\  
\ \ \ \ \ {}^1({}^0A_1{}^0A_2):\ \ (y^2+t\pm q_1\pm q_2^2-z,y^3+q_1y-z);\\
\ \ \ \ \ \ {}^0({}^0A_1{}^0B_2): \ \  (y^2+q_1-z,x^2+q_2x-z);\\  
\ \ \ \ \ {}^1({}^0A_1{}^0A_3):\ \ (y^2+t\pm q_1-z,y^4+q_1y^2+q_2y-z);\\
\ \ \ \ \ {}^1({}^0A_2{}^0A_2):\  \ (y^3+(t\pm q_1)y+q_2-z,y^3+q_1y-z);\\
\ \ \ \ \ {}^1({}^0A_1{}^0B_2):\ \  (y^2+t\pm q_1\pm q_2^2-z,x^2+q_1x-z);\\
\ \ \ \ \ {}^1({}^0A_1B_3):\ \ (y^2+t\pm q_1-z,x^3+q_1x^2+q_2x-z);\\
\ \ \ \ \ {}^1({}^0A_1{}^0C^\pm_3):\ \ 
(y^2+t+q_1-z,\pm xy+y^3+q_1y^2+q_2y-z);\\
\ \ \ \ \ {}^1({}^0B_2{}^0B_2):\ \ (x^2+(t\pm q_1)x+q_2-z,x^2+q_1x-z).\vspace{1mm}\\
\ \ \  m=3;\vspace{1mm}\\
\ \ \ \ \ {}^0({}^0A_1{}^0A_1{}^0A_1):\ \ (y^2+q_1-z,y^2+q_2-z,y^2-z);\\
\ \ \ \ \ {}^1({}^0A_1{}^0A_1{}^0A_1):\ \ (y^2+t\pm q_1\pm q_2^2-z,
y^2+q_1-z,y^2-z);\\
\ \ \ \ \ {}^1({}^0A_1{}^0A_1{}^0A_2):\ \ (y^2+t\pm q_1-z,
y^2-z,y^3+q_1y+q_2-z);\\
\ \ \ \ \ {}^1({}^0A_1{}^0A_1{}^0B_2):\ \ (y^2+t\pm q_1-z,
y^2-z,x^2+q_1x+q_2-z).\vspace{1mm}\\
\ \ \  m=4;\vspace{1mm}\\
\ \ \ \ \ {}^1({}^0A_1{}^0A_1{}^0A_1{}^0A_1):\ \ 
(y^2+t\pm q_1\pm q_2-z,y^2+q_1-z,y^2+q_2-z,y^2-z)$.\\

\begin{tth}\label{genericclass}
Let $r_i=0$ or $1$ for $i=1,\ldots,m$ and $n\leq 2$.
Let  $U$ be a neighborhood of $0$ in 
$J^1({\mathbb R}\times{\mathbb R}^n,{\mathbb R})$.
Then there exists a residual set $O\subset  C^\Theta_T(U,J^1({\mathbb R}\times{\mathbb R}^n,{\mathbb R}))^{(m)}$ 
such that for any $\tilde{C}=(\tilde{C_1},\ldots,\tilde{C_m})\in O$ and 
$x=(x_1,\ldots,x_m)\in U^{(m)}$,
the multi-reticular Legendrian unfolding $
(\tilde{C_1}_{x_1}|_{{\mathbb L}^{r_1}},\ldots,
\tilde{C_m}_{x_m}|_{{\mathbb L}^{r_m}})$ 
is stable and have a multi-generating family which is stably reticular 
$t$-$({\cal P}$-${\cal K})_{(m)}$-equivalent for one of the types in the 
above list.
\end{tth}
{\it Proof}.
Let $X_i=(X_{i,1},\cdots,X_{i,m})\ (i=1,\ldots,s)$ be all simple singularities with 
$\sum_{j=1}^m r$-${\cal K}$-codim $X_{i,j}\leq n+2$,
that is each $X_{i,j}$ is one of simple singularities $A,B,C$.
For $i=1,\ldots,s$ let  $F_{X_i}$ be ${\cal P}_{(m)}$-$C$-non-degenerate map germ which is unfolding of $X_i$.
We choose an extension $C_{{}^0X_i}$ and $C_{{}^1X_i}$ of multi-reticular Legendrian unfoldings with 
multi-generating families $F_{{}^0X_i}$ and $F_{{}^1X_i}$ respectively.
Let $l>16$. We define that 
\begin{eqnarray*}
 O'=\{ \tilde{C}\in  
C^\Theta_T(U,J^1({\mathbb R}\times{\mathbb R}^n,{\mathbb R}))^{(m)}\ | 
j^{l}_0 \tilde{C} \mbox{ is transversal to} \hspace{1cm}\\  
\ [(j^{l}C_{{}^jX_i}(0)]
\mbox{ for all } i=1,\ldots,s \mbox{ and } j=0,1\}.
\end{eqnarray*}
Let $X'_i=(X_{i,1},\ldots,X_{i,m-1})\ (i=1,\ldots,s')$ be all simple singularities with 
$\sum_{j=1}^{m-1} r$-${\cal K}$-codim $X_{i,j}\leq n+1$.
We choose $C_{{}^0X'_i}\in
 C^\Theta_T(U,J^1({\mathbb R}\times{\mathbb R}^n,{\mathbb R}))^{(m-1)}$
analogous way.
Then we define that
\begin{eqnarray*}
 O''=\{ \tilde{C}\in  
C^\Theta_T(U,J^1({\mathbb R}\times{\mathbb R}^n,{\mathbb R}))^{(m)}\ | 
(j^{l}_0 \tilde{C_1},\ldots,\check{j^{l}_0 \tilde{C_j}},\ldots,j^{l}_0 \tilde{C_m}) \mbox{ is}\\
\mbox{transversal to }  [(j^{l}C_{{}^0X'_{i}}(0)]
\mbox{ for all } i=1\cdots,s' \}.
\end{eqnarray*}
Then $O'$ and $O''$ are residual sets. We set 
\[ Y=\{ j^{l}C(0)\in tC^{l}(n,m)\ |\ 
\mbox{the codimension of }[j^{l}C(0)]>2n+4\}.\]
Then $Y$ is an algebraic set in $tC^{l}(n,m)$.
Therefore we can define that 
\[ O'''=\{ \tilde{C}\in  C^\Theta_T(U,J^1({\mathbb R}\times{\mathbb R}^n,{\mathbb R}))^{(m)}\ |\ 
j^{l}_0\tilde{C}\mbox{ is transversal to } Y\}.\]
Then $Y$ has codimension $>2n+4$ because all 
${\cal P}_{(m)}$-contact diffeomorphism germ with 
$j^{l}C(0)\in Y$ adjoin to the above list which are simple.
Therefore 
\[ O'''=\{ \tilde{C}\in  C^\Theta_T(U,J^1({\mathbb R}\times{\mathbb R}^n,{\mathbb R}))^{(m)}\ |\ 
j^{l}_0\tilde{C}(U^{(m)})\cap Y=\emptyset \}.\]
Then the set $O=O'\cap O''\cap O'''$ has the required condition.\hfill q.e.d.
\begin{figure}[htbp]
 \begin{minipage}{0.48\hsize}
  \begin{center}
    \includegraphics*[width=5cm,height=5cm]{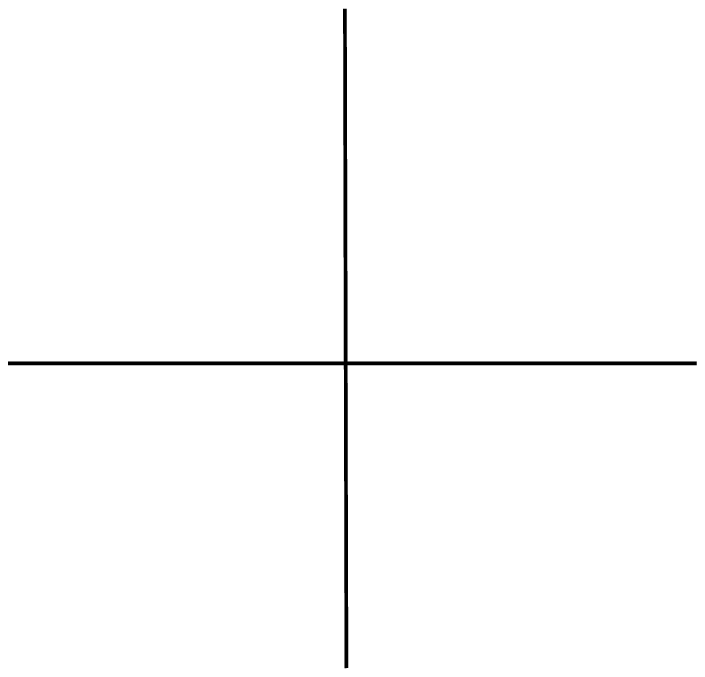}
  \end{center}
\caption{${}^0({}^0A_1{}^0A_1)$}
 \end{minipage}
 \begin{minipage}{0.48\hsize}
   \begin{center}
     \includegraphics*[width=5cm,height=5cm]{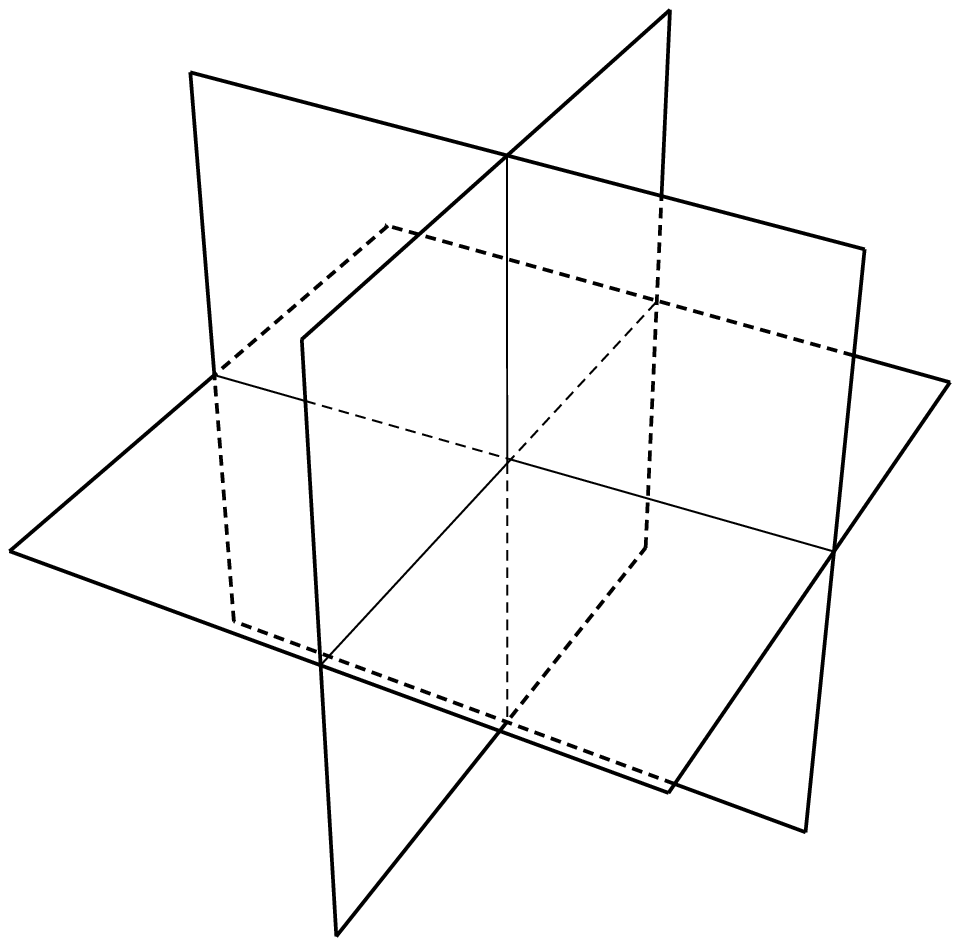}
   \end{center}
\caption{${}^0({}^0A_1{}^0A_1{}^0A_1)$}
 \end{minipage}
\end{figure}
\begin{figure}[htbp]
 \begin{minipage}{0.48\hsize}
  \begin{center}
    \includegraphics*[width=5cm,height=5cm]{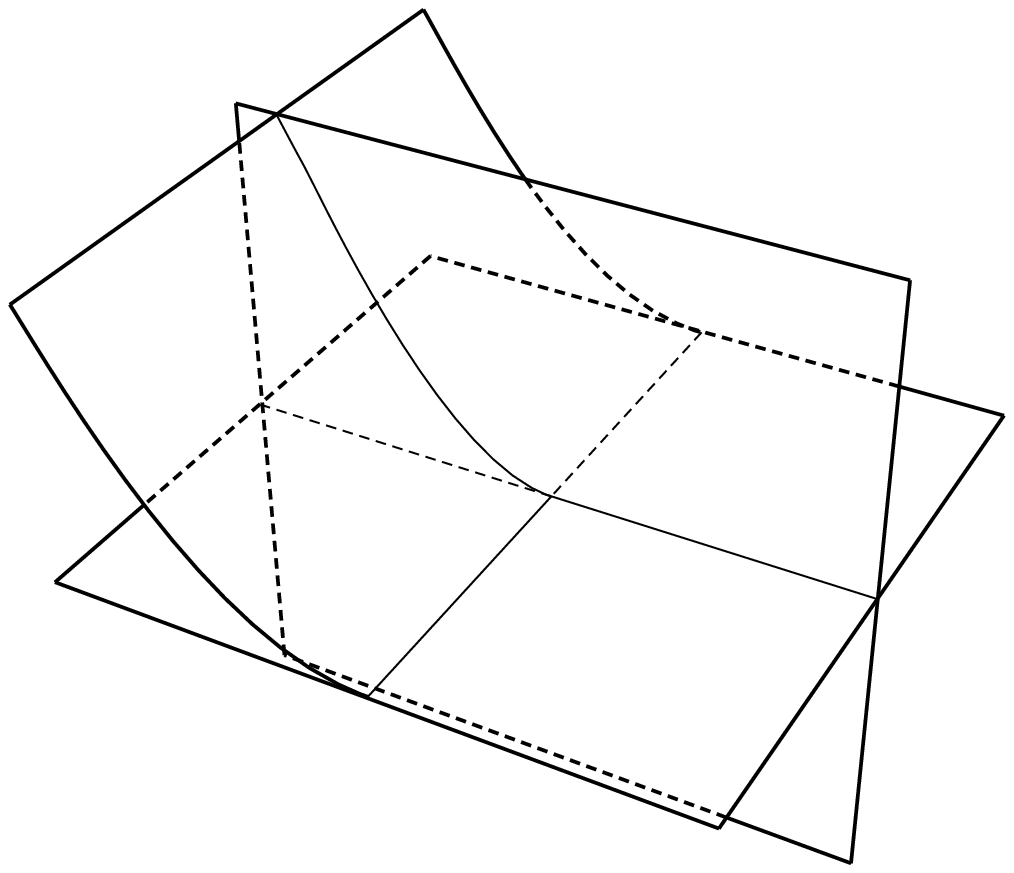}
  \end{center}
\caption{${}^0({}^0A_1{}^0B_2)$}
 \end{minipage}
 \begin{minipage}{0.48\hsize}
   \begin{center}
     \includegraphics*[width=5cm,height=5cm]{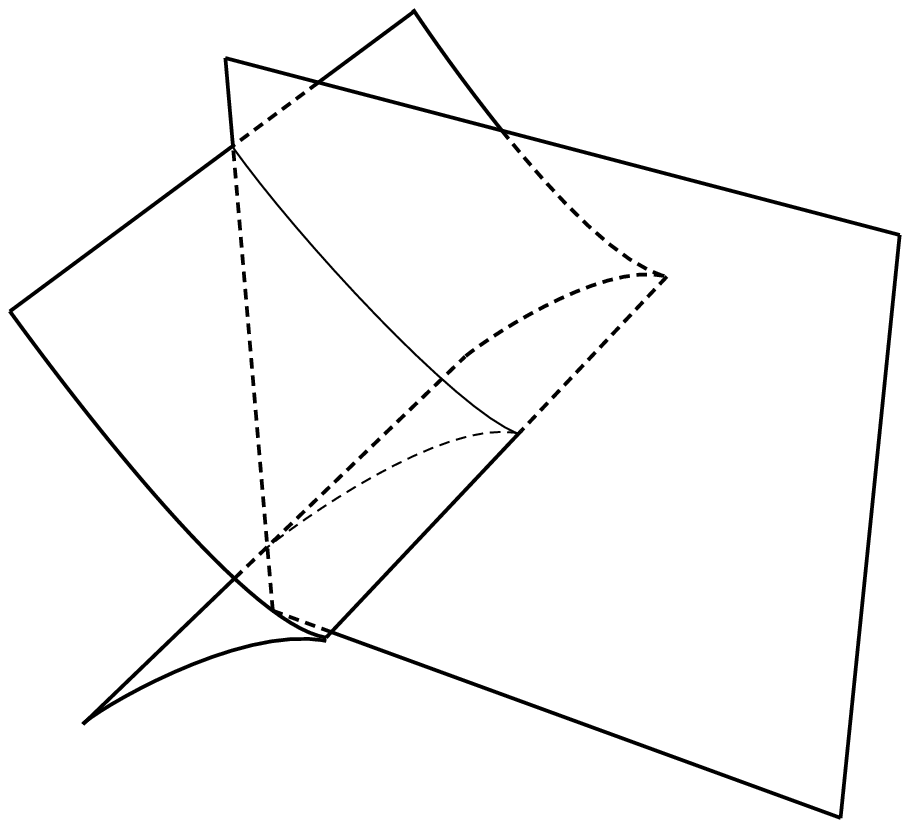}
   \end{center}
\caption{${}^0({}^0A_1{}^0A_2)$}
 \end{minipage}
\end{figure}
\begin{figure}[htbp]
 \begin{minipage}{0.30\hsize} 
  \begin{center}
    \includegraphics*[width=3cm,height=3cm]{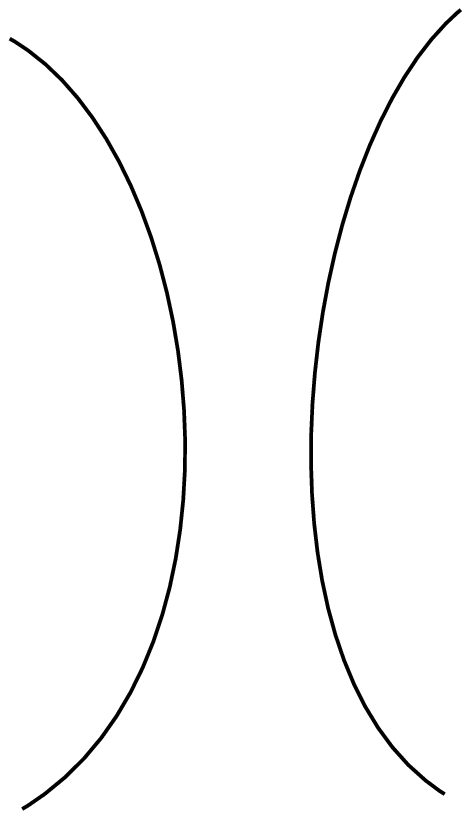}
  \end{center}
 \end{minipage}
$\leftrightarrow $
 \begin{minipage}{0.30\hsize}
   \begin{center}
     \includegraphics*[width=3cm,height=3cm]{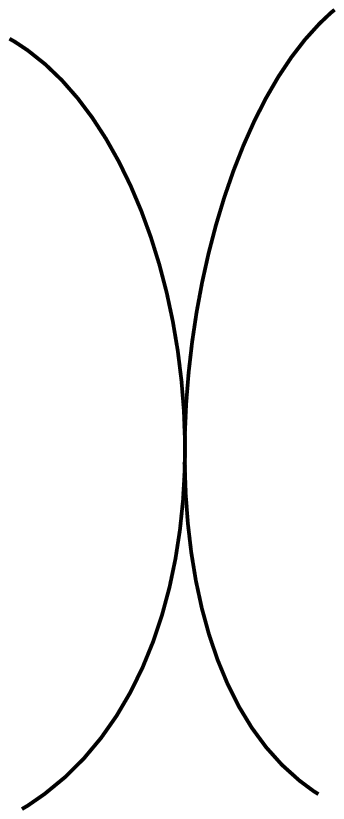}
   \end{center}
 \end{minipage}
$\leftrightarrow $
 \begin{minipage}{0.30\hsize}
  \begin{center}
  \includegraphics*[width=3cm,height=3cm]{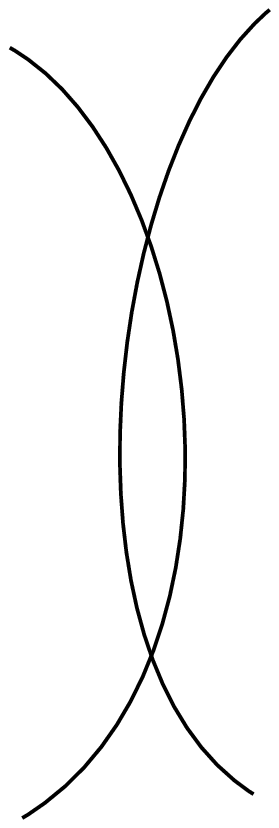}
  \end{center}
\end{minipage}\\
 \begin{minipage}{0.30\hsize} 
  \begin{center}
    \includegraphics*[width=3cm,height=3cm]{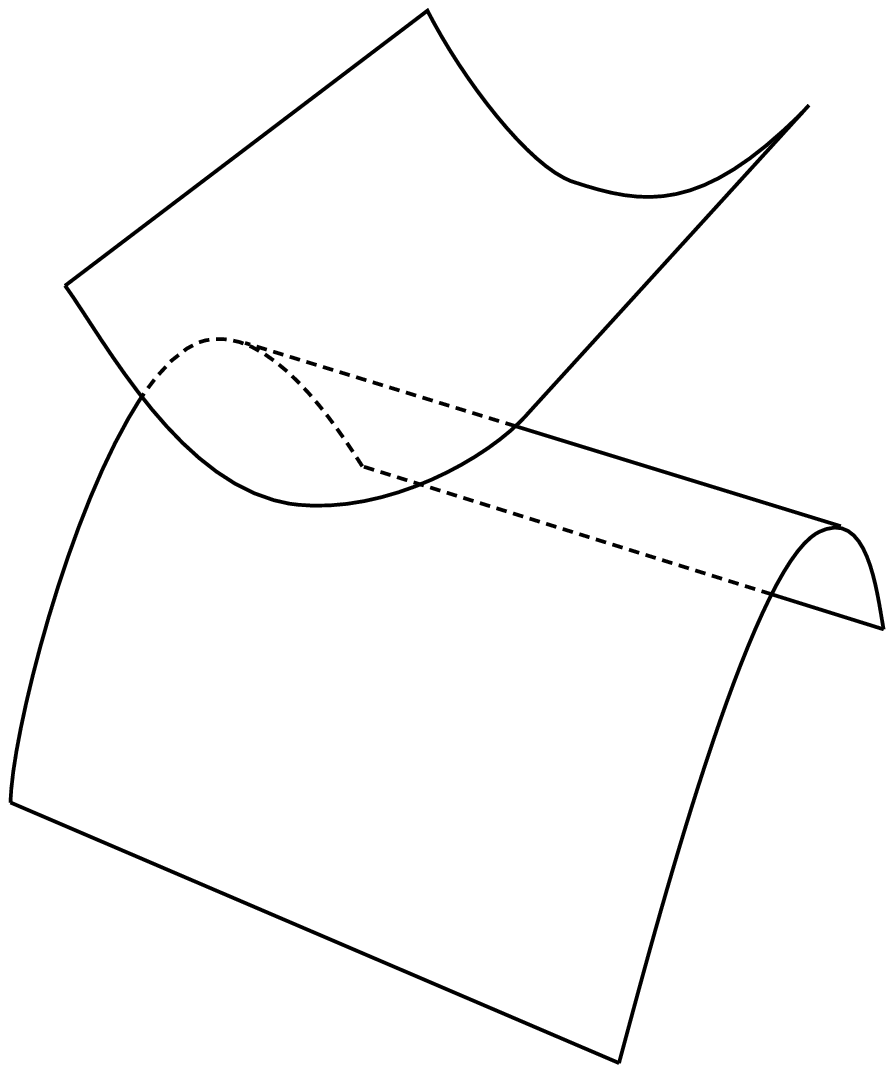}
  \end{center}
 \end{minipage}
$\leftrightarrow $
 \begin{minipage}{0.30\hsize}
   \begin{center}
     \includegraphics*[width=3cm,height=3cm]{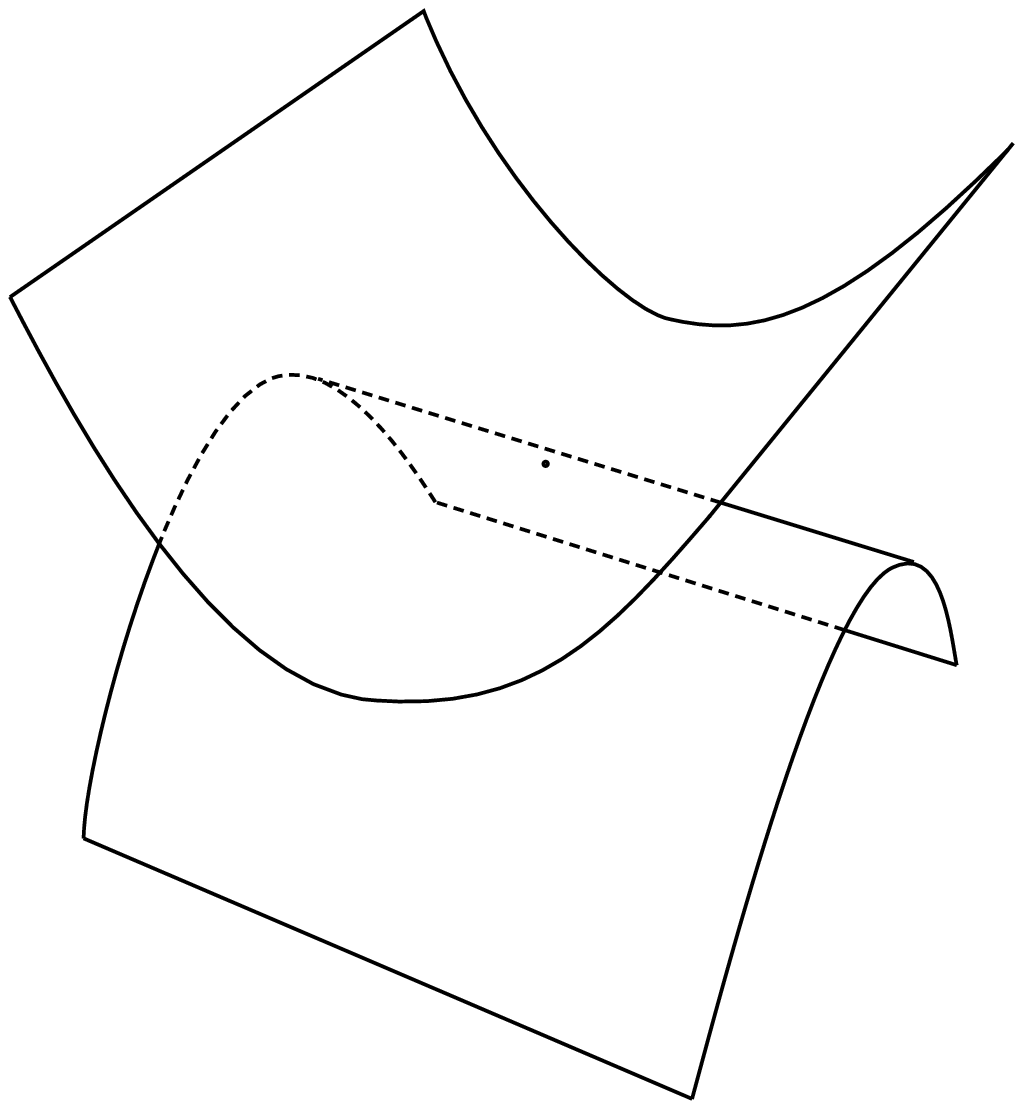}
   \end{center}
 \end{minipage}
$\leftrightarrow $
 \begin{minipage}{0.30\hsize}
  \begin{center}
  \includegraphics*[width=3cm,height=3cm]{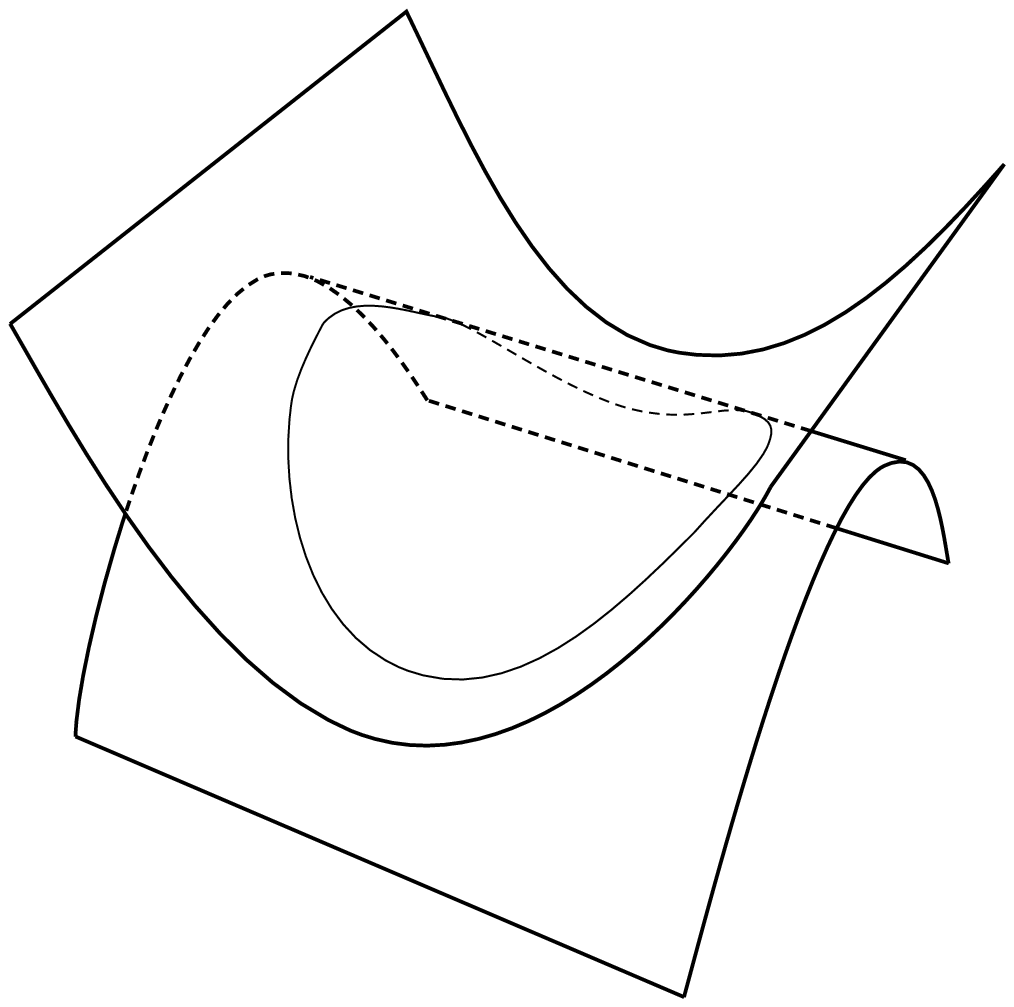}
  \end{center}
\end{minipage} \\
\begin{minipage}{0.30\hsize} 
  \begin{center}
    \includegraphics*[width=3cm,height=3cm]{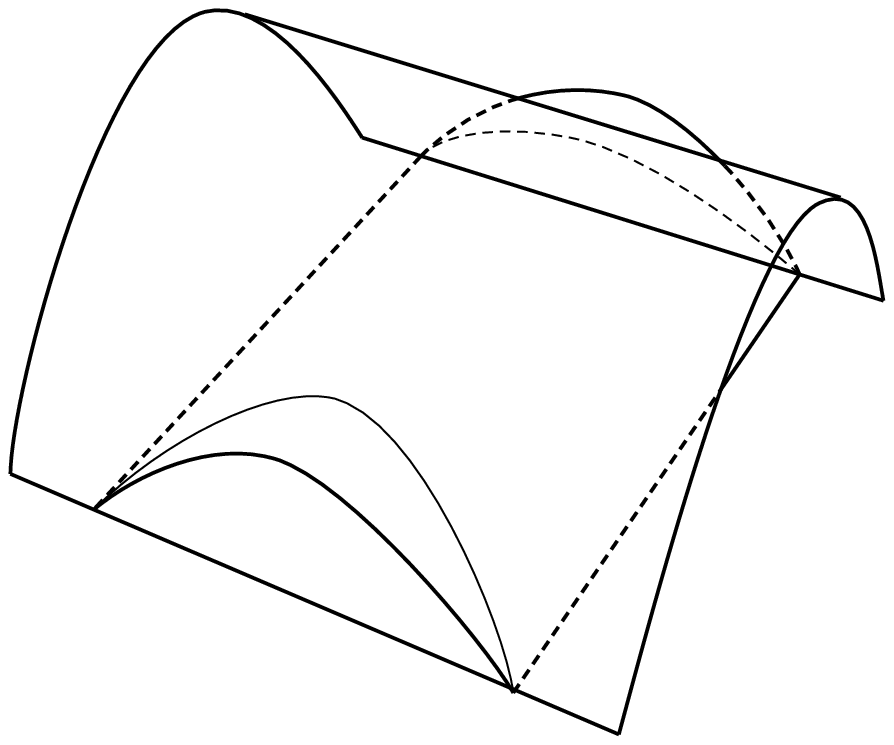}
  \end{center}
 \end{minipage}
$\leftrightarrow $
 \begin{minipage}{0.30\hsize}
   \begin{center}
     \includegraphics*[width=3cm,height=3cm]{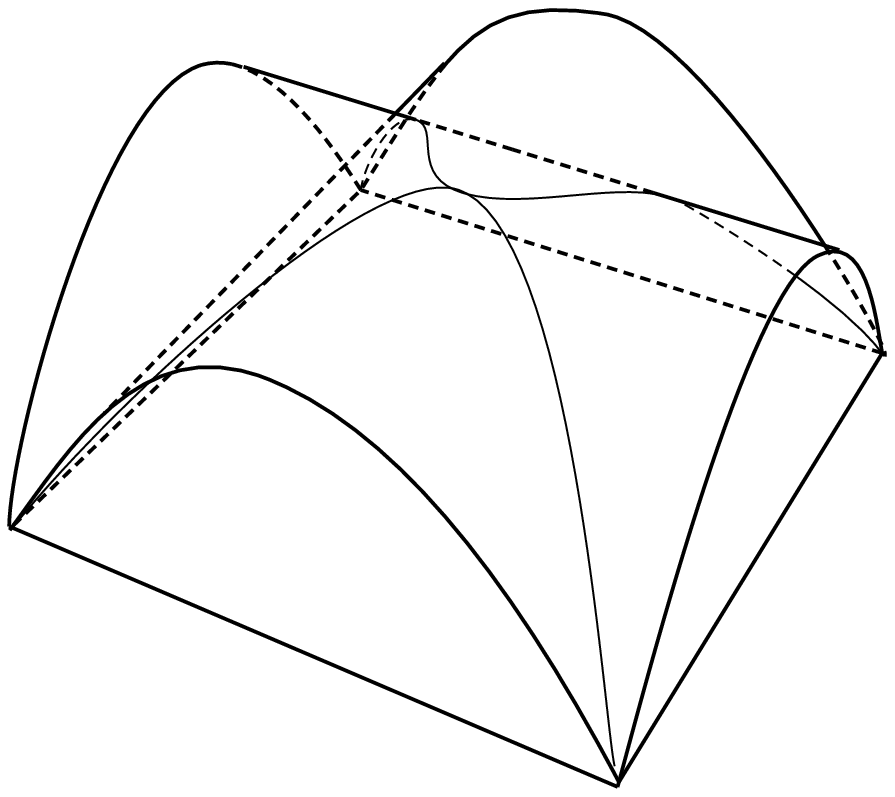}
   \end{center}
 \end{minipage}
$\leftrightarrow $
 \begin{minipage}{0.30\hsize}
  \begin{center}
\includegraphics*[width=3cm,height=3cm]{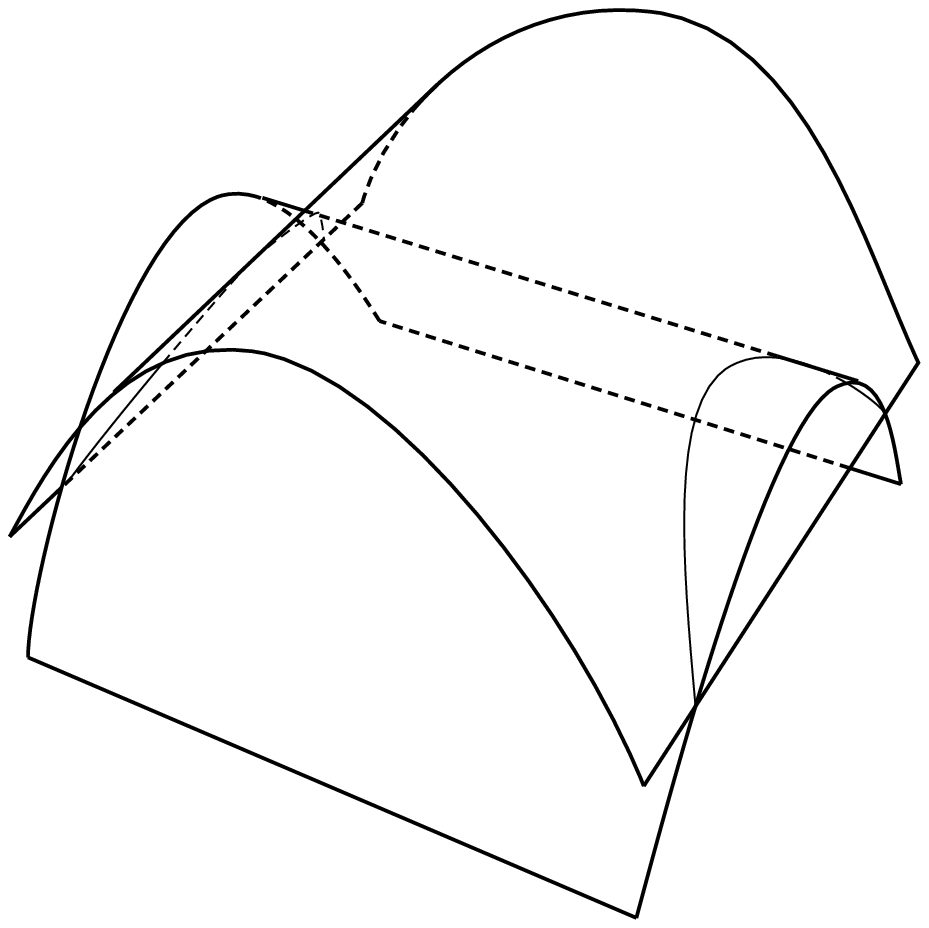}
  \end{center}
\end{minipage}
\caption{${}^1({}^0A_1{}^0A_1)$}
\end{figure}
\begin{figure}[htbp]
 \begin{minipage}{0.30\hsize} 
  \begin{center}
    \includegraphics*[width=3cm,height=3cm]{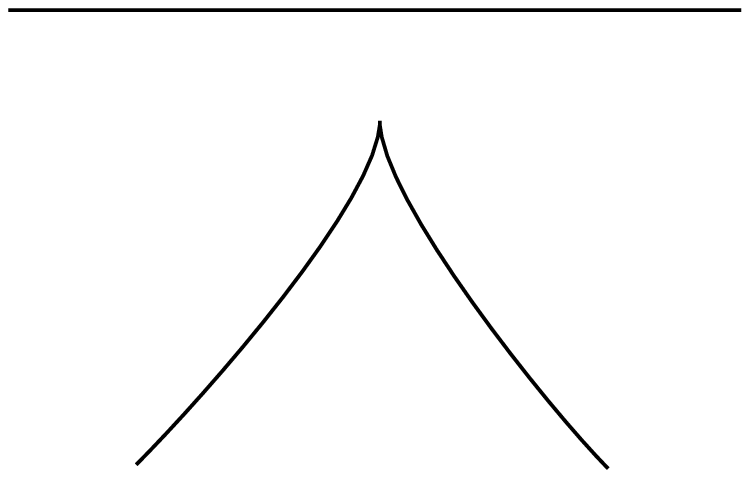}
  \end{center}
 \end{minipage}
$\leftrightarrow $
 \begin{minipage}{0.30\hsize}
   \begin{center}
     \includegraphics*[width=3cm,height=3cm]{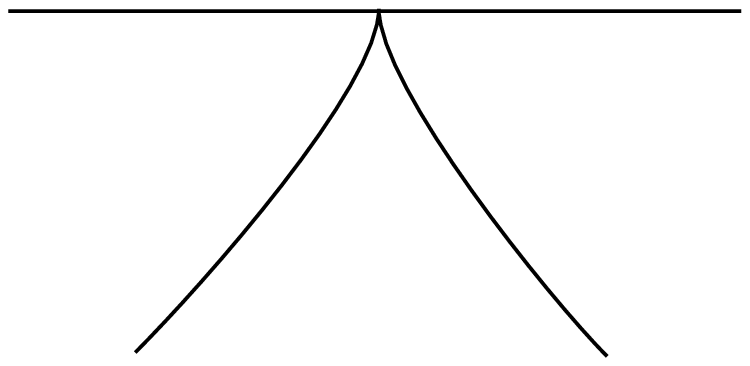}
   \end{center}
 \end{minipage}
$\leftrightarrow $
 \begin{minipage}{0.30\hsize}
  \begin{center}
  \includegraphics*[width=3cm,height=3cm]{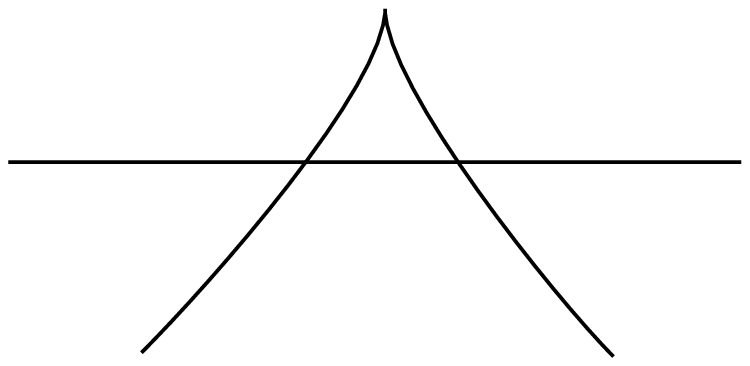}
  \end{center}
\end{minipage}\\
 \begin{minipage}{0.30\hsize} 
  \begin{center}
    \includegraphics*[width=3cm,height=3cm]{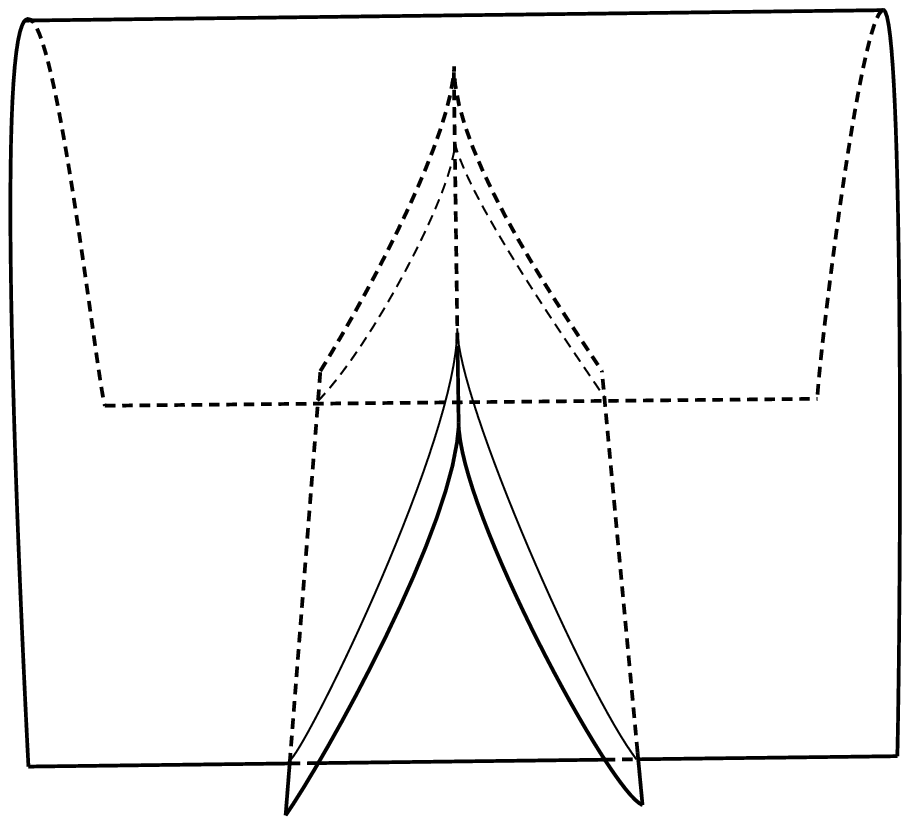}
  \end{center}
 \end{minipage}
$\leftrightarrow $
 \begin{minipage}{0.30\hsize}
   \begin{center}
     \includegraphics*[width=3cm,height=3cm]{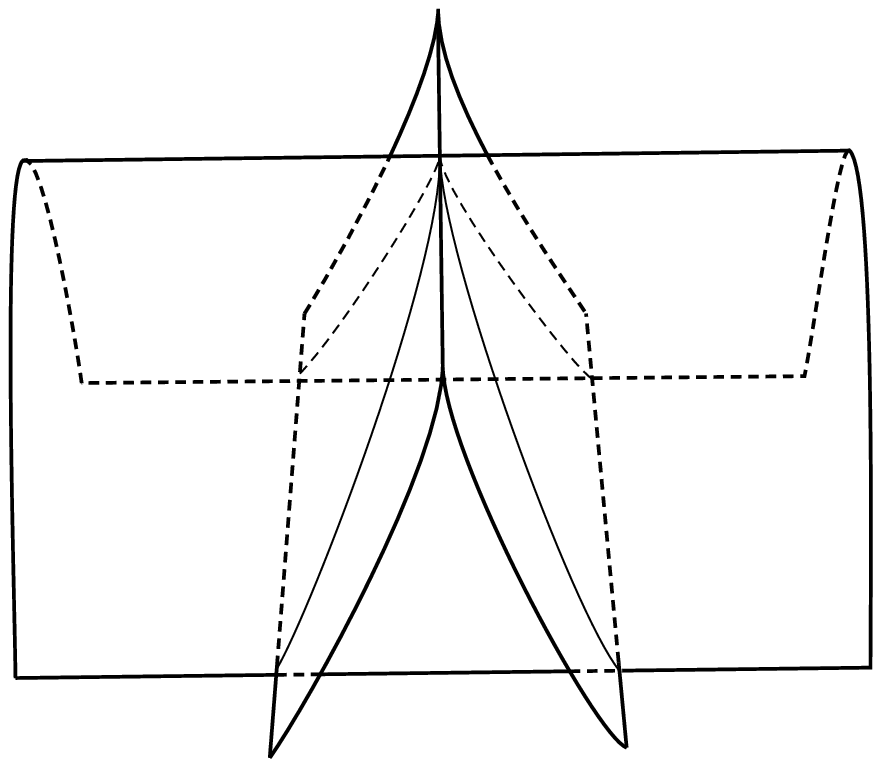}
   \end{center}
 \end{minipage}
$\leftrightarrow $
 \begin{minipage}{0.30\hsize}
  \begin{center}
  \includegraphics*[width=3cm,height=3cm]{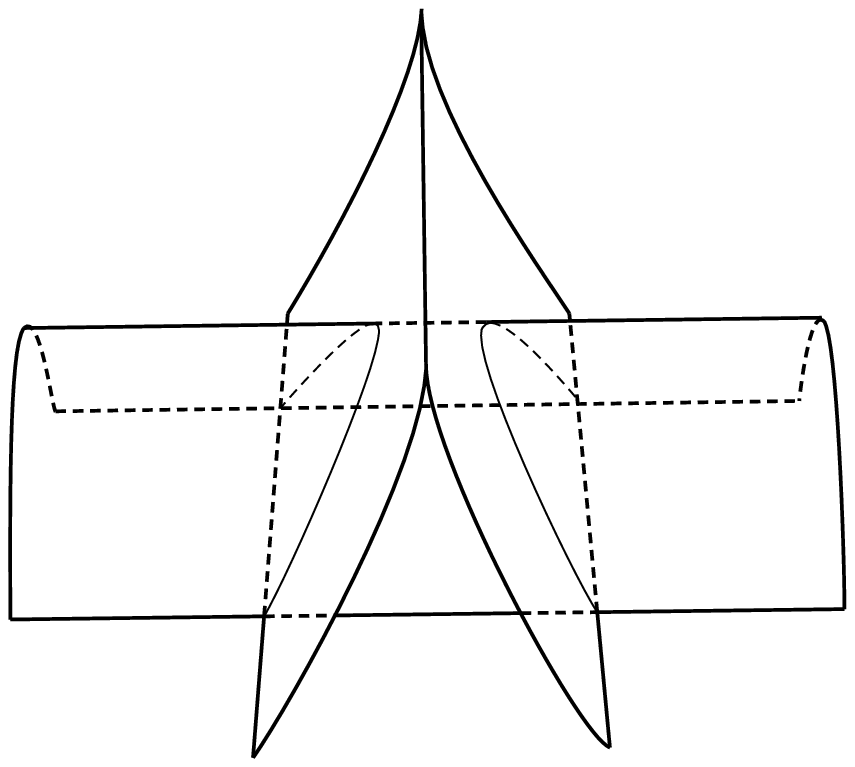}
  \end{center}
\end{minipage}\\
 \begin{minipage}{0.30\hsize} 
  \begin{center}
    \includegraphics*[width=3cm,height=3cm]{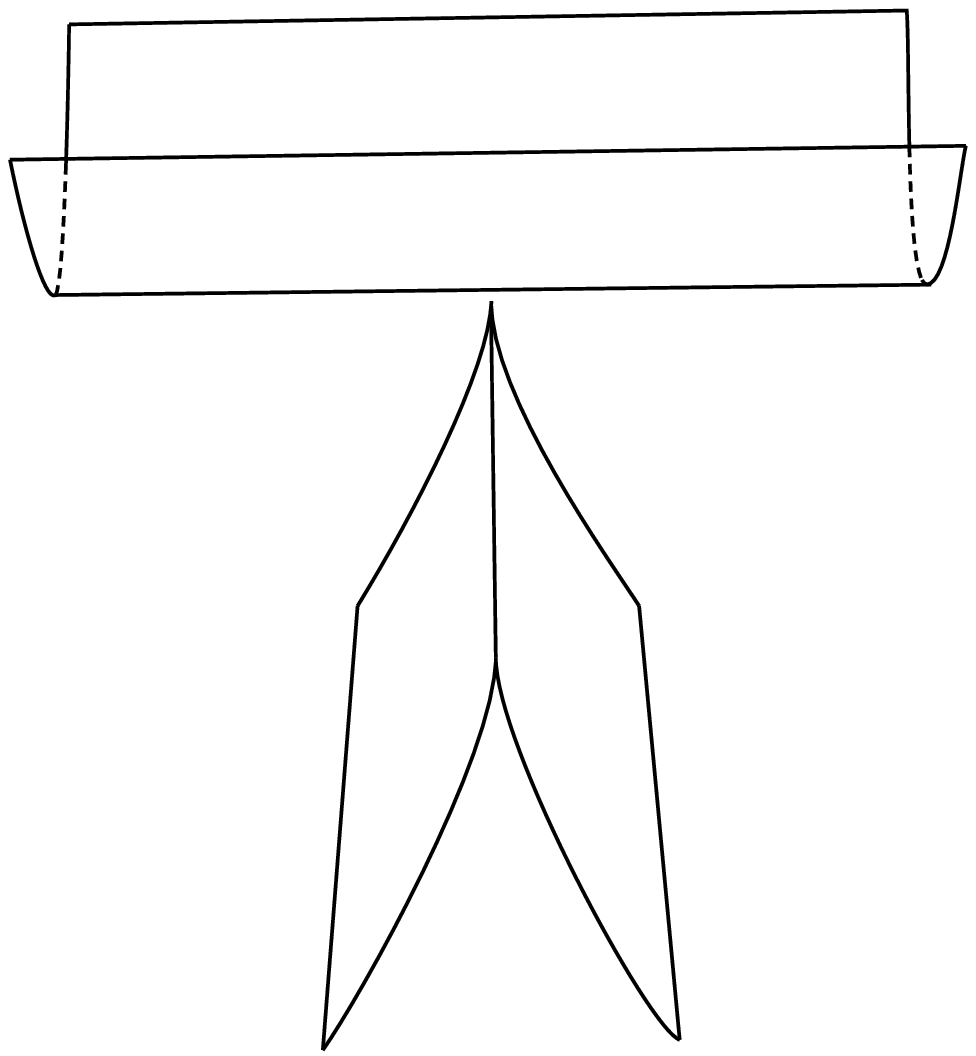}
  \end{center}
 \end{minipage}
$\leftrightarrow $
 \begin{minipage}{0.30\hsize}
   \begin{center}
     \includegraphics*[width=3cm,height=3cm]{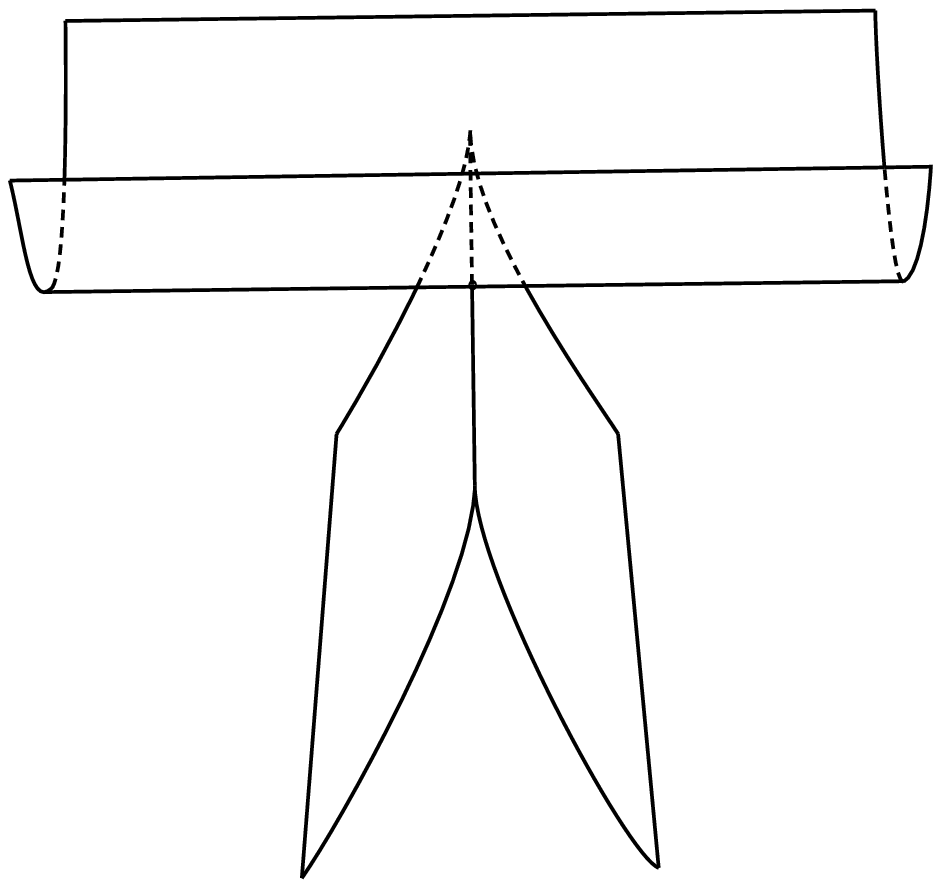}
   \end{center}
 \end{minipage}
$\leftrightarrow $
 \begin{minipage}{0.30\hsize}
  \begin{center}
  \includegraphics*[width=3cm,height=3cm]{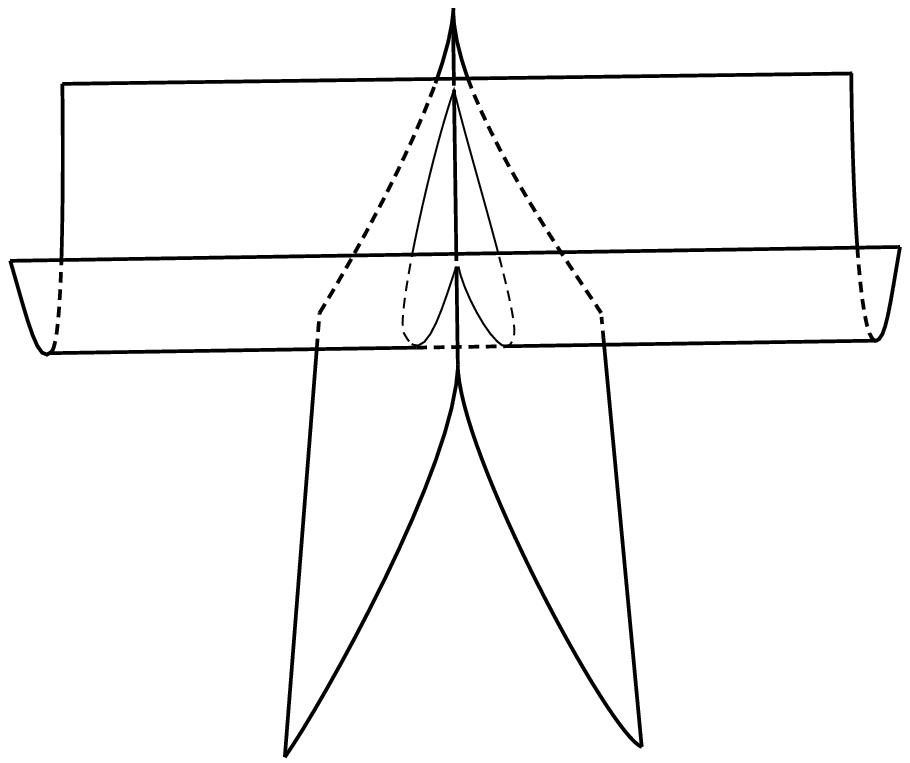}
  \end{center}
\end{minipage}
\caption{${}^1({}^0A_1{}^0A_2)$}
\end{figure}
\begin{figure}[htbp]
 \begin{minipage}{0.30\hsize} 
  \begin{center}
    \includegraphics*[width=3cm,height=3cm]{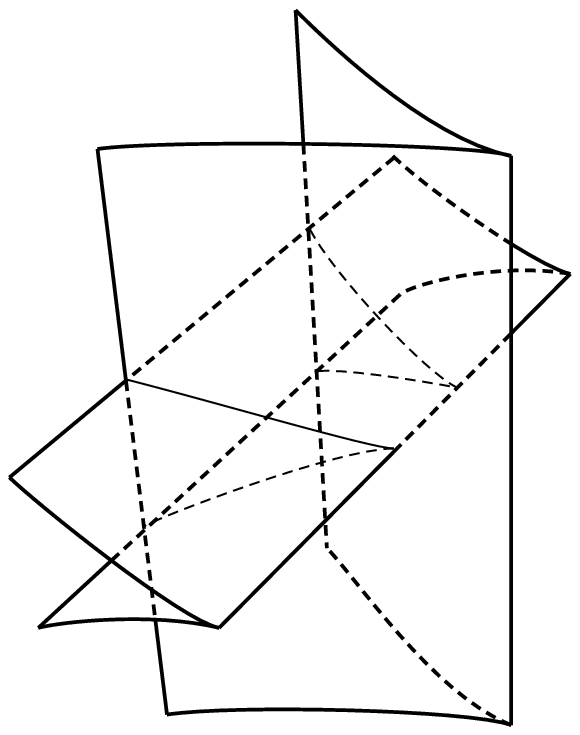}
  \end{center}
 \end{minipage}
$\leftrightarrow $
 \begin{minipage}{0.30\hsize}
   \begin{center}
     \includegraphics*[width=3cm,height=3cm]{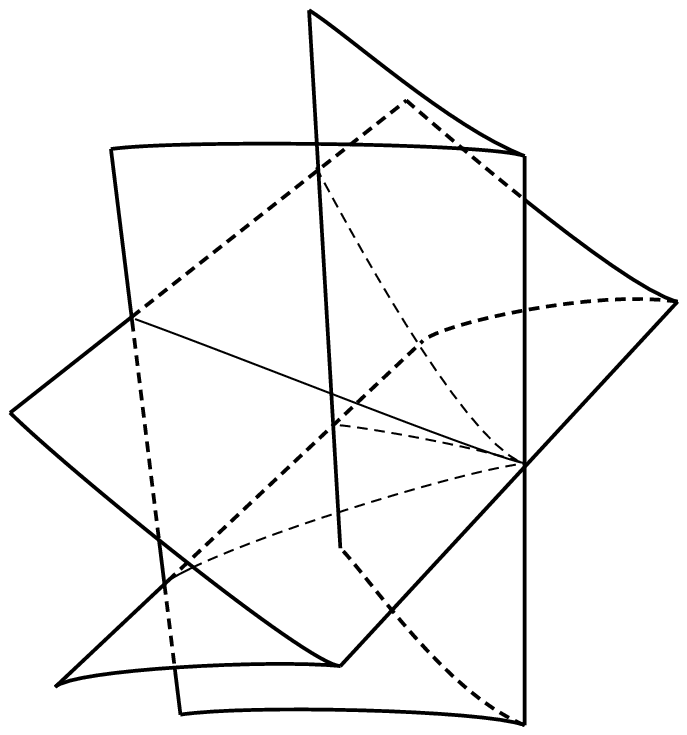}
   \end{center}
 \end{minipage}
$\leftrightarrow $
 \begin{minipage}{0.30\hsize}
  \begin{center}
  \includegraphics*[width=3cm,height=3cm]{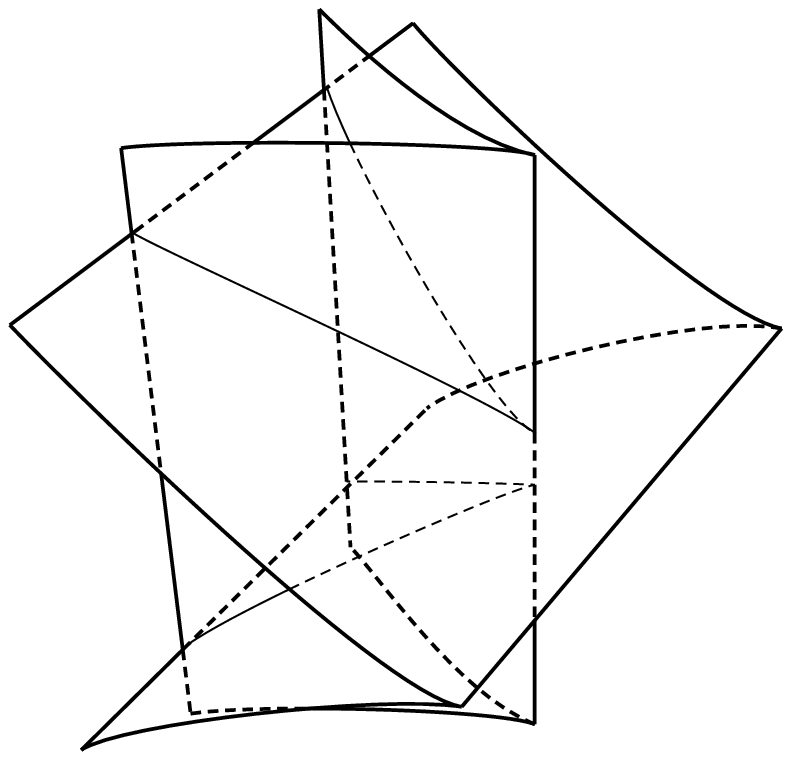}
  \end{center}
\end{minipage}\\
 \begin{minipage}{0.30\hsize} 
  \begin{center}
    \includegraphics*[width=3cm,height=3cm]{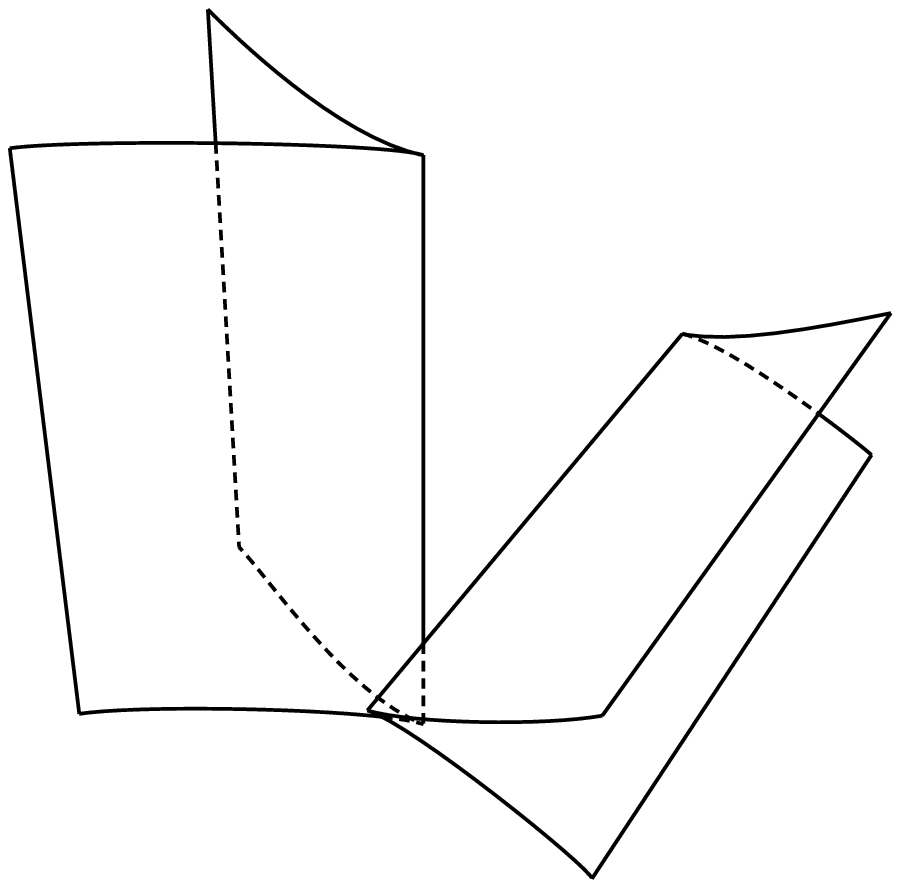}
  \end{center}
 \end{minipage}
$\leftrightarrow $
 \begin{minipage}{0.30\hsize}
   \begin{center}
     \includegraphics*[width=3cm,height=3cm]{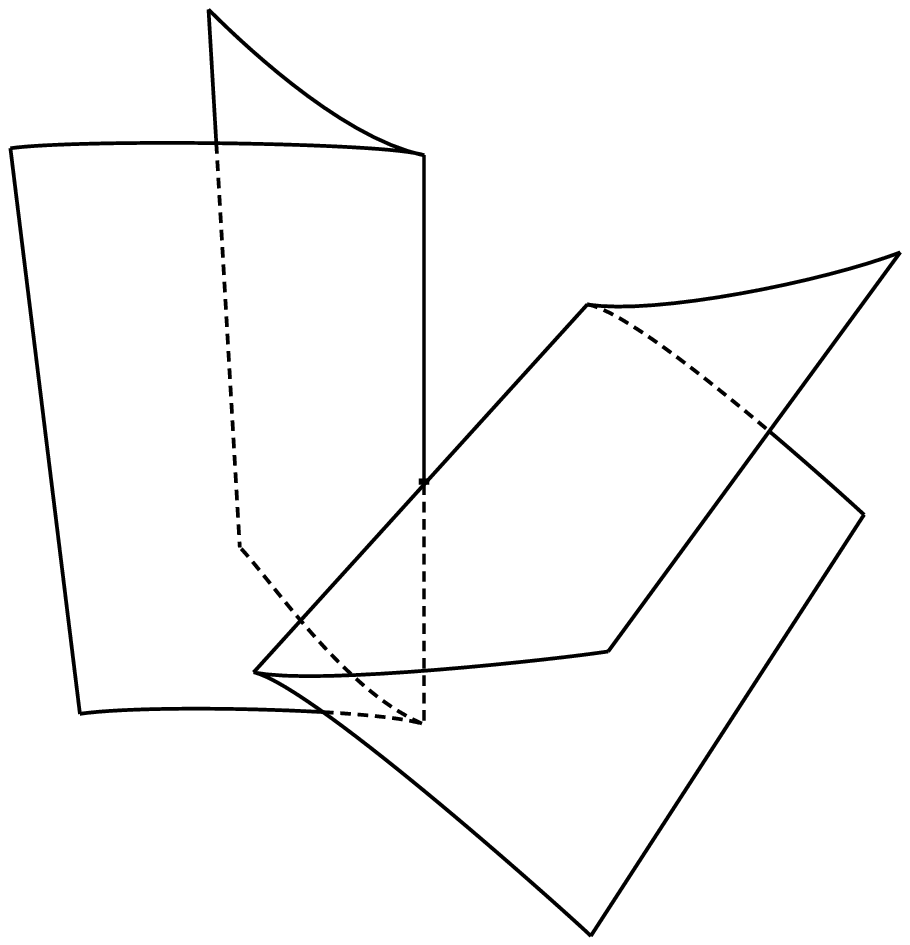}
   \end{center}
 \end{minipage}
$\leftrightarrow $
 \begin{minipage}{0.30\hsize}
  \begin{center}
  \includegraphics*[width=3cm,height=3cm]{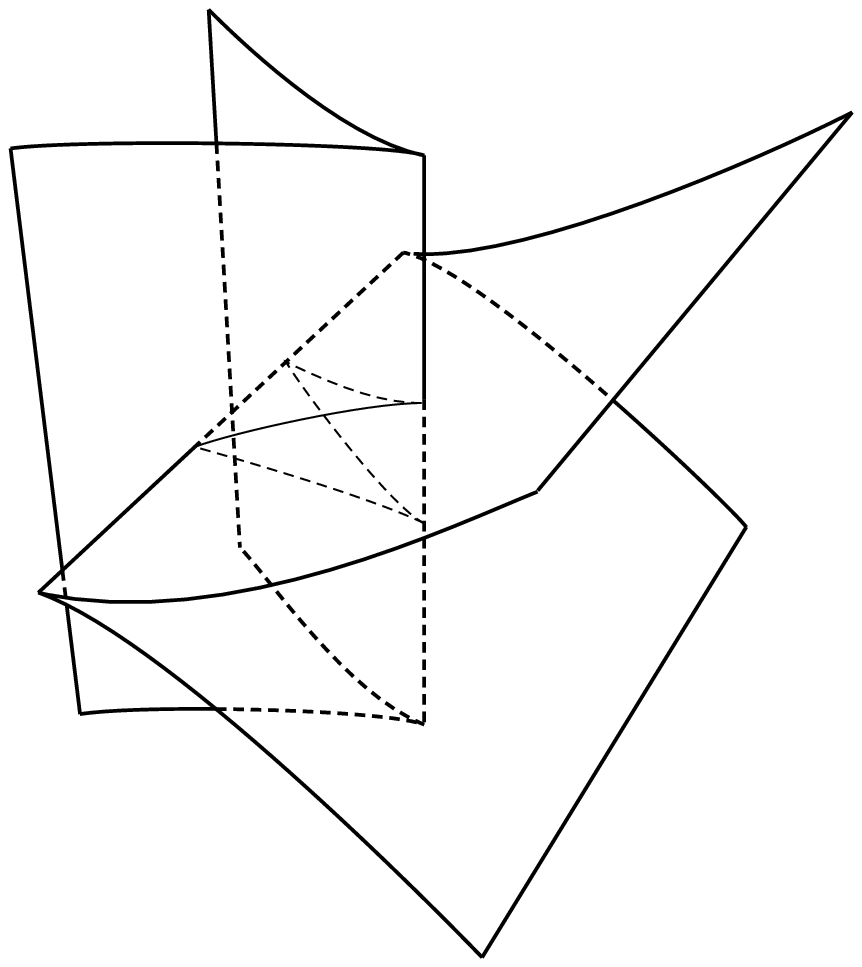}
  \end{center}
\end{minipage}
\caption{${}^1({}^0A_2{}^0A_2)$}
\end{figure}
\begin{figure}[htbp]
 \begin{minipage}{0.30\hsize} 
  \begin{center}
    \includegraphics*[width=3cm,height=3cm]{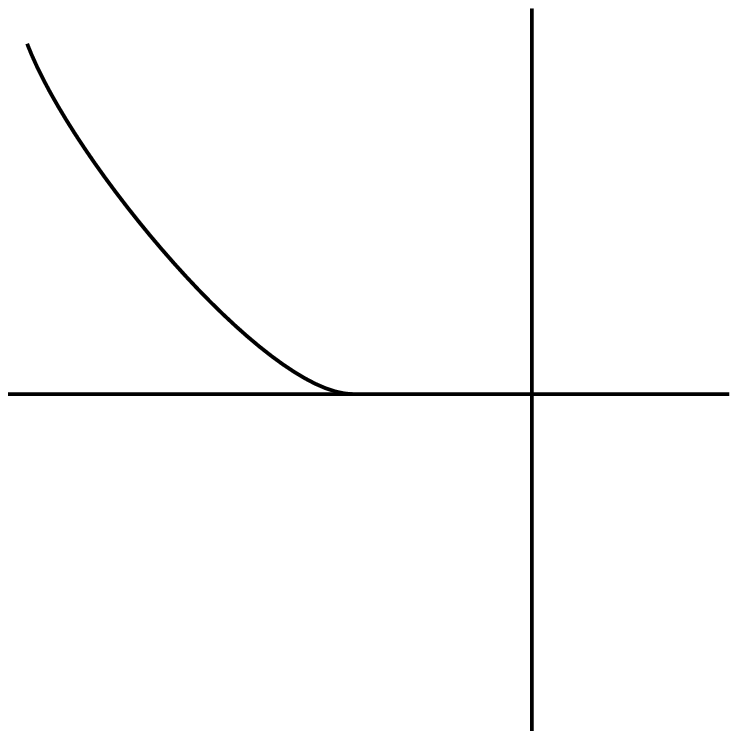}
  \end{center}
 \end{minipage}
$\leftrightarrow $
 \begin{minipage}{0.30\hsize}
   \begin{center}
     \includegraphics*[width=3cm,height=3cm]{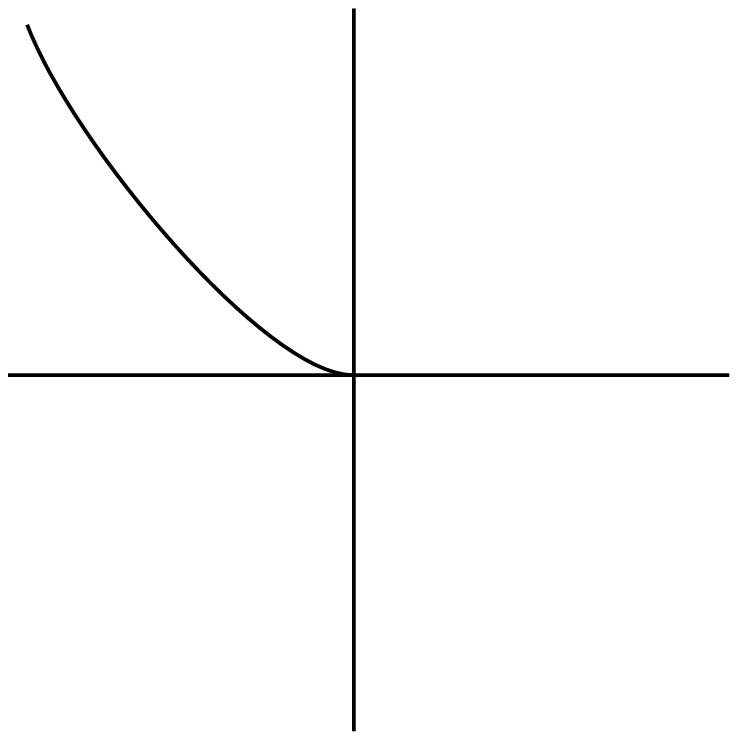}
   \end{center}
 \end{minipage}
$\leftrightarrow $
 \begin{minipage}{0.30\hsize}
  \begin{center}
  \includegraphics*[width=3cm,height=3cm]{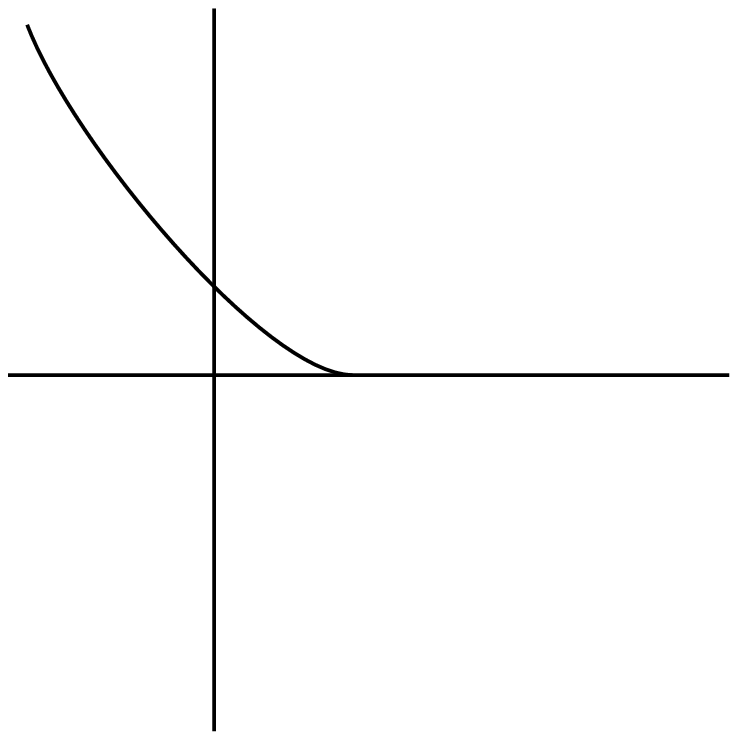}
  \end{center}
\end{minipage}\\
 \begin{minipage}{0.30\hsize} 
  \begin{center}
    \includegraphics*[width=3cm,height=3cm]{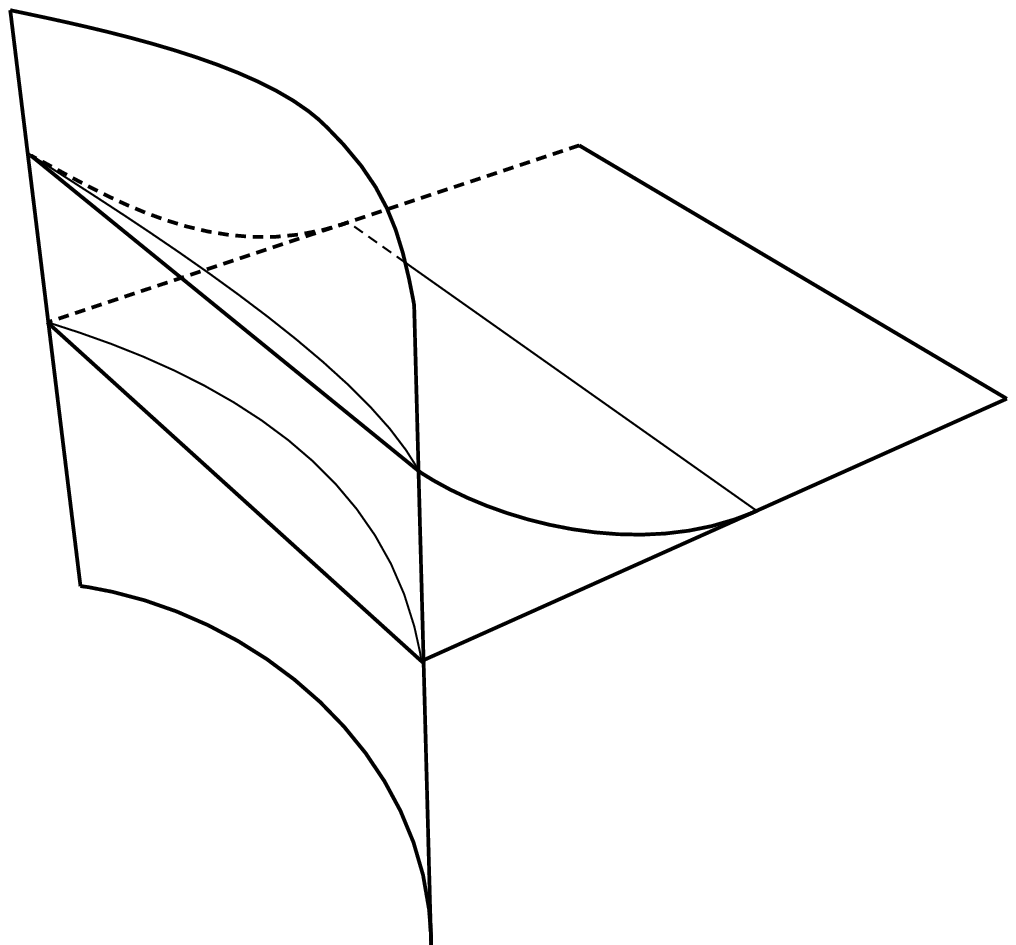}
  \end{center}
 \end{minipage}
$\leftrightarrow $
 \begin{minipage}{0.30\hsize}
   \begin{center}
     \includegraphics*[width=3cm,height=3cm]{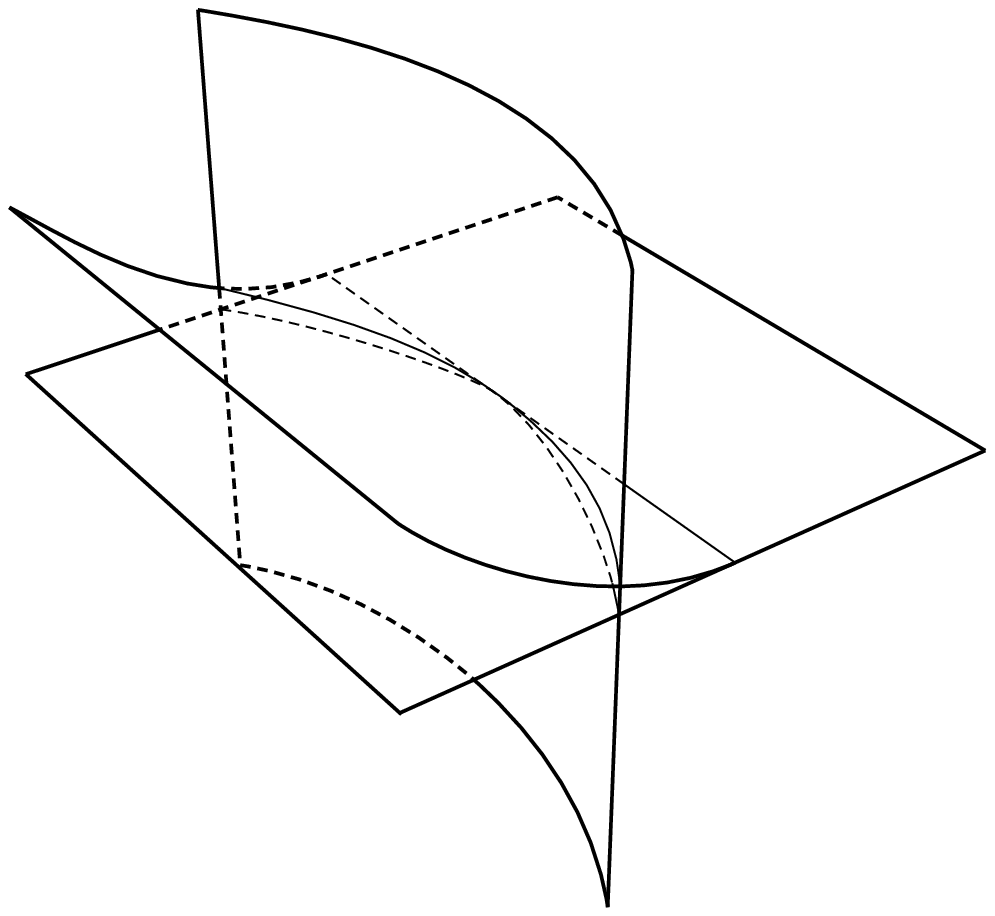}
   \end{center}
 \end{minipage}
$\leftrightarrow $
 \begin{minipage}{0.30\hsize}
  \begin{center}
  \includegraphics*[width=3cm,height=3cm]{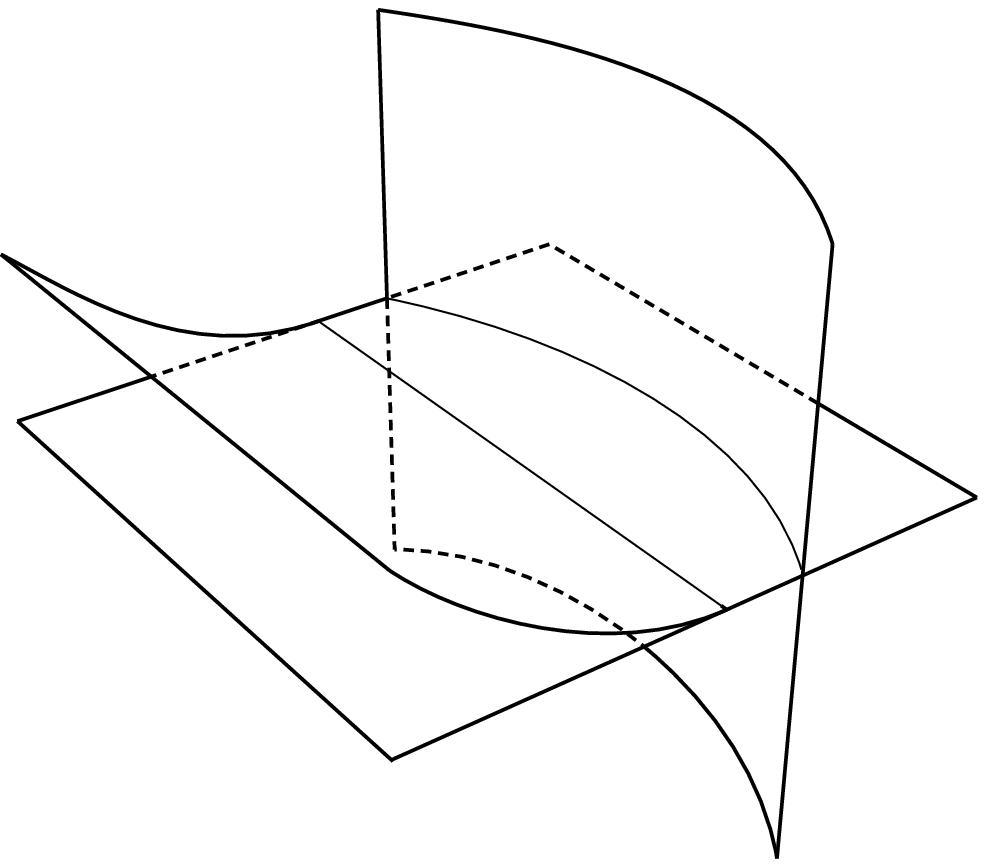}
  \end{center}
\end{minipage}\\
 \begin{minipage}{0.30\hsize} 
  \begin{center}
    \includegraphics*[width=3cm,height=3cm]{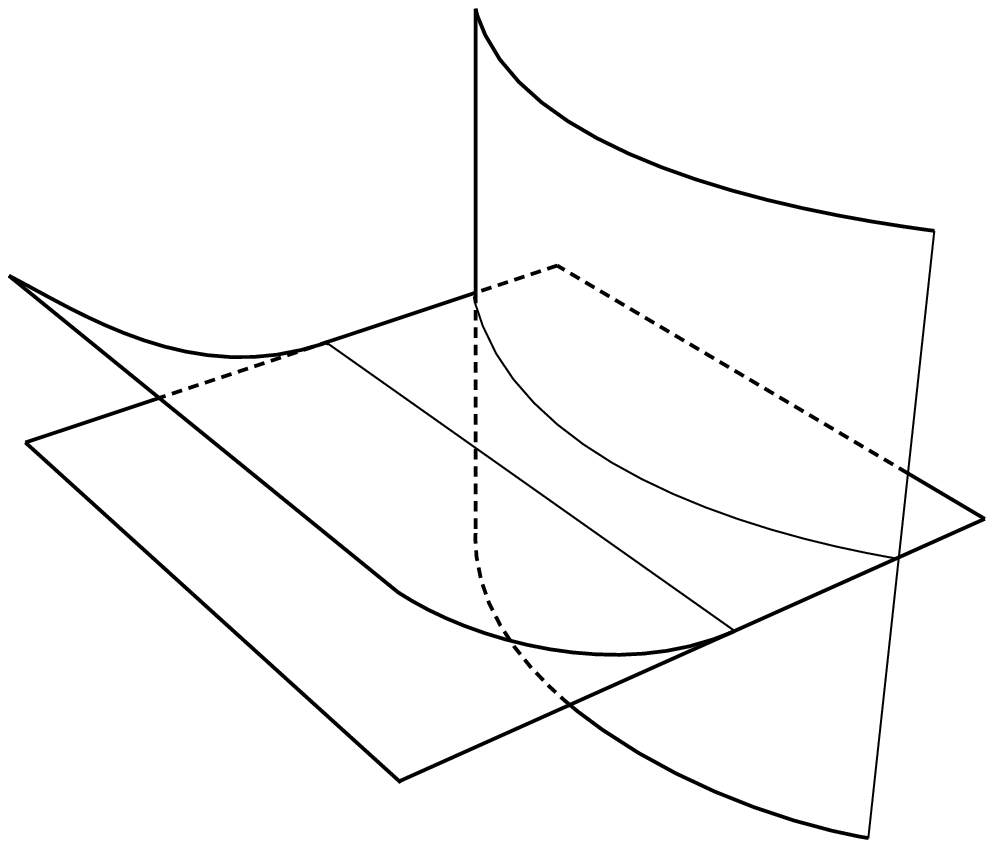}
  \end{center}
 \end{minipage}
$\leftrightarrow $
 \begin{minipage}{0.30\hsize}
   \begin{center}
     \includegraphics*[width=3cm,height=3cm]{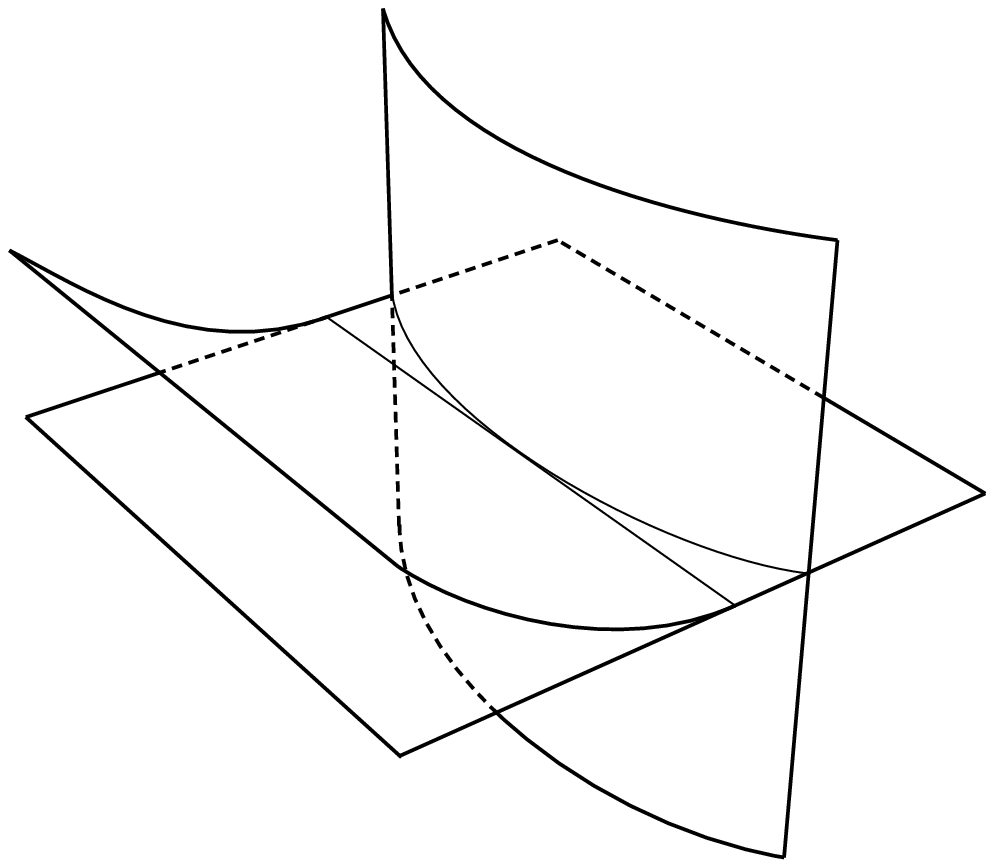}
   \end{center}
 \end{minipage}
$\leftrightarrow $
 \begin{minipage}{0.30\hsize}
  \begin{center}
  \includegraphics*[width=3cm,height=3cm]{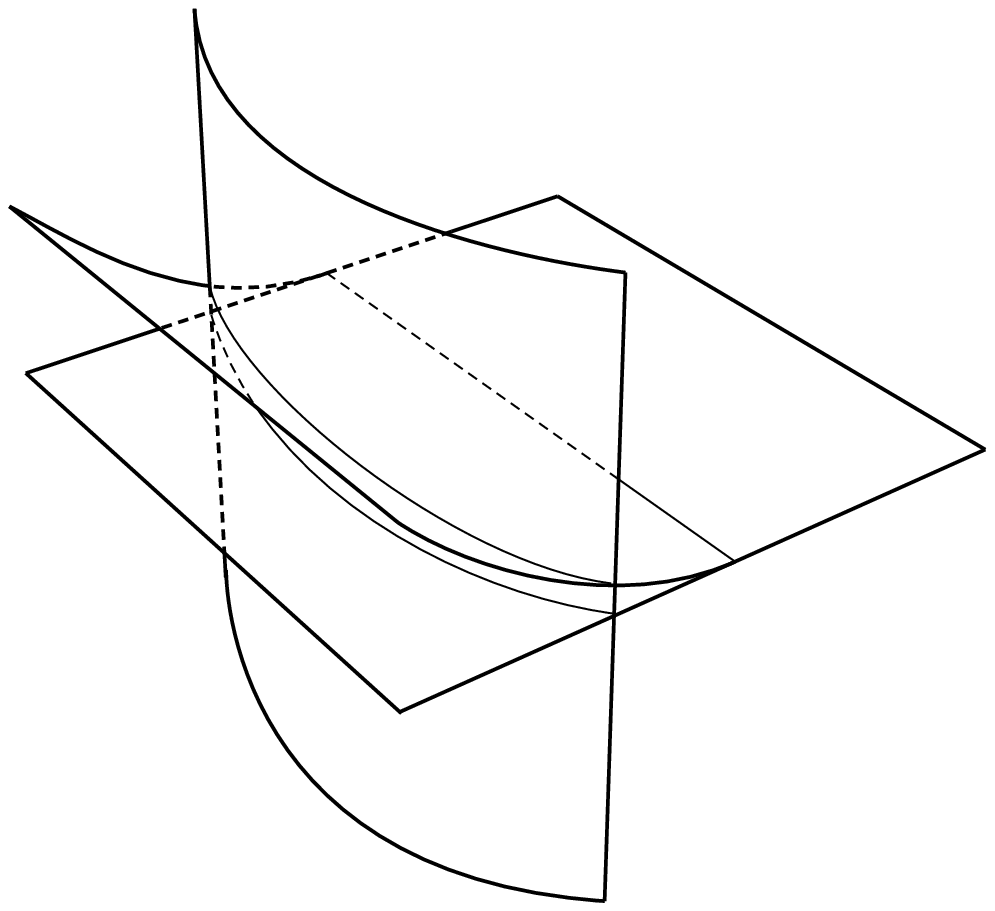}
  \end{center}
\end{minipage}
\caption{${}^1({}^0A_1{}^0B_2)$}
\end{figure}
\begin{figure}[htbp]
 \begin{minipage}{0.30\hsize} 
  \begin{center}
    \includegraphics*[width=3cm,height=3cm]{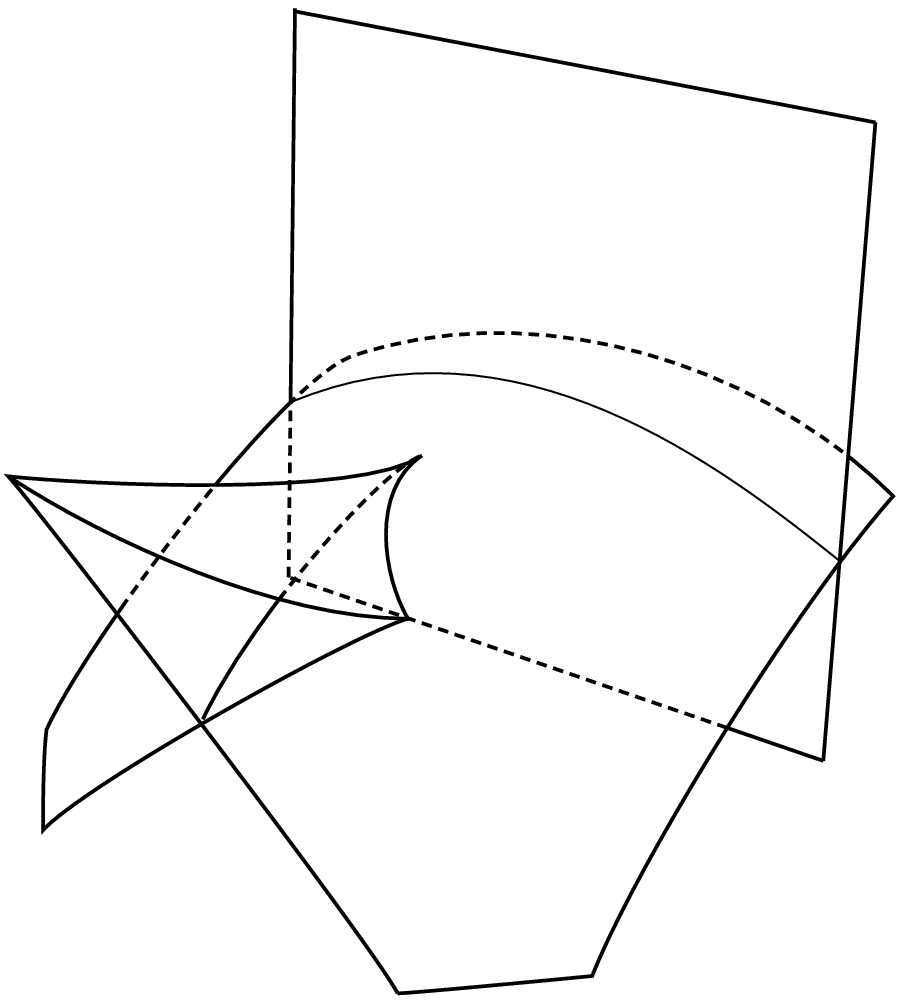}
  \end{center}
 \end{minipage}
$\leftrightarrow $
 \begin{minipage}{0.30\hsize}
   \begin{center}
     \includegraphics*[width=3cm,height=3cm]{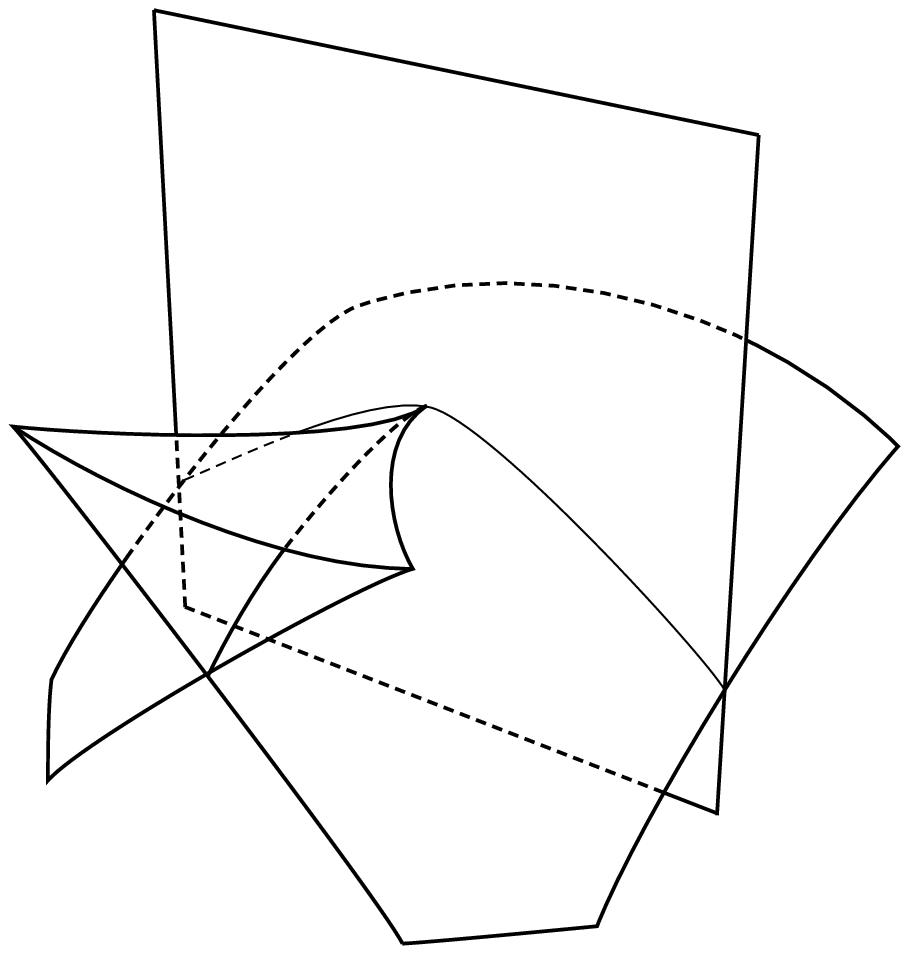}
   \end{center}
 \end{minipage}
$\leftrightarrow $
 \begin{minipage}{0.30\hsize}
  \begin{center}
  \includegraphics*[width=3cm,height=3cm]{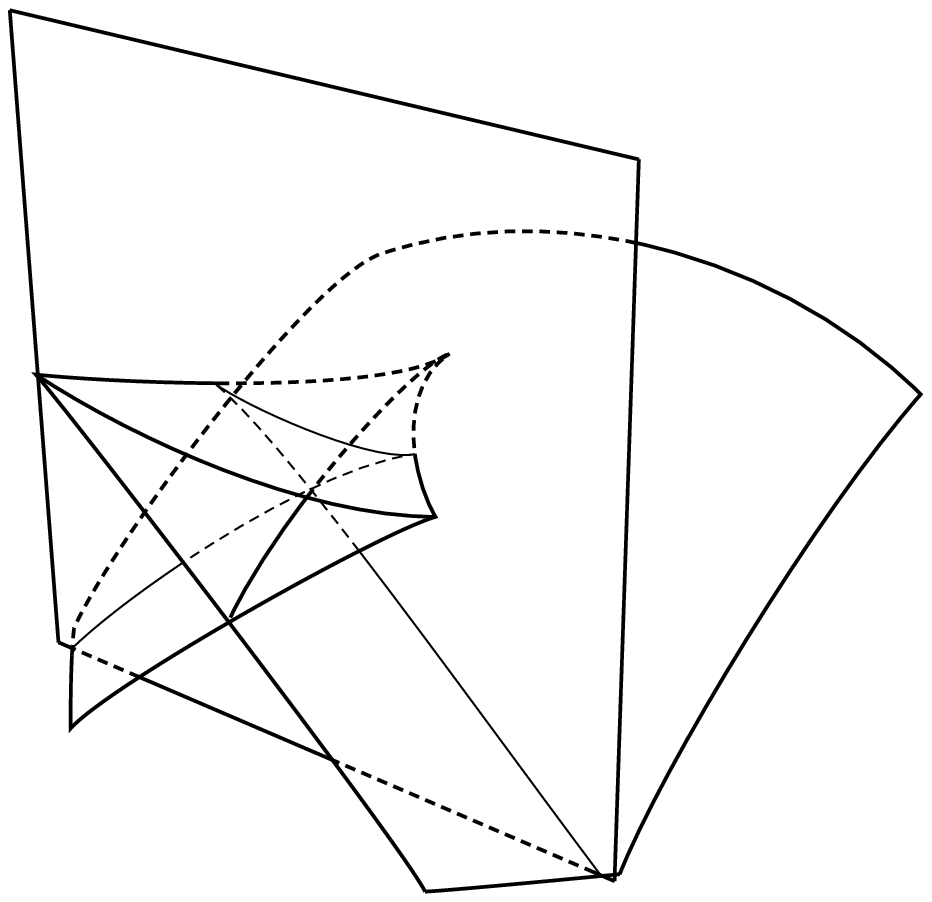}
  \end{center}
\end{minipage}
\caption{${}^1({}^0A_1{}^0A_3)$}
\end{figure}
\begin{figure}[htbp]
 \begin{minipage}{0.30\hsize} 
  \begin{center}
    \includegraphics*[width=3cm,height=3cm]{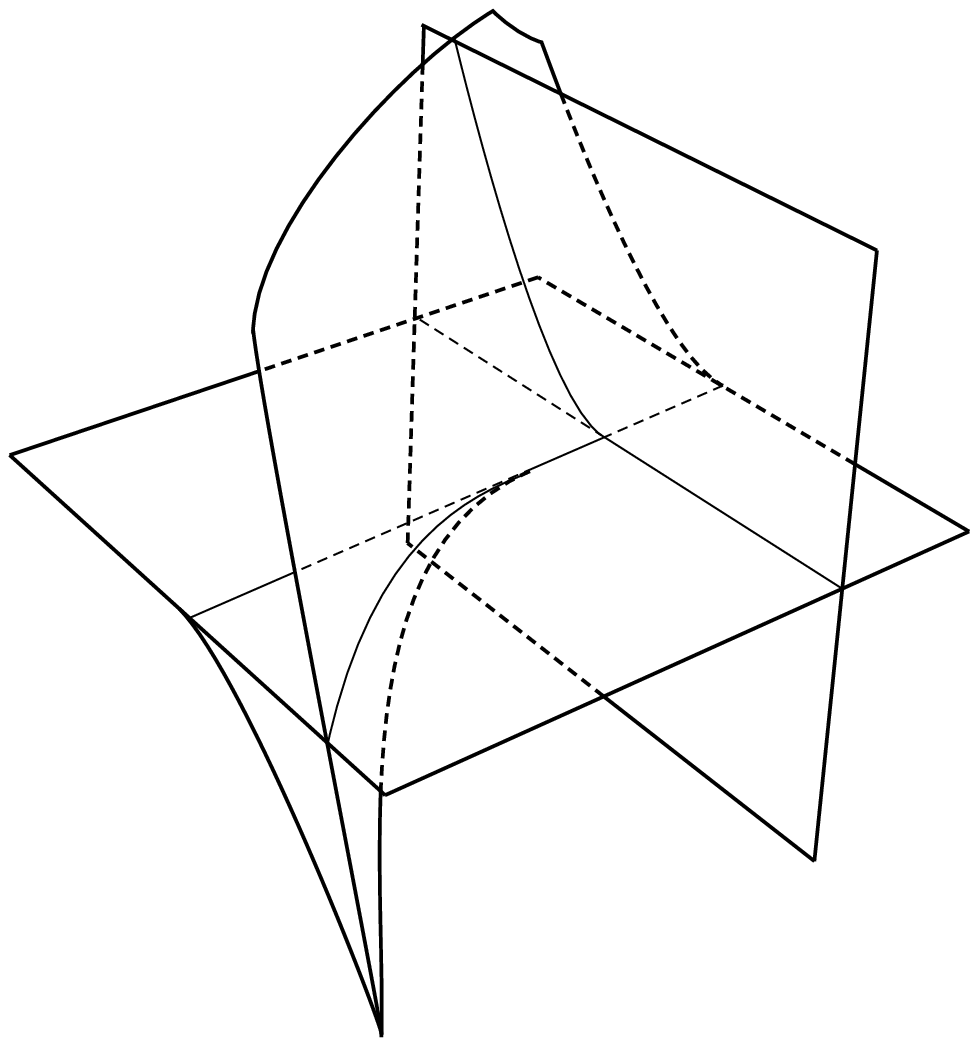}
  \end{center}
 \end{minipage}
$\leftrightarrow $
 \begin{minipage}{0.30\hsize}
   \begin{center}
     \includegraphics*[width=3cm,height=3cm]{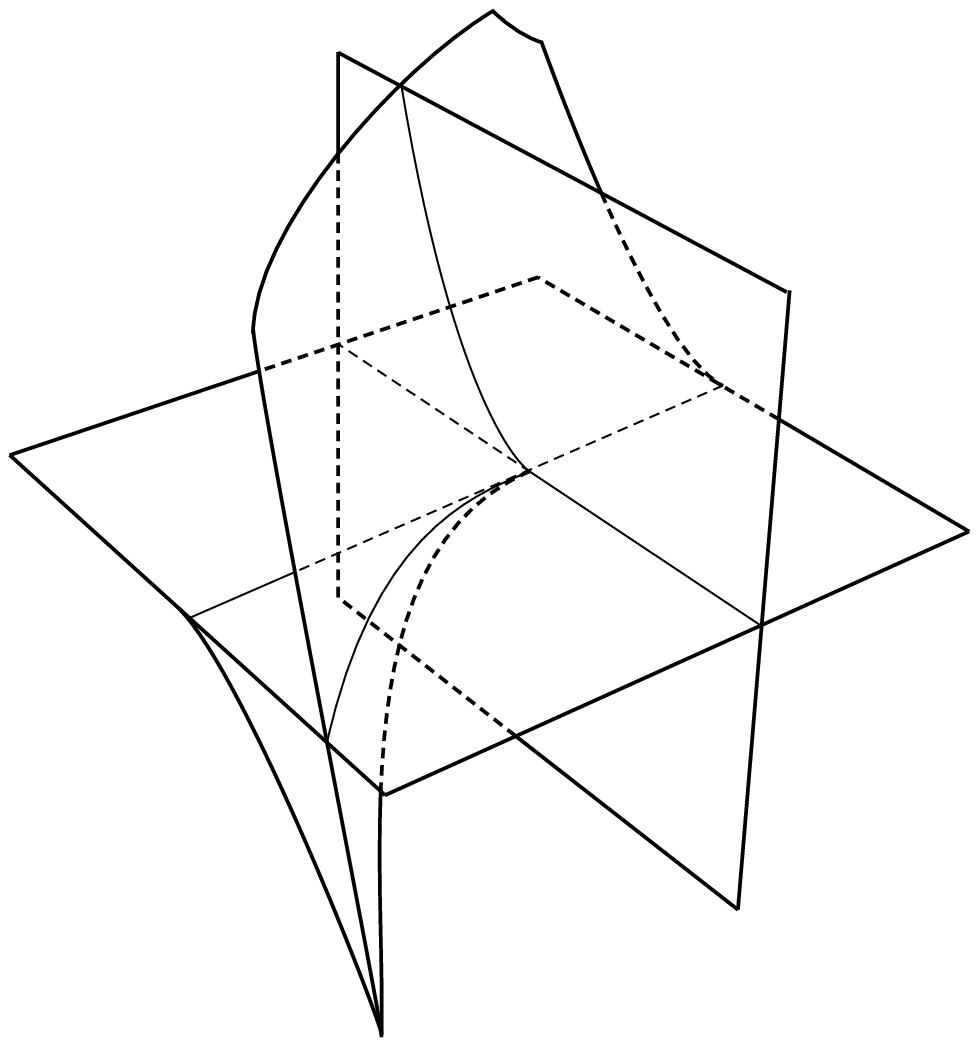}
   \end{center}
 \end{minipage}
$\leftrightarrow $
 \begin{minipage}{0.30\hsize}
  \begin{center}
  \includegraphics*[width=3cm,height=3cm]{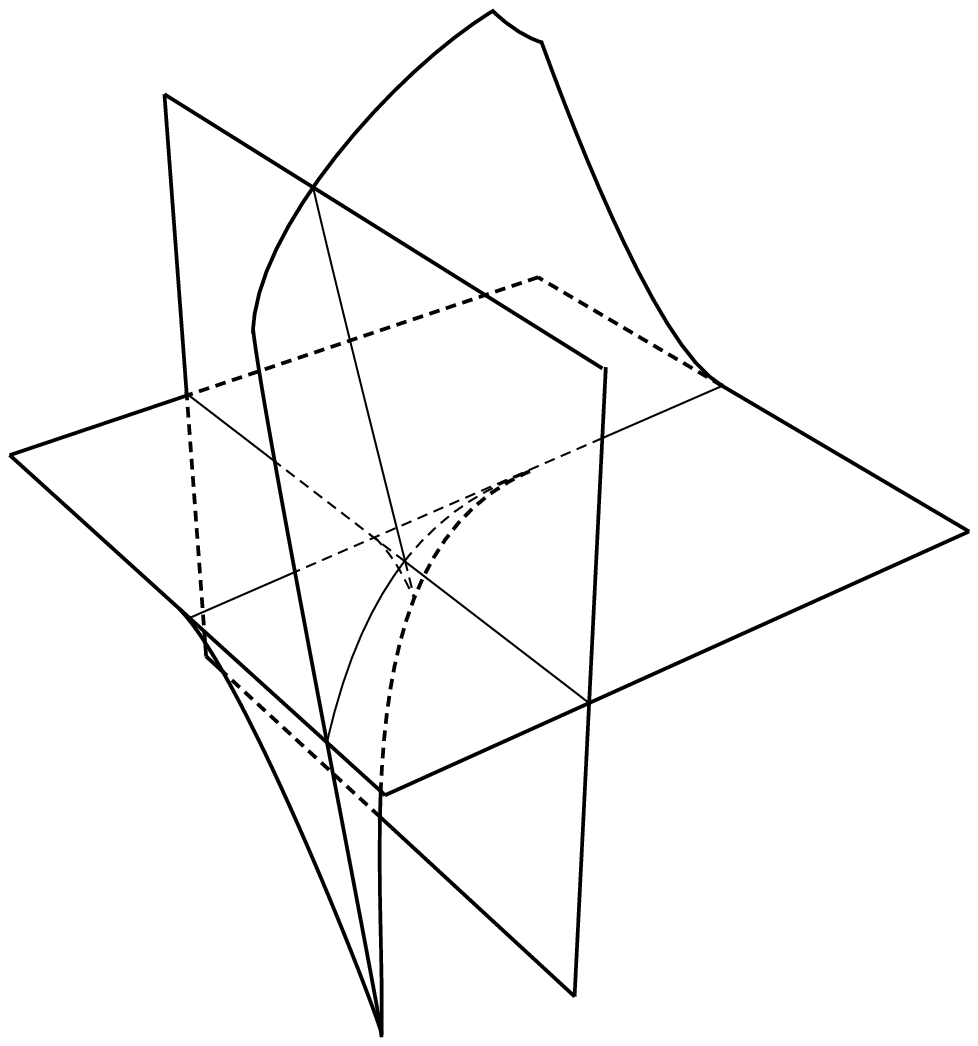}
  \end{center}
\end{minipage}
\caption{${}^1({}^0A_1{}^0B_3)$}
\end{figure}
\begin{figure}[htbp]
 \begin{minipage}{0.30\hsize} 
  \begin{center}
    \includegraphics*[width=3cm,height=3cm]{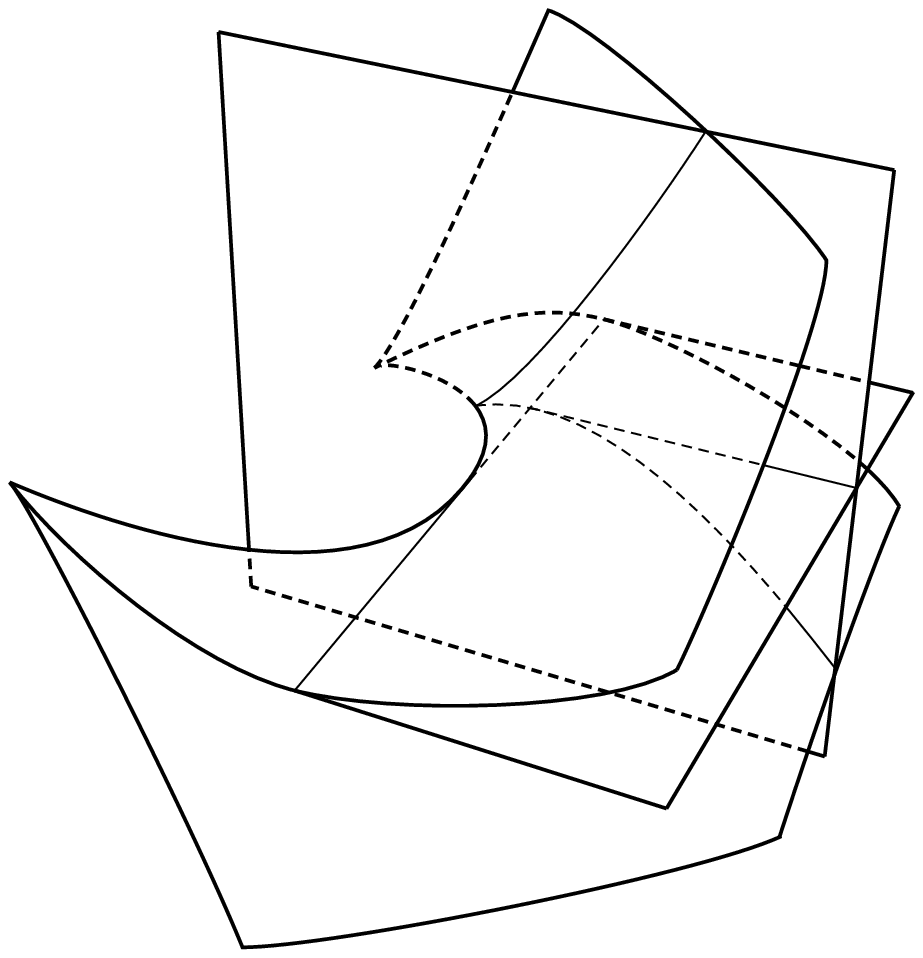}
  \end{center}
 \end{minipage}
$\leftrightarrow $
 \begin{minipage}{0.30\hsize}
   \begin{center}
     \includegraphics*[width=3cm,height=3cm]{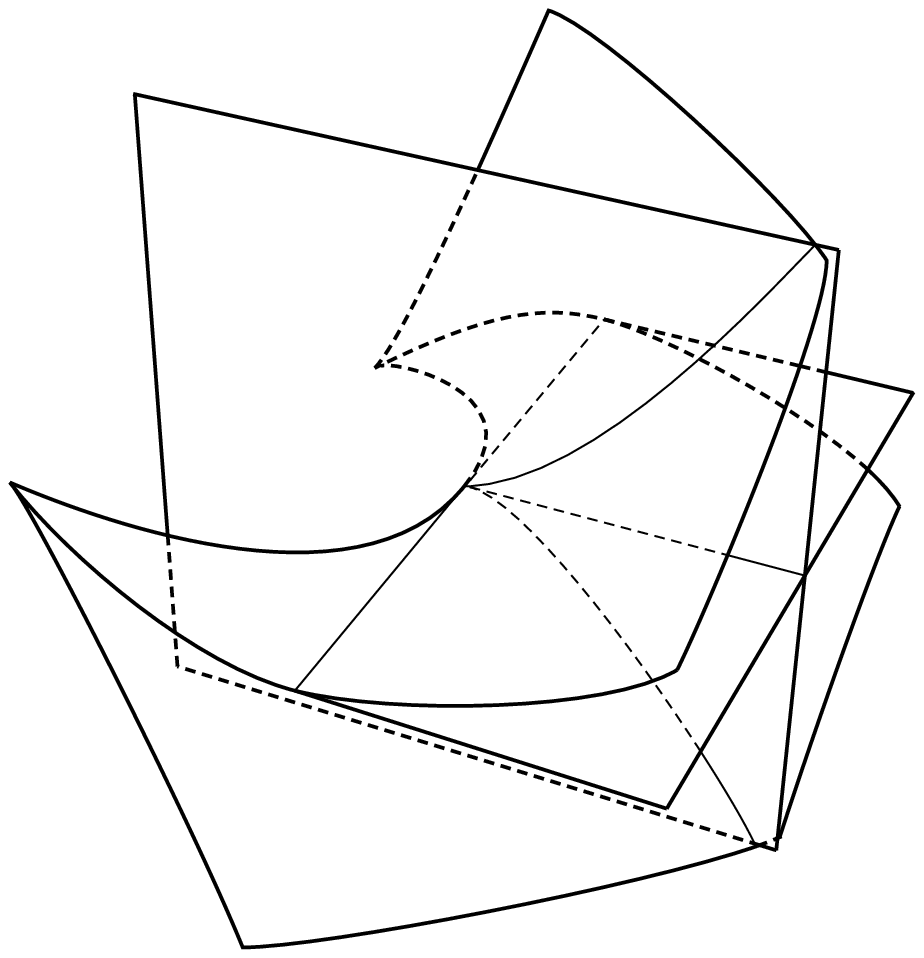}
   \end{center}
 \end{minipage}
$\leftrightarrow $
 \begin{minipage}{0.30\hsize}
  \begin{center}
  \includegraphics*[width=3cm,height=3cm]{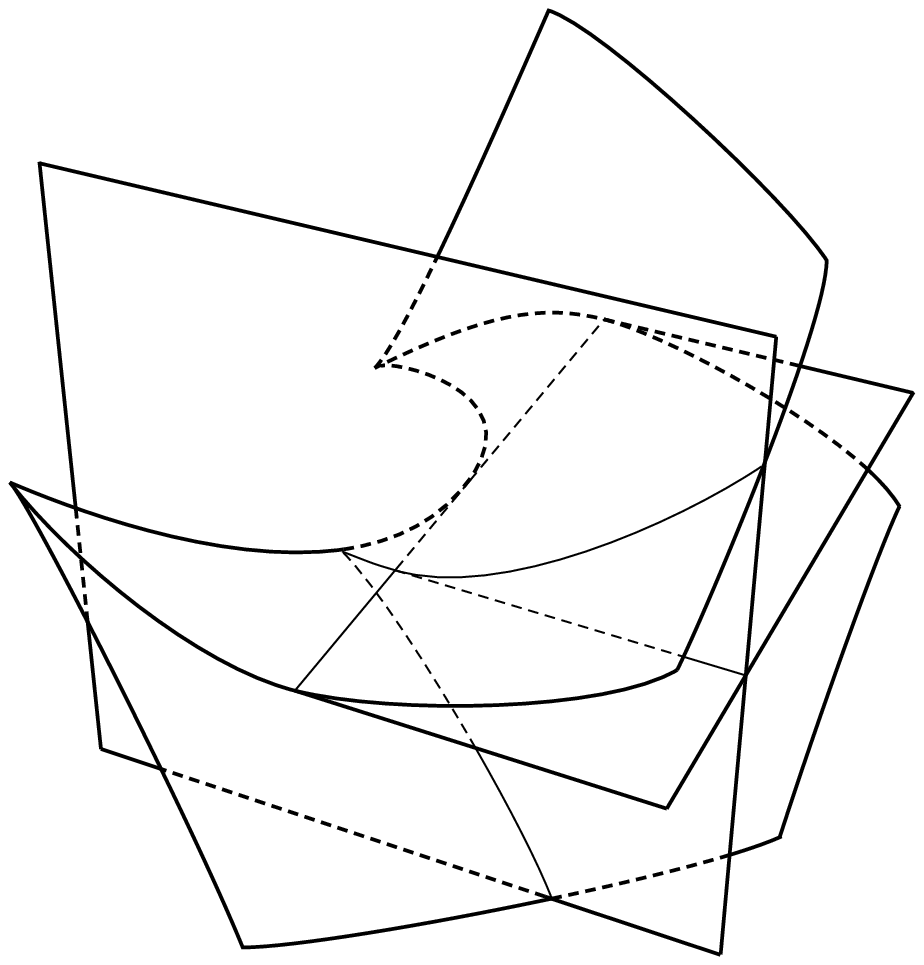}
  \end{center}
\end{minipage}
\caption{${}^1({}^0A_1{}^0C^+_3)$}
\end{figure}
\begin{figure}[htbp]
 \begin{minipage}{0.30\hsize} 
  \begin{center}
    \includegraphics*[width=3cm,height=3cm]{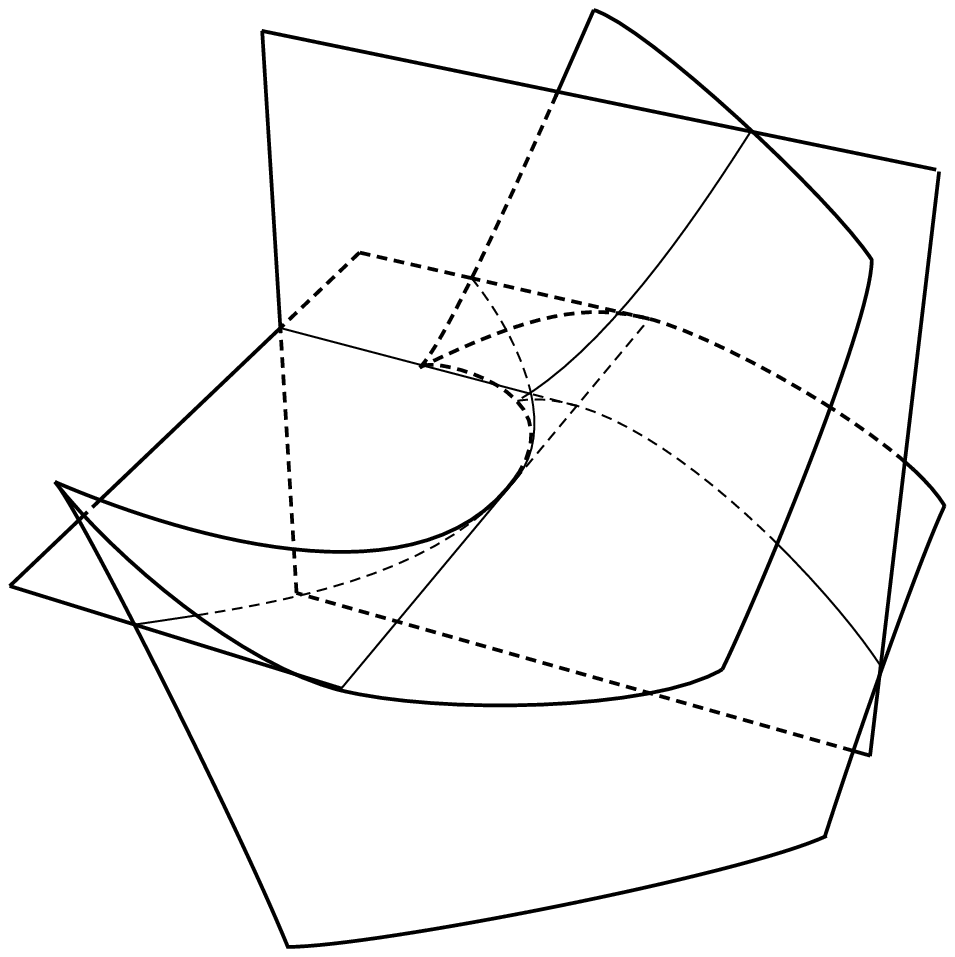}
  \end{center}
 \end{minipage}
$\leftrightarrow $
 \begin{minipage}{0.30\hsize}
   \begin{center}
     \includegraphics*[width=3cm,height=3cm]{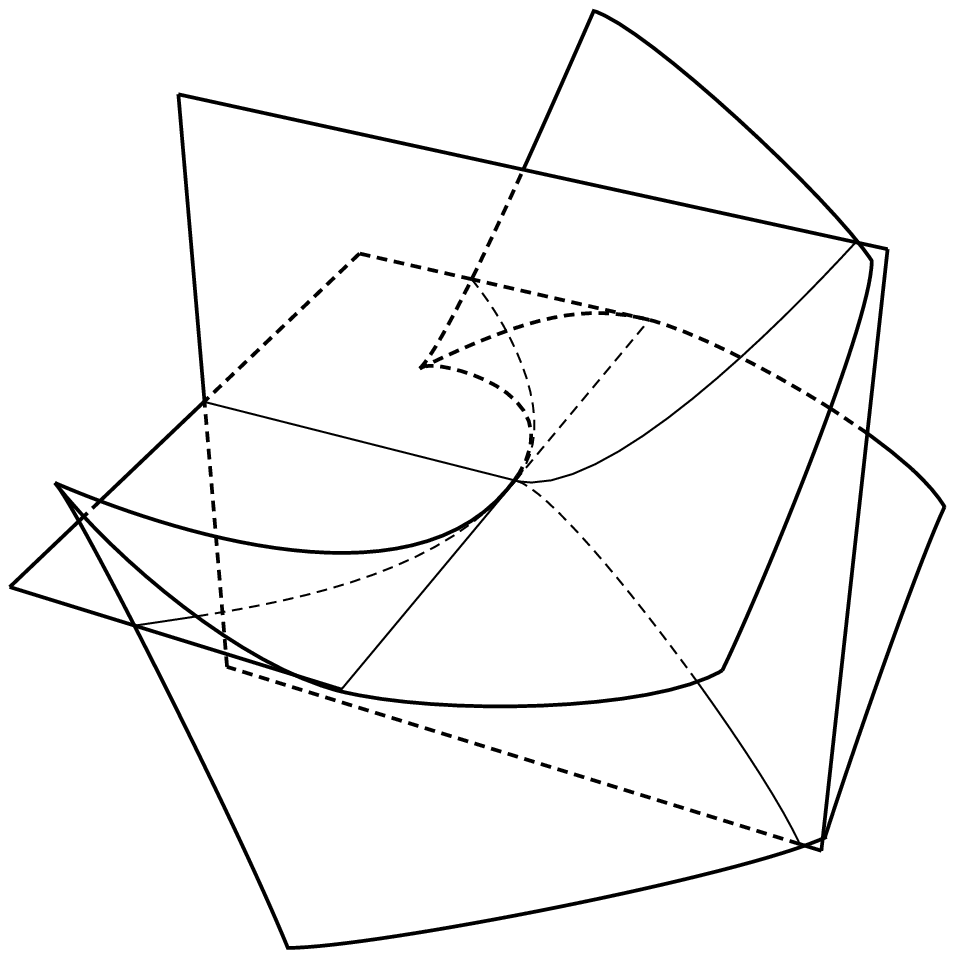}
   \end{center}
 \end{minipage}
$\leftrightarrow $
 \begin{minipage}{0.30\hsize}
  \begin{center}
  \includegraphics*[width=3cm,height=3cm]{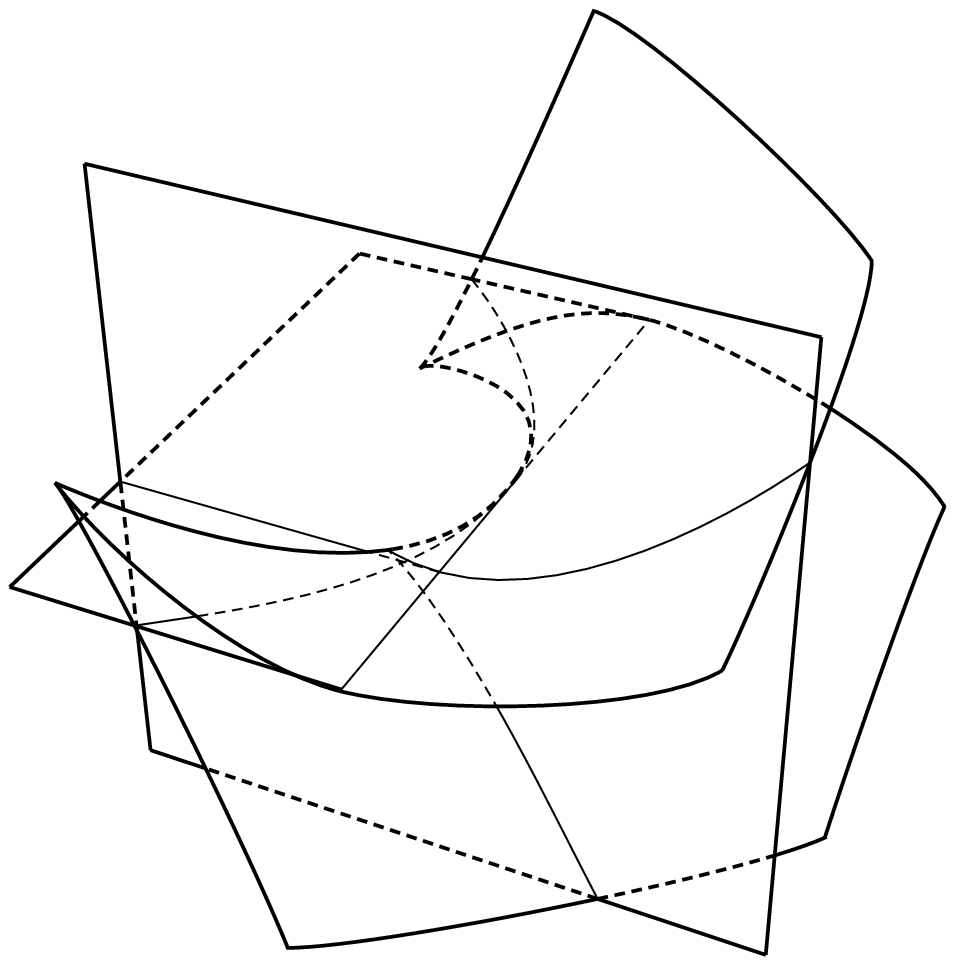}
  \end{center}
\end{minipage}
\caption{${}^1({}^0A_1{}^0C^-_3)$}
\end{figure}
\begin{figure}[htbp]
 \begin{minipage}{0.30\hsize} 
  \begin{center}
    \includegraphics*[width=3cm,height=3cm]{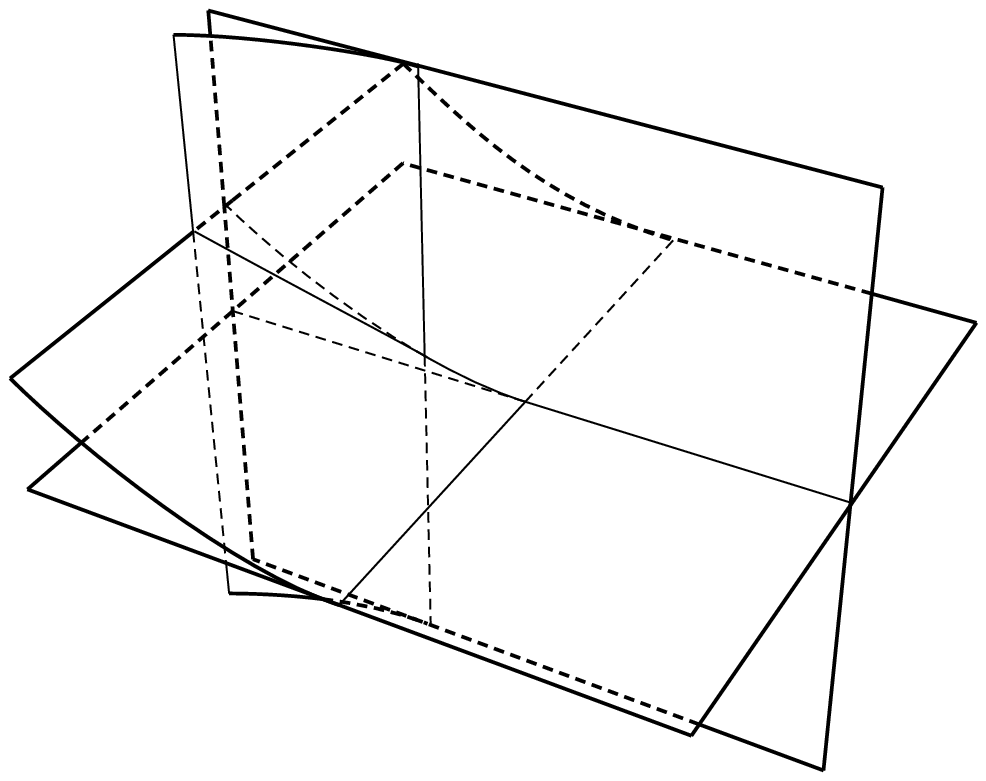}
  \end{center}
 \end{minipage}
$\leftrightarrow $
 \begin{minipage}{0.30\hsize}
   \begin{center}
     \includegraphics*[width=3cm,height=3cm]{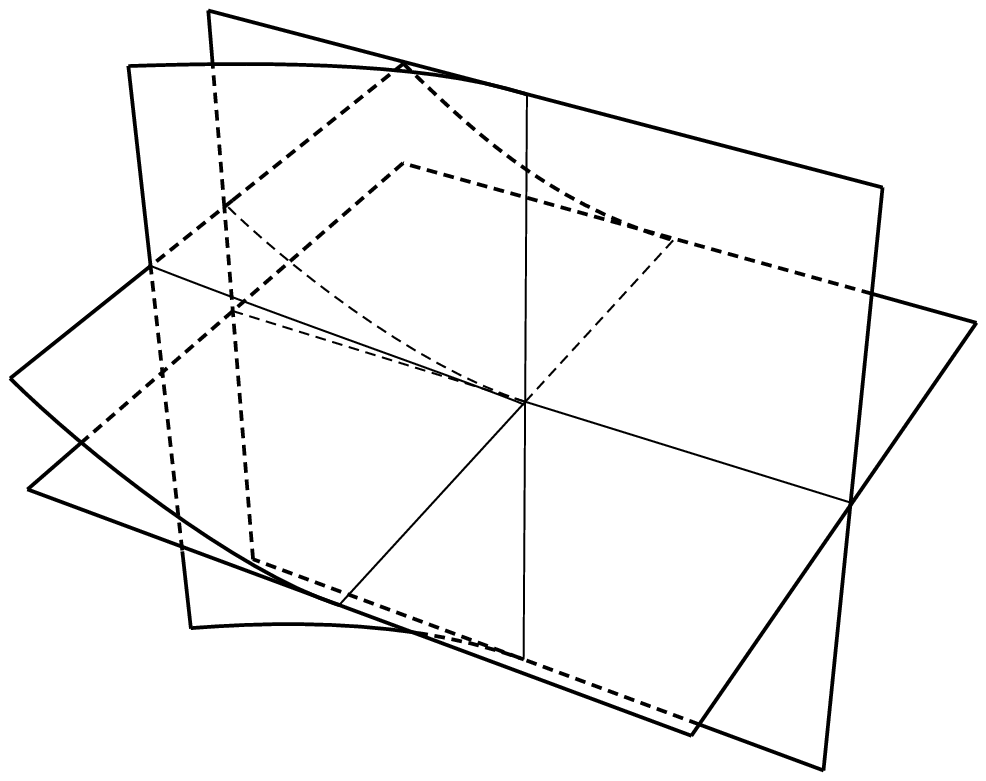}
   \end{center}
 \end{minipage}
$\leftrightarrow $
 \begin{minipage}{0.30\hsize}
  \begin{center}
  \includegraphics*[width=3cm,height=3cm]{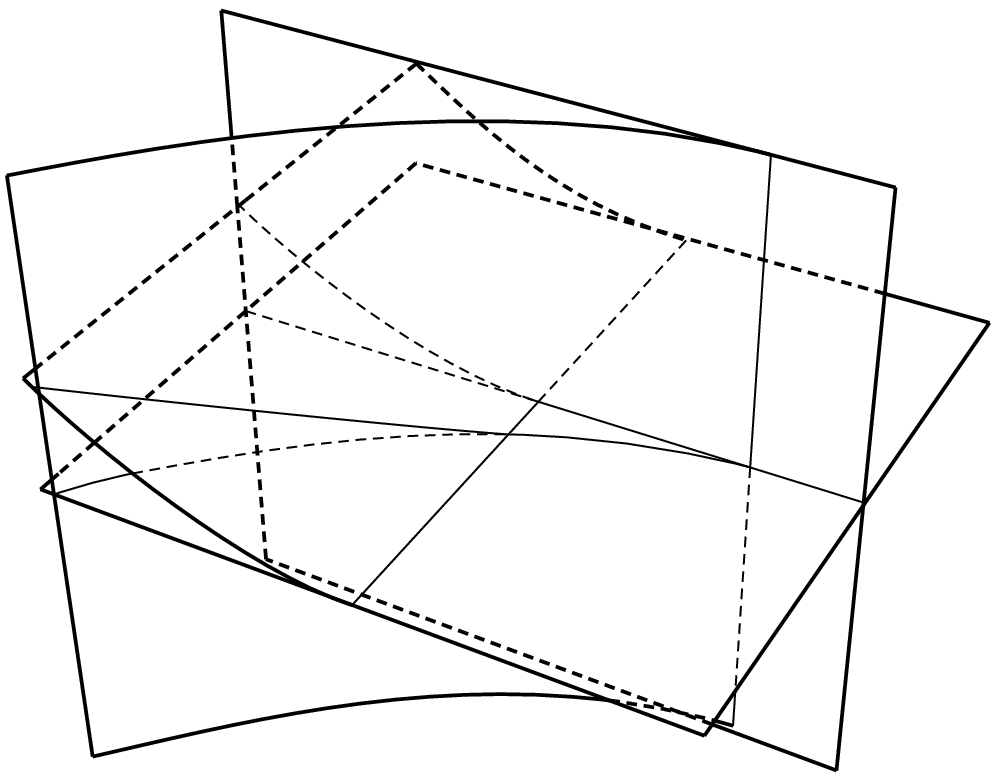}
  \end{center}
\end{minipage}\\
 \begin{minipage}{0.30\hsize} 
  \begin{center}
    \includegraphics*[width=3cm,height=3cm]{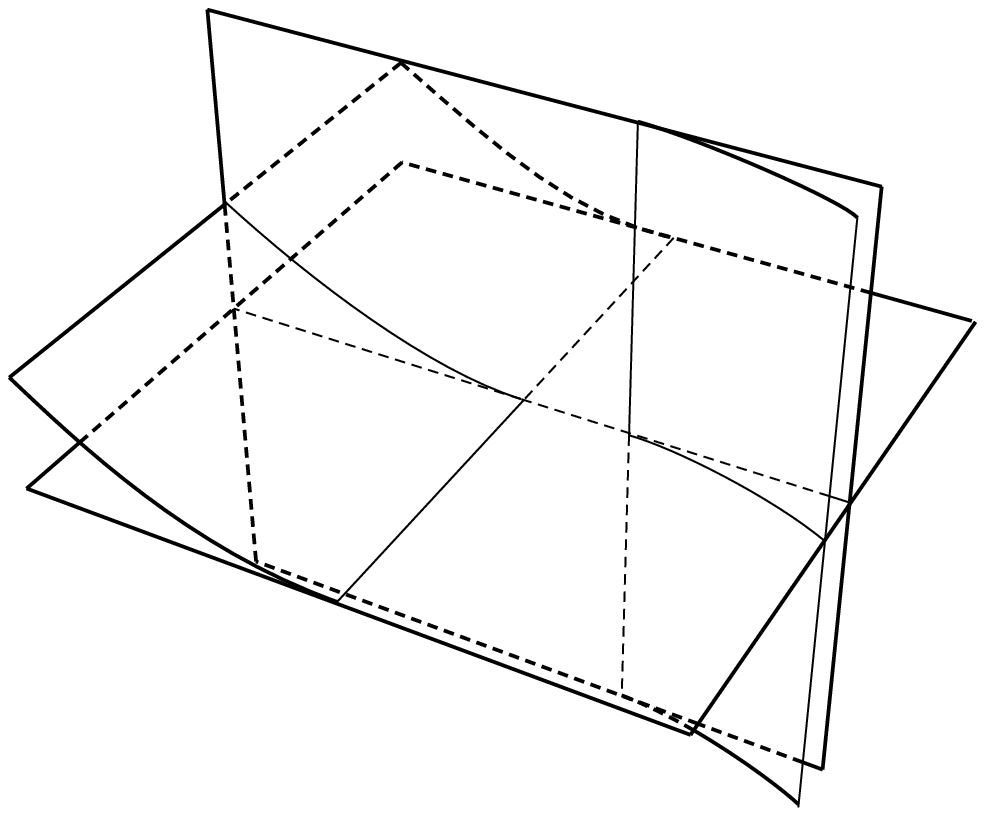}
  \end{center}
 \end{minipage}
$\leftrightarrow $
 \begin{minipage}{0.30\hsize}
   \begin{center}
     \includegraphics*[width=3cm,height=3cm]{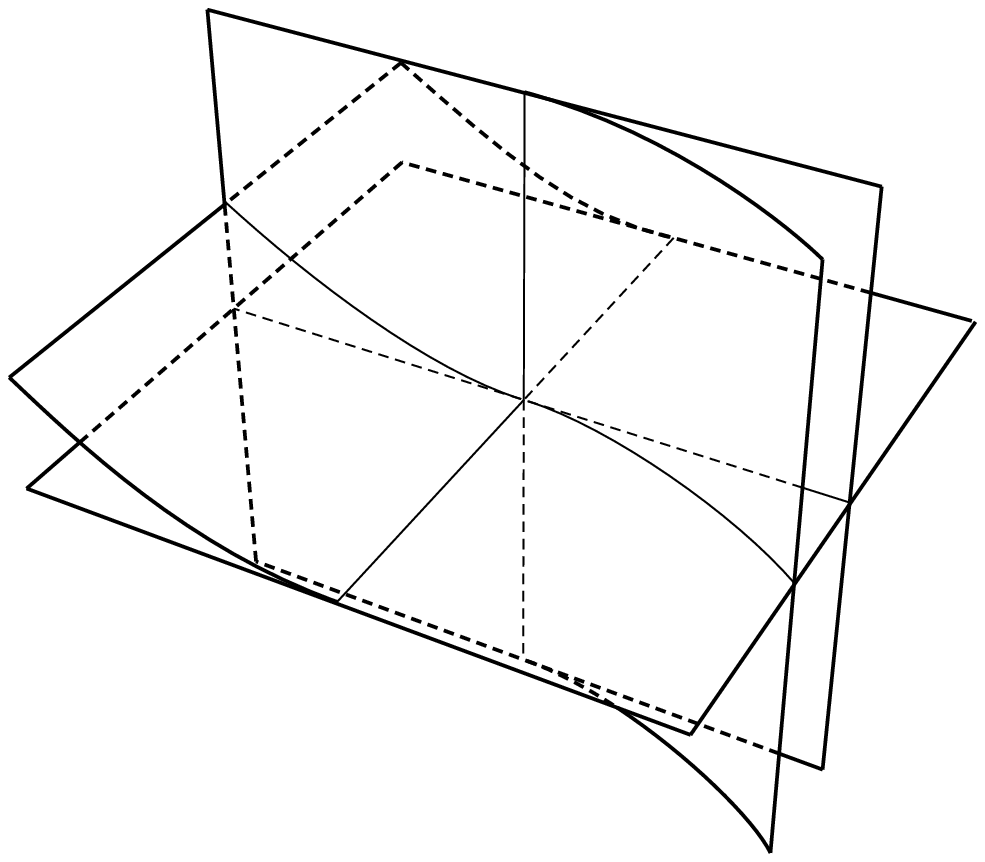}
   \end{center}
 \end{minipage}
$\leftrightarrow $
 \begin{minipage}{0.30\hsize}
  \begin{center}
  \includegraphics*[width=3cm,height=3cm]{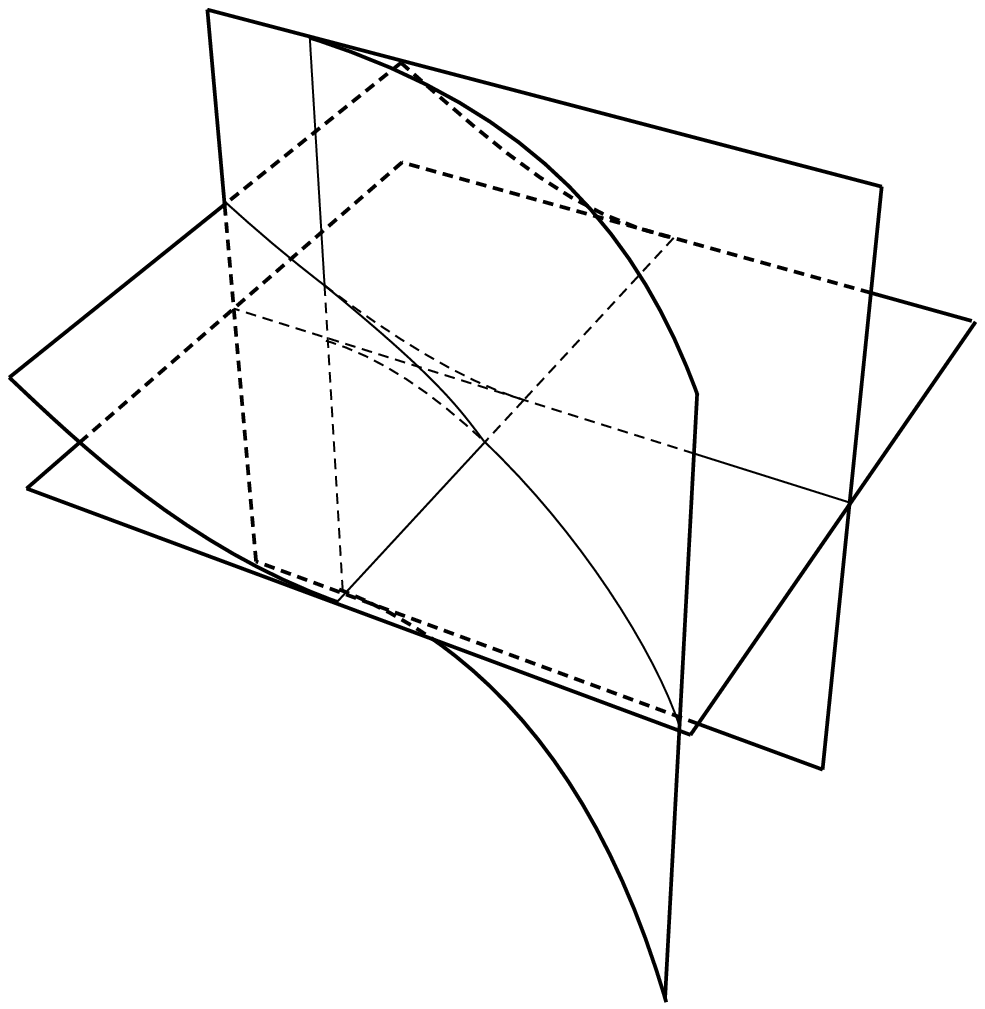}
  \end{center}
\end{minipage}
\caption{${}^1({}^0B_2{}^0B_2)$}
\end{figure}  
\begin{figure}[htbp]
 \begin{minipage}{0.30\hsize} 
  \begin{center}
    \includegraphics*[width=3cm,height=3cm]{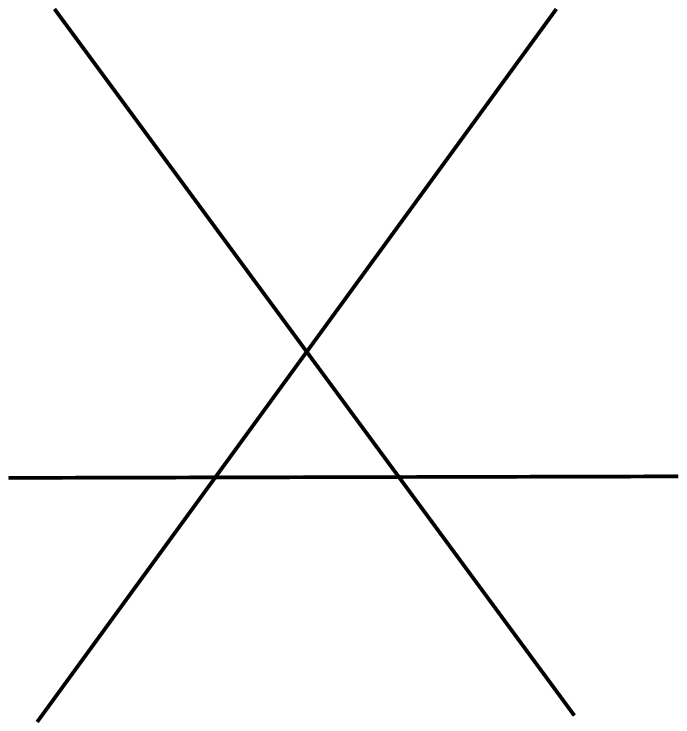}
  \end{center}
 \end{minipage}
$\leftrightarrow $
 \begin{minipage}{0.30\hsize}
   \begin{center}
     \includegraphics*[width=3cm,height=3cm]{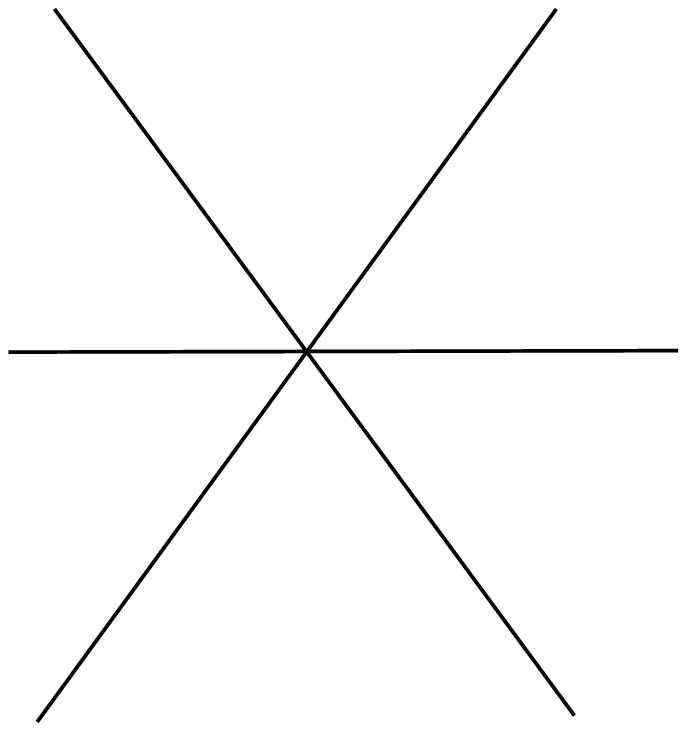}
   \end{center}
 \end{minipage}
$\leftrightarrow $
 \begin{minipage}{0.30\hsize}
  \begin{center}
  \includegraphics*[width=3cm,height=3cm]{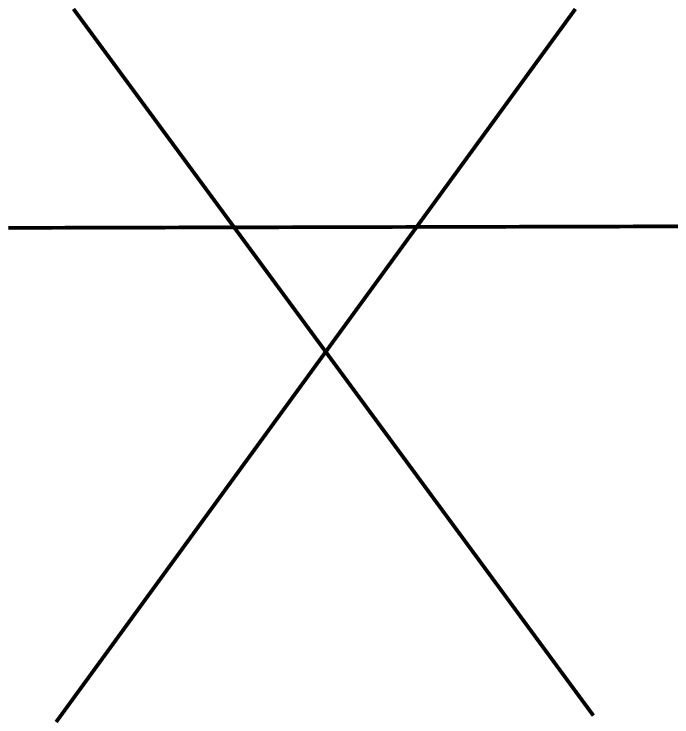}
  \end{center}
\end{minipage}\\
 \begin{minipage}{0.30\hsize} 
  \begin{center}
    \includegraphics*[width=3cm,height=3cm]{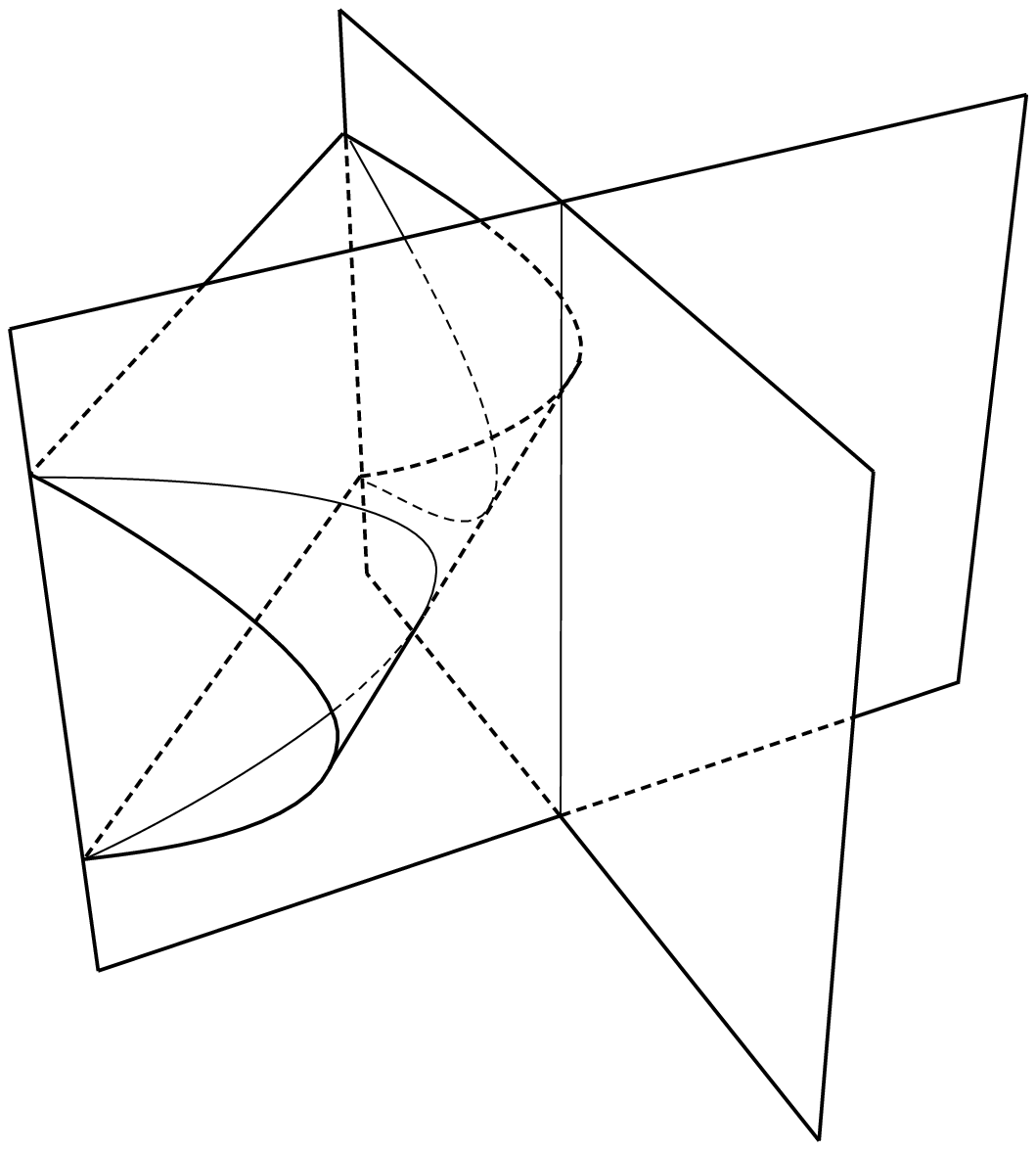}
  \end{center}
 \end{minipage}
$\leftrightarrow $
 \begin{minipage}{0.30\hsize}
   \begin{center}
     \includegraphics*[width=3cm,height=3cm]{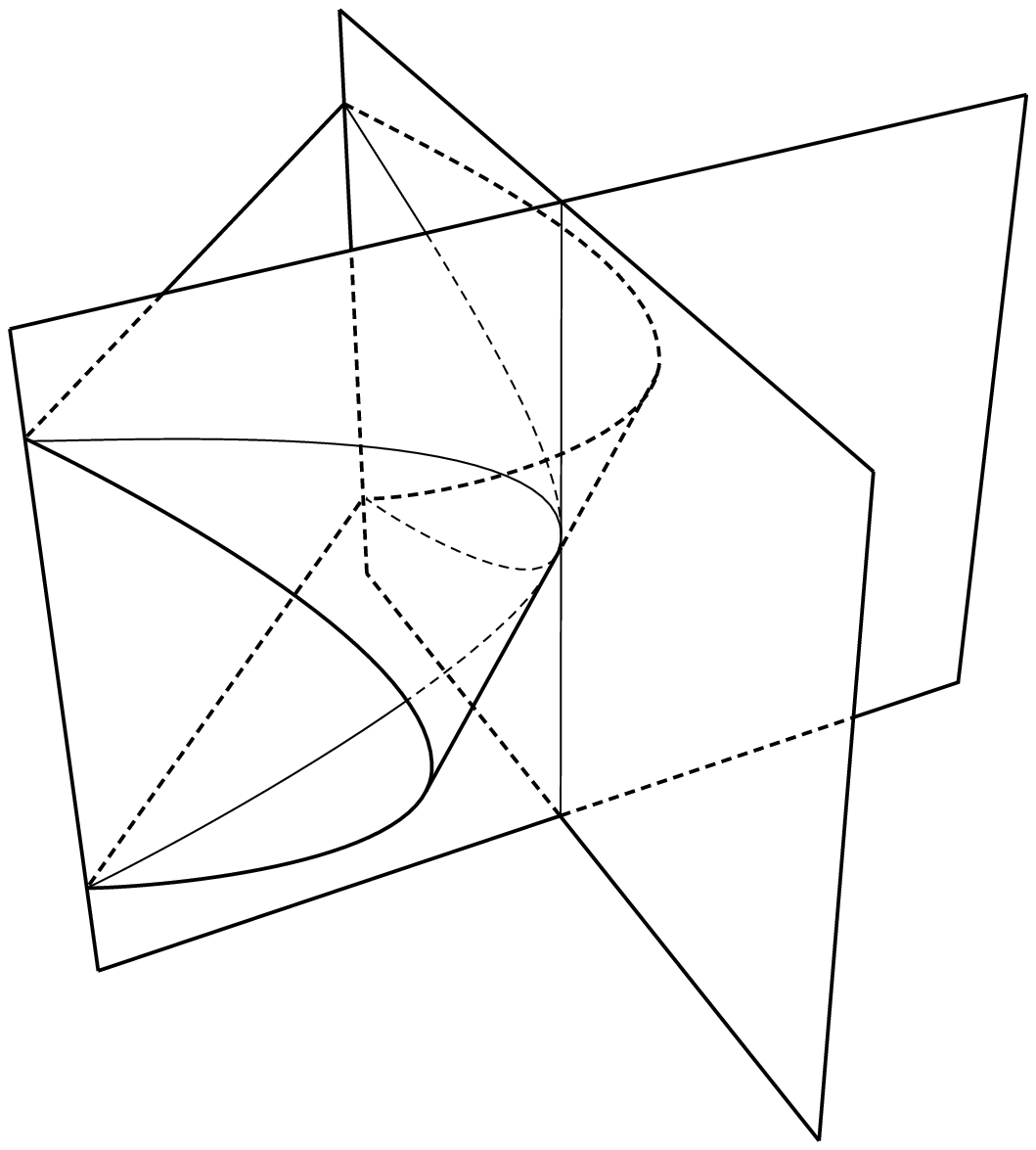}
   \end{center}
 \end{minipage}
$\leftrightarrow $
 \begin{minipage}{0.30\hsize}
  \begin{center}
  \includegraphics*[width=3cm,height=3cm]{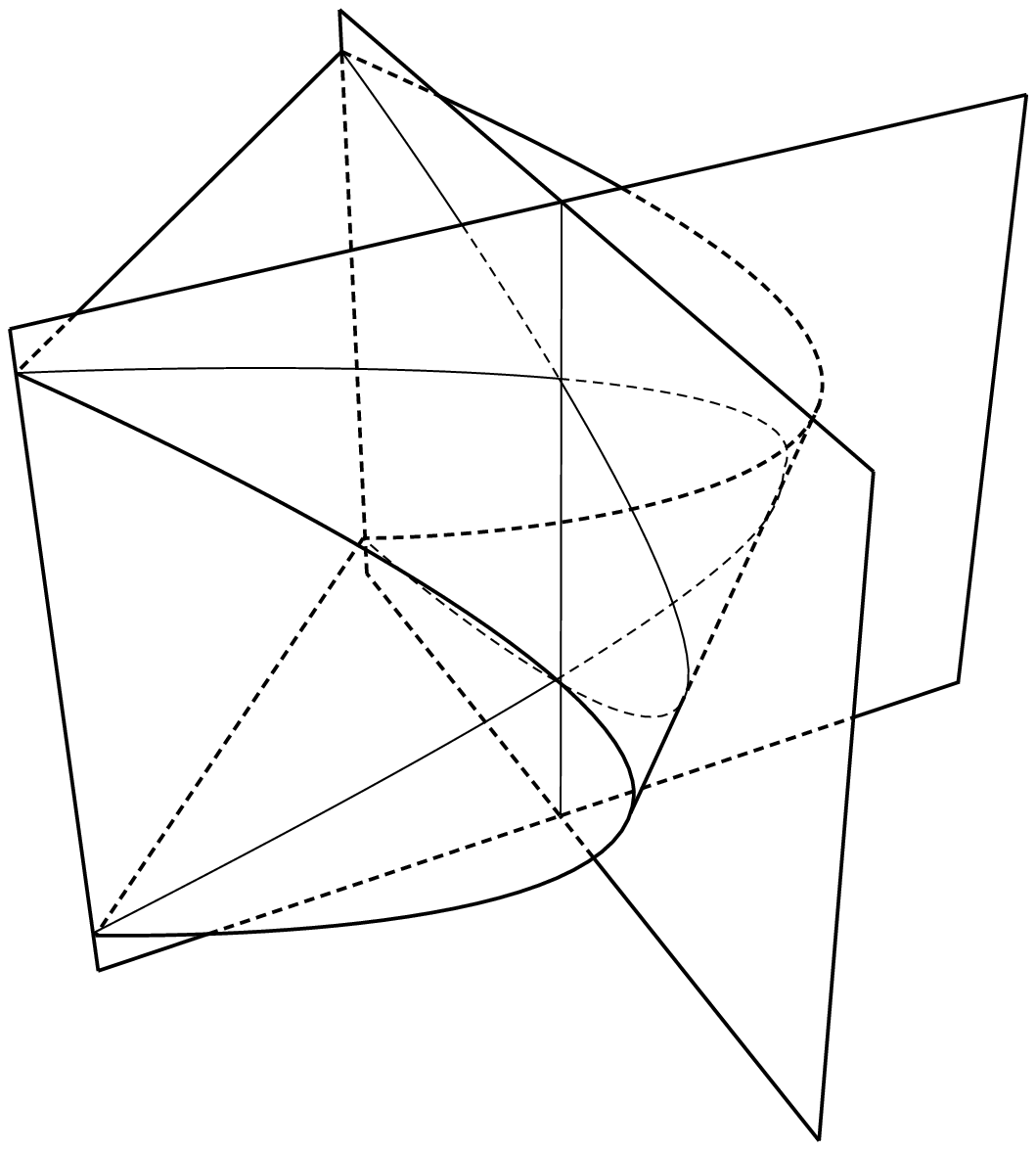}
  \end{center}
\end{minipage}
\caption{${}^1({}^0A_1{}^0A_1{}^0A_1)$}
\end{figure}  
\begin{figure}
 \begin{minipage}{0.30\hsize} 
  \begin{center}
    \includegraphics*[width=3cm,height=3cm]{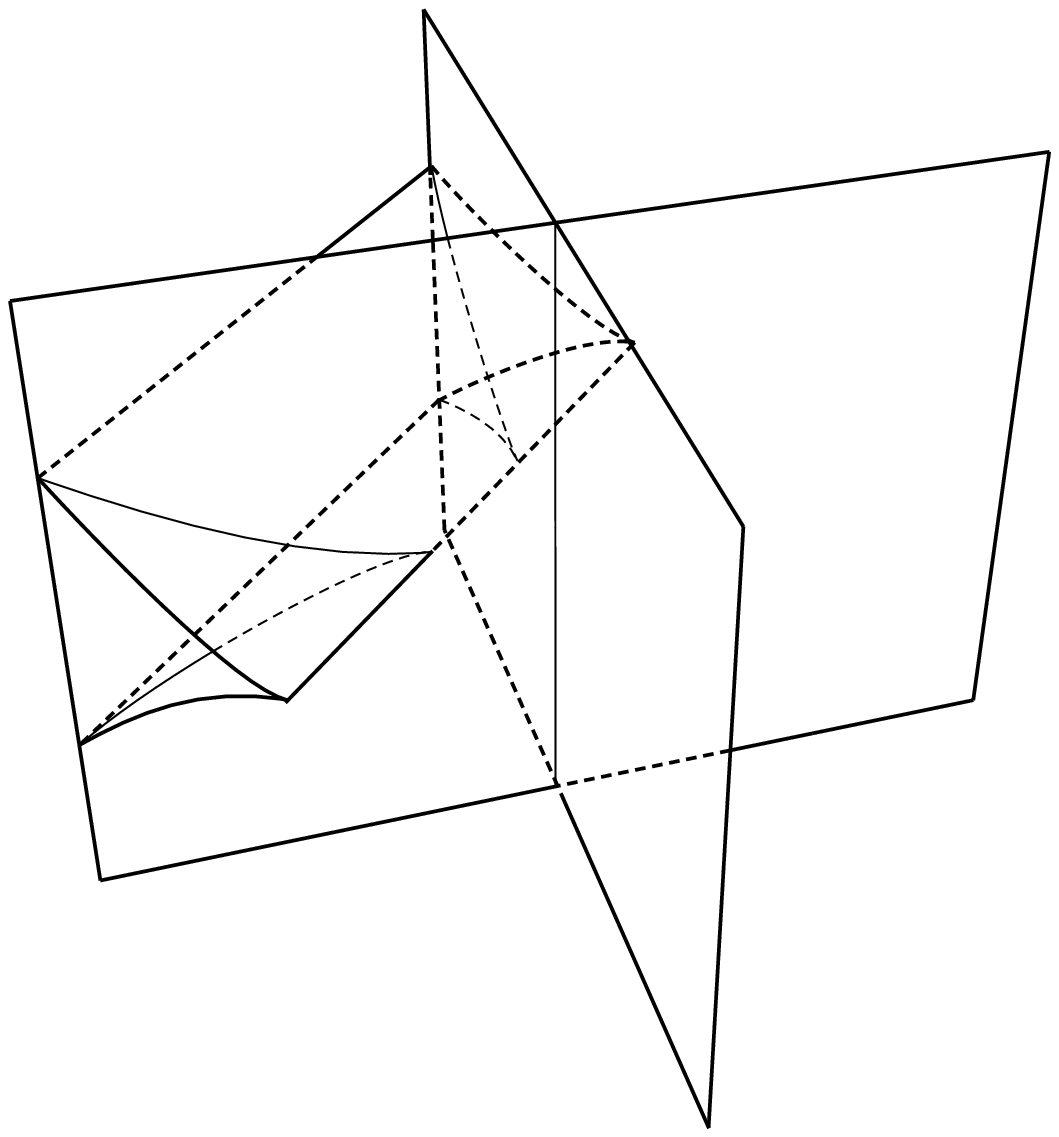}
  \end{center}
 \end{minipage}
$\leftrightarrow $
 \begin{minipage}{0.30\hsize}
   \begin{center}
     \includegraphics*[width=3cm,height=3cm]{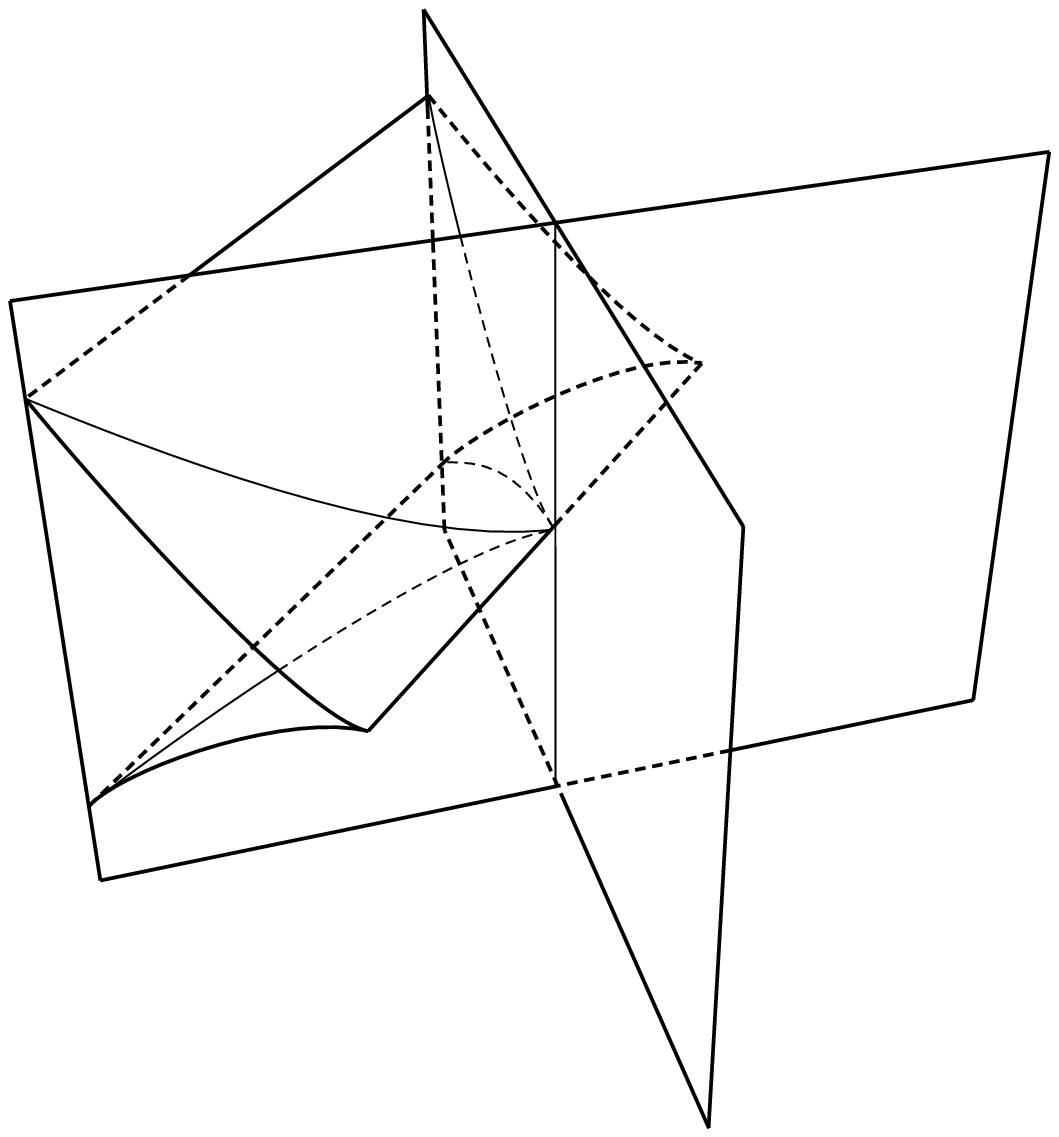}
   \end{center}
 \end{minipage}
$\leftrightarrow $
 \begin{minipage}{0.30\hsize}
  \begin{center}
  \includegraphics*[width=3cm,height=3cm]{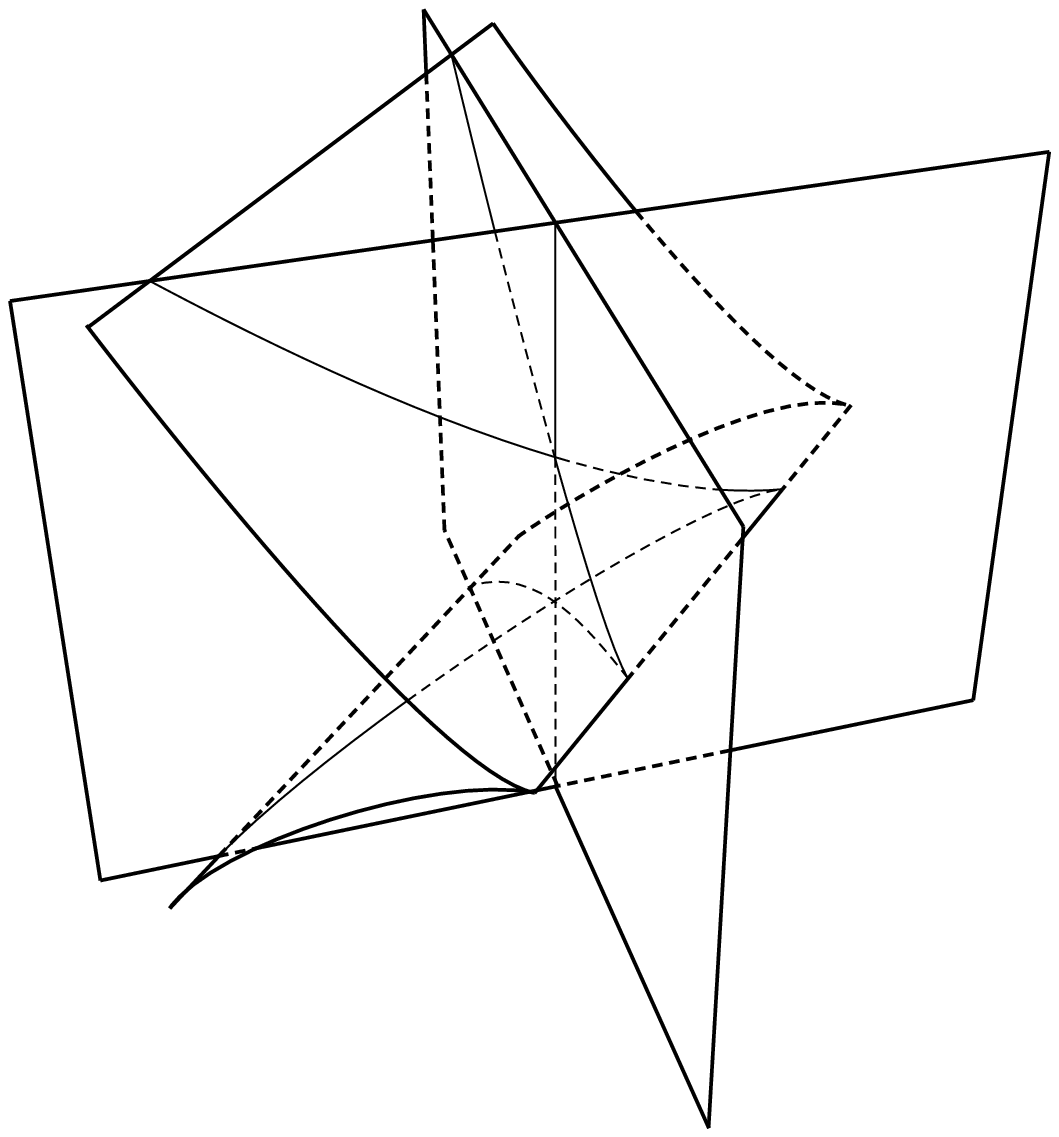}
  \end{center}
\end{minipage}
\caption{${}^1({}^0A_1{}^0A_1{}^0A_2)$}
\end{figure}  
\begin{figure}
 \begin{minipage}{0.30\hsize} 
  \begin{center}
    \includegraphics*[width=3cm,height=3cm]{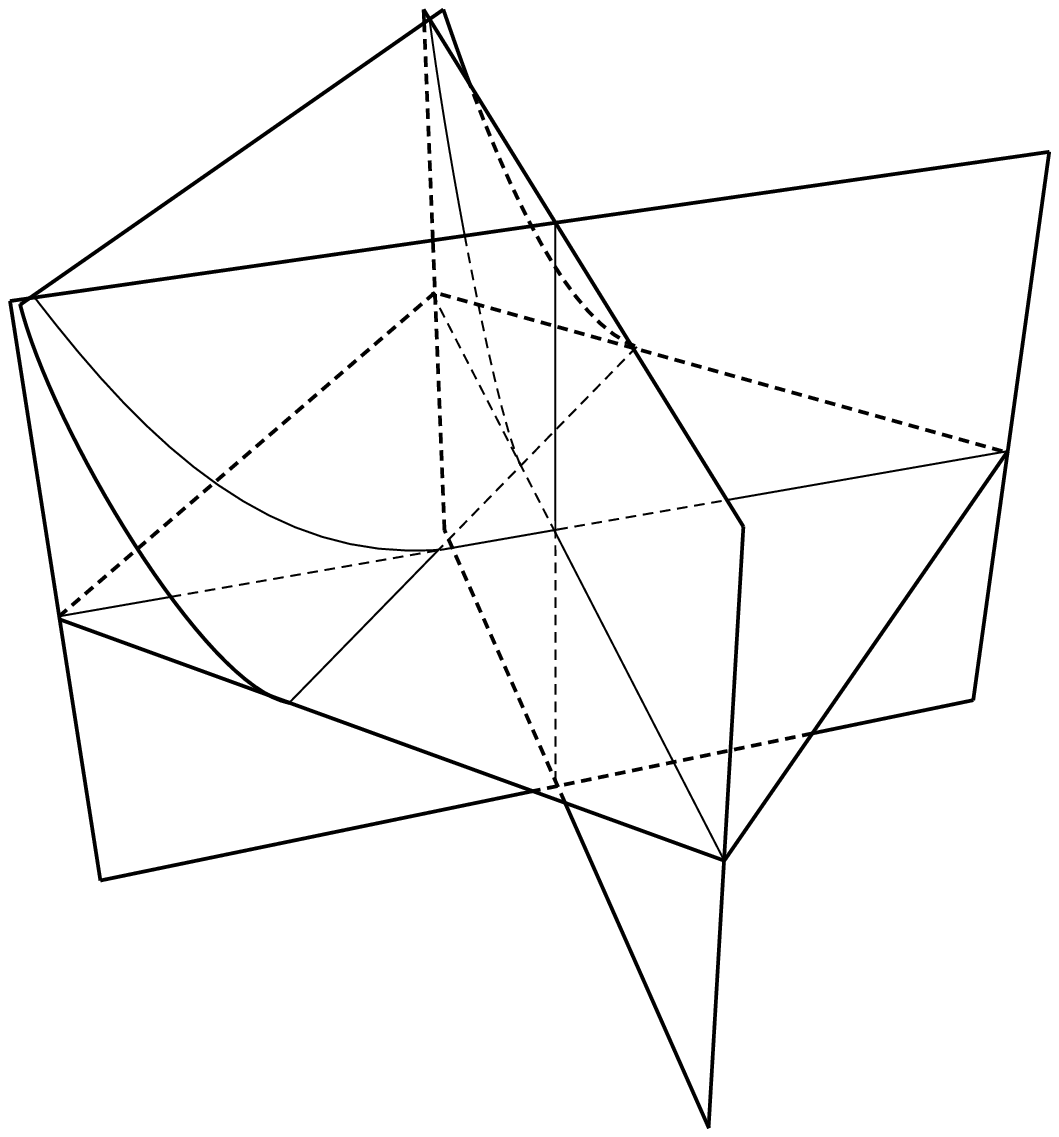}
  \end{center}
 \end{minipage}
$\leftrightarrow $
 \begin{minipage}{0.30\hsize}
   \begin{center}
     \includegraphics*[width=3cm,height=3cm]{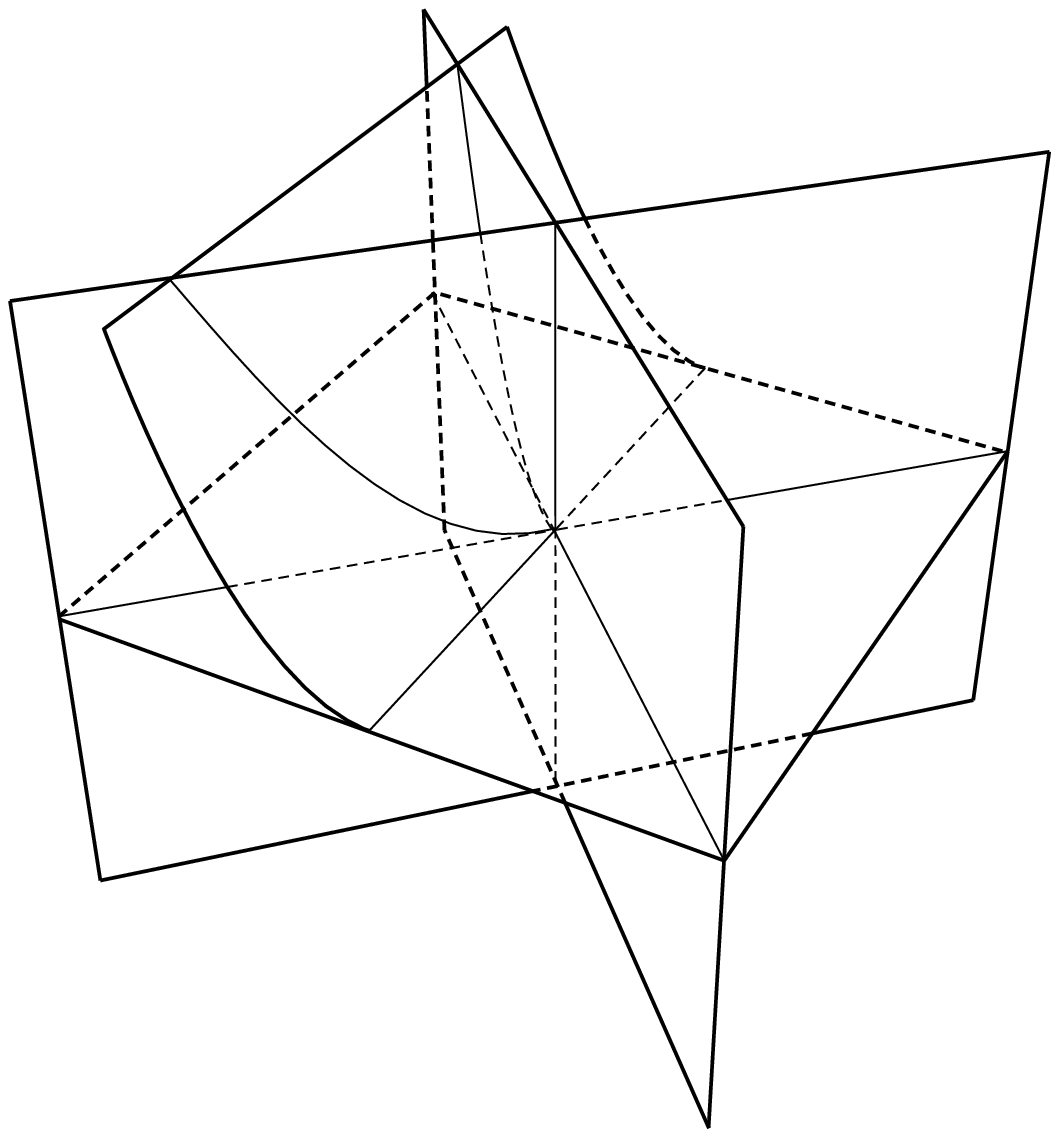}
   \end{center}
 \end{minipage}
$\leftrightarrow $
 \begin{minipage}{0.30\hsize}
  \begin{center}
  \includegraphics*[width=3cm,height=3cm]{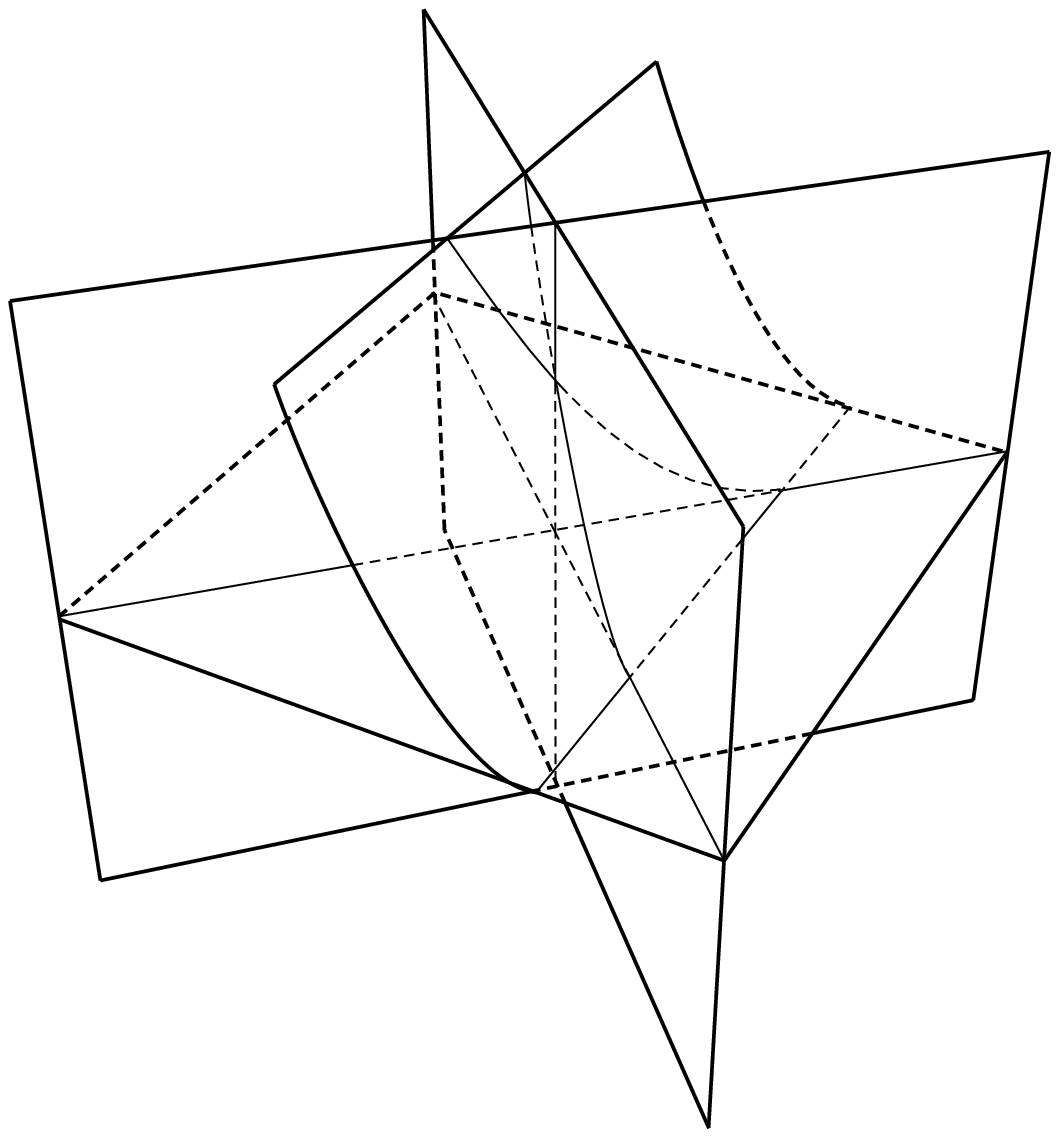}
  \end{center}
\end{minipage}
\caption{${}^1({}^0A_1{}^0A_1{}^0B_2)$}
\end{figure}  
\begin{figure}
 \begin{minipage}{0.30\hsize} 
  \begin{center}
    \includegraphics*[width=3cm,height=3cm]{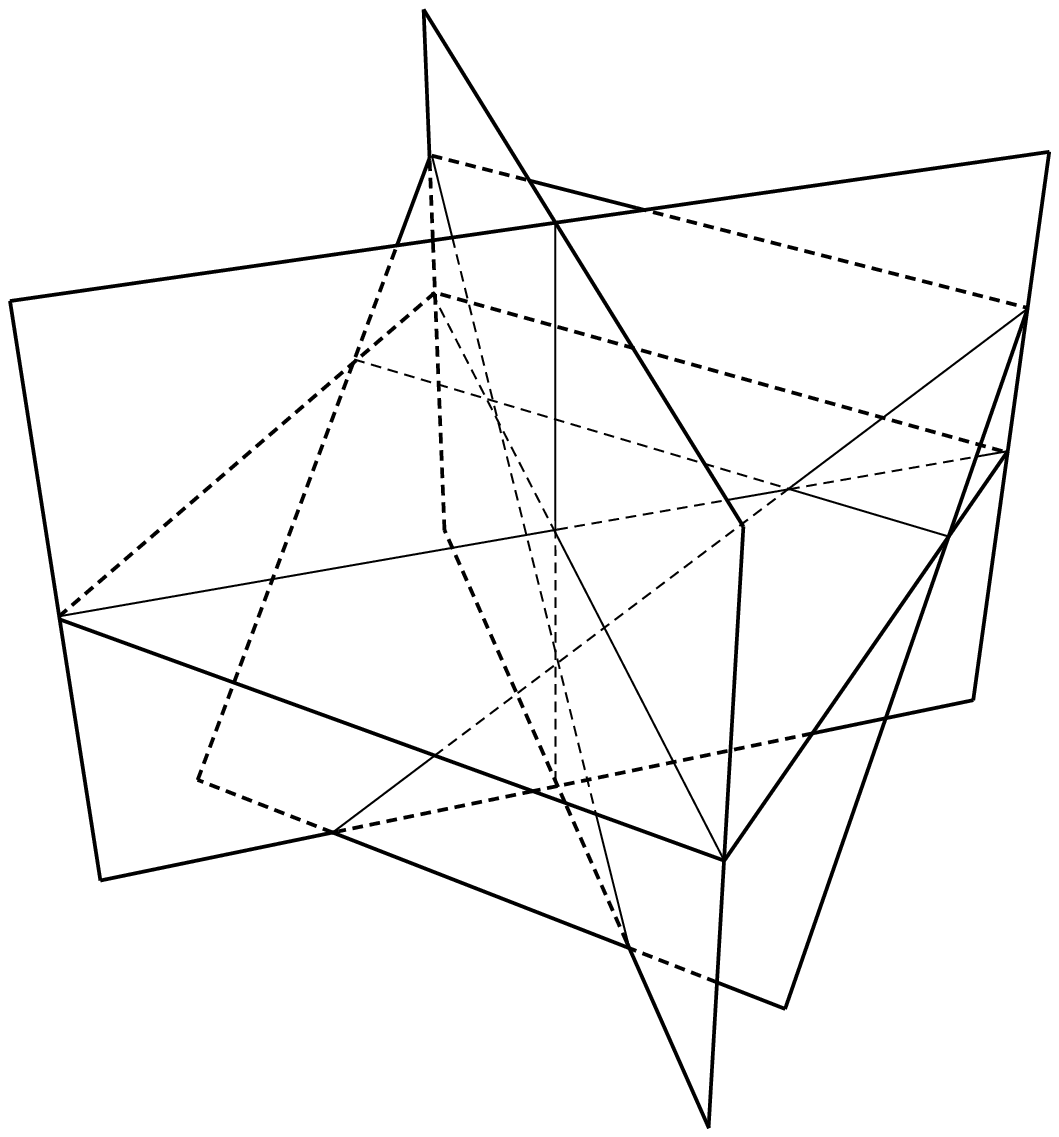}
  \end{center}
 \end{minipage}
$\leftrightarrow $
 \begin{minipage}{0.30\hsize}
   \begin{center}
     \includegraphics*[width=3cm,height=3cm]{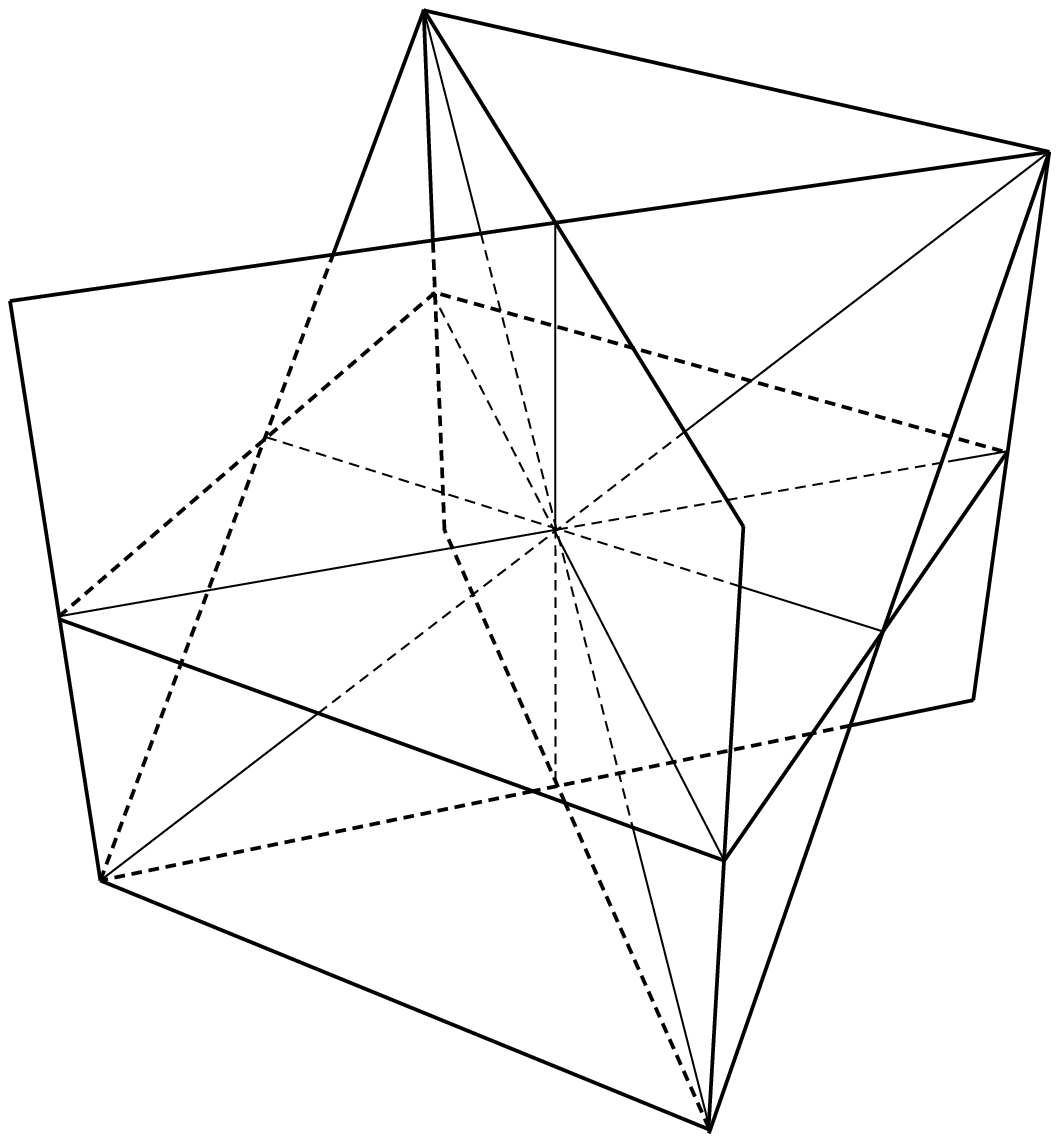}
   \end{center}
 \end{minipage}
$\leftrightarrow $
 \begin{minipage}{0.30\hsize}
  \begin{center}
  \includegraphics*[width=3cm,height=3cm]{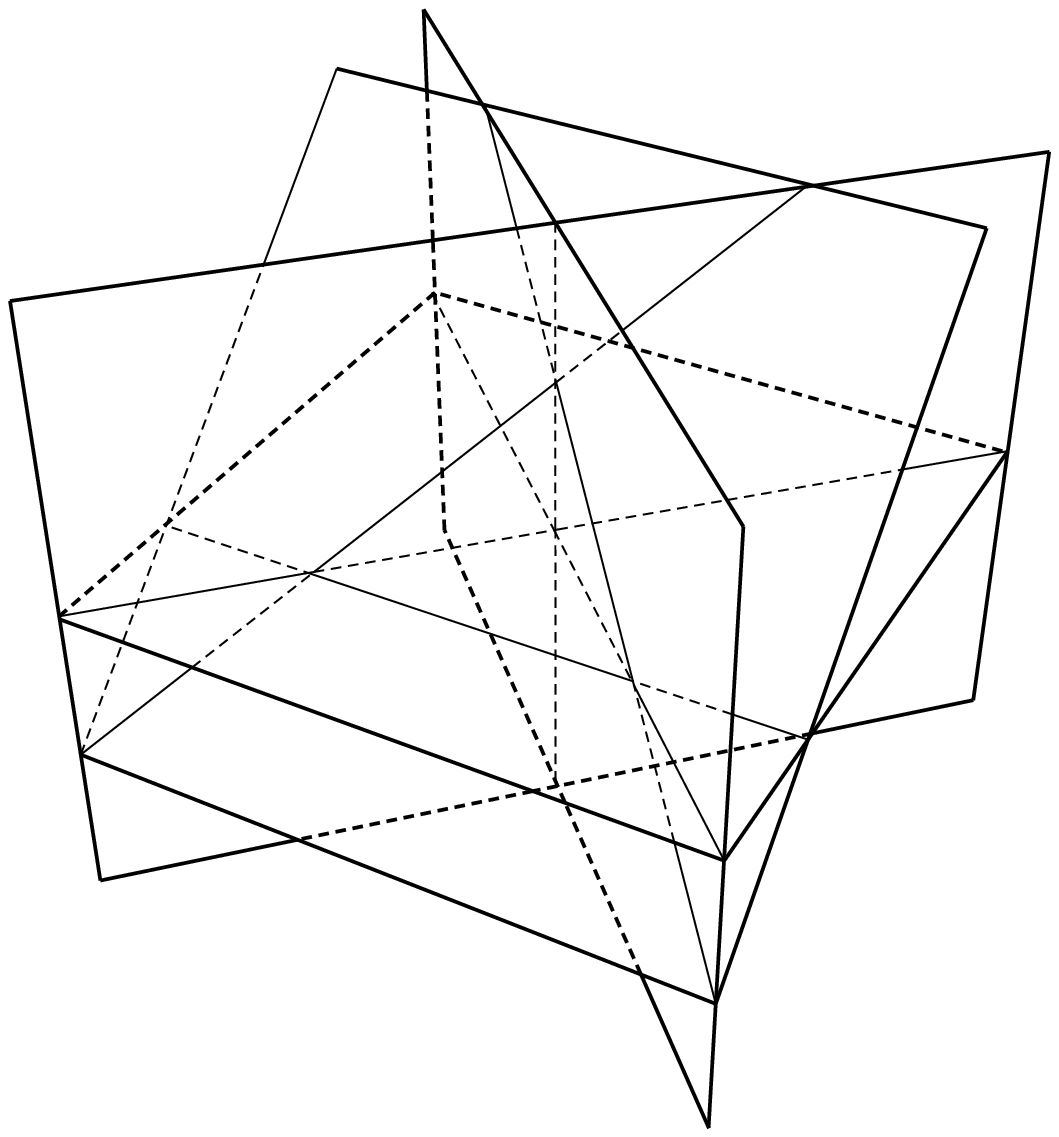}
  \end{center}
\end{minipage}
\caption{${}^1({}^0A_1{}^0A_1{}^0A_1{}^0A_1)$}
\end{figure}

\end{document}